\documentclass[11pt]{amsart}
\usepackage{amsmath}
\usepackage{amsxtra}
\usepackage{amscd}
\usepackage{amsthm}
\usepackage{amsfonts}
\usepackage{amssymb}
\usepackage{eucal}

\input xy
\xyoption{all}

\newtheorem{thm}{Theorem}
\newtheorem{conj}[thm]{Conjecture}
\newtheorem{cor}[thm]{Corollary}
\newtheorem{lem}[thm]{Lemma}
\newtheorem{prop}[thm]{Proposition}

\theoremstyle{remark}
\newtheorem{remark}{Remark}

\theoremstyle{definition}
\newtheorem{definition}[thm]{Definition}

\numberwithin{thm}{section}
\numberwithin{equation}{section}

\newcommand{\secref}[1]{Section~\ref{#1}}
\newcommand{\lemref}[1]{Lemma~\ref{#1}}

\newcommand{\conjref}[1]{Conjecture~\ref{#1}}
\newcommand{\nc}{\newcommand}
\nc{\on}{\operatorname}
\nc{\ch}{\mbox{ch}}
\nc{\Z}{{\mathbb Z}}
\nc{\C}{{\mathbb C}}
\nc{\pone}{{\mathbb P}^1}
\nc{\pa}{\partial}
\nc{\F}{{\mathcal F}}
\nc{\arr}{\rightarrow}
\nc{\larr}{\longrightarrow}
\nc{\al}{\alpha}
\nc{\ri}{\rangle}
\nc{\lef}{\langle}
\nc{\W}{{\mathbb W}}
\nc{\la}{\lambda}
\nc{\ep}{\epsilon}
\nc{\su}{\widehat{{\mathfrak s}{\mathfrak l}}_2}
\nc{\sw}{{\mathfrak s}{\mathfrak l}}
\nc{\g}{{\mathfrak g}}
\nc{\h}{{\mathfrak h}}
\nc{\n}{{\mathfrak n}}
\nc{\N}{\widehat{\n}}
\nc{\G}{\widehat{\g}}
\nc{\De}{\Delta}
\nc{\gt}{\widetilde{\g}}
\nc{\Ga}{\Gamma}
\nc{\one}{{\mathbf 1}}
\nc{\z}{{\mathfrak Z}}
\nc{\La}{\Lambda}
\nc{\wt}{\widetilde}
\nc{\wh}{\widehat}
\nc{\cri}{_{\kappa_c}}
\nc{\kk}{\underline{\mathbf C}}
\nc{\sun}{\widehat{\sw}_N}
\nc{\si}{\sigma}
\nc{\el}{\ell}
\nc{\bi}{\bibitem}
\nc{\om}{\omega}
\nc{\ol}{\overline}
\nc{\ds}{\displaystyle}
\nc{\dzz}{\frac{dz}{z}}
\nc{\Res}{\on{Res}}
\nc{\mc}{\mathcal}
\nc{\Cal}{\mathcal}
\nc{\bb}{{\mathfrak b}}
\nc{\ot}{\otimes}
\nc{\R}{{\mc R}}
\nc{\yy}{{\mc Y}}
\nc{\ga}{\gamma}

\nc{\us}{\underset}
\nc{\opl}{\oplus}

\nc{\Fq}{{\mathbb F}_q}
\nc{\Mq}{{\mathcal M}}
\nc{\Rep}{\on{Rep}}
\nc{\sssec}{\subsubsection}
\nc{\ssec}{\subsection}
\nc{\lan}{\langle}
\nc{\ran}{\rangle}

\nc{\D}{\mathcal D}
\nc{\Vect}{\on{Vect}}
\nc{\ghat}{\G}
\nc{\T}{\mc T}
\nc{\Tloc}{\T^\g_{\on{loc}}}
\nc{\vac}{|0\ran}
\nc{\Wick}{{\mb :}}
\nc{\mb}{\mathbf}
\nc{\delz}{\partial_z}
\nc{\K}{{\mathbb K}}
\nc{\cali}{\mathcal}
\nc{\li}{\mathfrak l}
\nc{\lt}{\widetilde{\li}}
\nc{\astar}{a^*}
\nc{\cA}{{\mc A}}
\nc{\ka}{\kappa}

\nc{\OO}{{\mc O}}
\nc{\AutO}{\on{Aut}\OO}
\nc{\DerO}{\on{Der}\OO}
\nc{\DerpO}{\on{Der}_+\OO}
\nc{\Au}{{\mc A}ut}
\nc{\mf}{\mathfrak}
\nc{\V}{{\mathcal V}}
\nc{\hh}{\wh{\h}}

\nc{\pp}{{\mathfrak p}}
\nc{\mm}{{\mathfrak m}}
\nc{\rr}{{\mathfrak r}}
\nc{\ket}{\rangle}
\nc{\zz}{{\mathfrak z}}
\nc{\gr}{\on{gr}}
\nc{\Spe}{\on{Spec}}
\nc{\rv}{\crho}
\nc{\can}{\on{can}}
\nc{\CC}{\on{Op}_G(D))}
\nc{\Op}{\on{Op}_G(D)}
\nc{\MOp}{\on{MOp}_G(D)}
\nc{\Db}{{\mathbb D}}
\nc{\ww}{w}

\nc{\af}{{\mathbb A}^1}
\nc{\bs}{\backslash}
\nc{\laa}{(\la_i)}
\nc{\zn}{(z_i)}

\nc{\cla}{\check{\la}}
\nc{\cmu}{\check{\mu}}
\nc{\crho}{\check{\rho}}
\nc{\chal}{\check{\al}}
\nc{\cc}{{\mathfrak c}}

\nc{\MM}{{\mathbb M}}

\nc{\ZZ}{{\mc Z}}

\nc{\UU}{{\mathbb U}}

\nc{\Conn}{\on{Conn}(\Omega^{\crho})}
\nc{\Con}{\on{Conn}(\Omega^{-\rho})}
\nc{\Co}{\on{Conn}(\Omega^{\rho})}

\nc{\ppart}{(\!(t)\!)}
\nc{\pparl}{(\!(\la)\!)}
\nc{\zpart}{(\!(z)\!)}
\nc{\ppzi}{(\!(t-z_i)\!)}
\nc{\ppinf}{(\!(t^{-1})\!)}
\nc{\Ind}{\on{Ind}}
\nc{\I}{{\mathbb I}}

\nc{\Bun}{\on{Bun}}

\newcommand {\IC}{\on{IC}}

\newcommand {\tr}{\operatorname{tr}}

\nc{\gtil}{\wt{\g}}
\nc{\ntil}{\wt{\n}}
\nc{\htil}{\wt{\h}}
\nc{\gbar}{\ol{\g}}
\nc{\nbar}{\ol{\n}}
\nc{\bbar}{\ol{\bb}}
\nc{\lhat}{\wh{\mf l}}

\nc{\ovc}{\overset{\circ}}
\nc{\Gr}{\on{Gr}}

\nc{\AD}{{\mathbb A}}
\nc{\Gm}{{\mathbb G}_m}
\nc{\Ql}{{\mathbb Q}_\ell}

\nc{\Loc}{\on{Loc}}

\nc\GG{\mathbb G}
\nc\Xb{\mathbf X}
\nc\inv{{\rm inv}}
\nc\isom{=}
\nc\Hc{{\mathcal H}}
\nc\ovl{\overline}
\newcommand{\Pic}{\on{Pic}}
\newcommand{\Jac}{\on{Jac}}
\newcommand{\BD}{{}}
\newcommand\Oc{\mathcal{O}}
\nc\Lc{{\mathcal L}}
\nc\Ec{{\mathcal E}}
\nc\wF{{\mc K}}
\nc\pr{{\rm pr}}
\nc\SL{{\rm SL}}
\nc\PGL{{\rm PGL}}
\nc\GL{{\rm GL}}
\nc\LG{{}^L G}
\nc\Pc{{\mc P}}
\nc\Ac{{\mc A}}
\nc\dv{{\rm div}}
\nc\Fc{{\mc F}}
\nc\M{{\mc M}}

\nc{\und}{\underline}

\def\H{{\mathbb H}}

\begin{document}

\title{Geometrization of trace formulas}

\author[Edward Frenkel]{Edward Frenkel$^1$}\thanks{$^1$Supported by
  DARPA under the grant HR0011-09-1-0015}

\address{Department of Mathematics, University of California,
Berkeley, CA 94720, USA}

\author[Ng\^o Bao Ch\^au]{Ng\^o Bao Ch\^au$^2$}\thanks{$^2$Supported
  by the Simonyi Foundation and NSF under the grant DMS-1000356}

\address{School of Mathematics, Institute for Advanced Study,
  Princeton, NJ 08540, USA}

\date{April 2010. Revised March 2011.}

\maketitle

\section{Introduction}

\setcounter{footnote}{1}

Recently, a new approach to functoriality of automorphic
representations \cite{L} has been proposed in \cite{FLN} following
earlier work \cite{L:BE,L:PR,L:Shaw} by Robert Langlands. The idea may
be roughly summarized as follows. Let $G$ and $H$ be two reductive
algebraic groups over a global field $F$ (which is either a number
field or a function field; that is, the field of functions on a smooth
projective curve $X$ over a finite field), and assume that $G$ is
quasi-split. Let $^L G$ and $^L H$ be the Langlands dual groups as
defined in \cite{L}. The {\em functoriality} principle states that for
each homomorphism $^L H \to {}^L G$ there exists a transfer of
automorphic representations from $H(\AD)$ to $G(\AD)$, where $\AD$ is
the ring of ad\'eles of $F$. In other words, to each $L$-packet of
automorphic representations of $H(\AD)$ should correspond an
$L$-packet of automorphic representations of $G(\AD)$, and this
correspondence should satisfy some natural properties. Functoriality
has been established in some cases, but is unknown in general.

In \cite{L:BE,L:PR,L:Shaw,FLN} the following strategy for proving
functoriality was proposed. In the space of automorphic functions on
$G(F)\bs G(\AD)$ one should construct a family of integral operators
which project onto those representations which come by functoriality
from other groups $H$. One should then use the trace formula to
decompose the traces of these operators as sums over these $H$ (there
should be sufficiently many operators to enable us to separate the
different $H$). The hope is that analyzing the orbital side of the
trace formula and comparing the corresponding orbital integrals for
$G$ and $H$, one should ultimately be able prove functoriality.

In the present paper we make the first steps in developing geometric
methods for analyzing these orbital integrals in the case of the
function field of a curve $X$ over a finite field $\Fq$. We also
suggest a conjectural framework of ``geometric trace formulas'' in the
case of curves defined over the complex field.

\medskip

In \cite{FLN} a family of averaging operators ${\mb K}_{d,\rho}$ has
been constructed. They depend on a positive integer $d$ and an
irreducible representation $\rho$ of $^L G$. The operator ${\mb
  K}_{d,\rho}$ is expected to have the following property: for
sufficiently large $d$, the representations of $G(\AD)$ that
contribute to the trace of ${\mb K}_{d,\rho}$ are those coming by
functoriality from the groups $H$ such that the dual group $^L H$ has
non-zero invariant vectors in $\rho$. To prove this, we wish to use
the {\em trace formula}
\begin{equation}    \label{tr}
\on{Tr} {\mb K}_{d,\rho} = \int K_{d,\rho}(x,x) \; dx,
\end{equation}
where $K_{d,\rho}(x,y)$ is the integral kernel of ${\mb K}_{d,\rho}$
on $G(F)\bs G(\AD) \times G(F)\bs G(\AD)$.\footnote{More precisely, we
  should consider the stabilized version of the trace formula, which
  is more efficient for comparing different groups.} Thus, we need to
express the right hand side of this formula as the sum of orbital
integrals and to compare these orbital integrals for $G$ and the
groups $H$ described above.\footnote{\label{foot1}We note that in
  general we should insert additional operators at the finite set $S$
  of closed points in $X$. However, in this paper, in order to
  simplify the initial discussion, we will restrict ourselves to the
  case when $S$ is either empty or the operators inserted belong to
  the spherical Hecke algebra (with respect to a particular choice of
  a maximal compact subgroup of $G(\AD)$). Thus, only the unramified
  representations of $G(\AD)$ (with respect to this subgroup) will
  contribute to \eqref{tr}.}

\medskip

Here an important remark is in order. There are two types of
automorphic representations of $G(\AD)$: tempered (or Ramanujan) and
non-tempered (or non-Ramanujan), see \secref{h eig} for more
details. As explained in \cite{FLN}, for the above strategy to work,
one first needs to remove the contribution of the non-tempered
representations from the trace of ${\mathbf K}_{d,\rho}$. This was the
focus of \cite{FLN}, where it was shown how to separate the
contribution of the trivial representation of $G(\AD)$ (which is in
some sense the most non-tempered).

\medskip

In the first part of this paper we show that the kernel $K_{d,\rho}$
may be obtained using the Grothendieck's {\em faisceaux--fonctions}
dictionary from a perverse sheaf ${\mc K}_{d,\rho}$ on a certain
algebraic stack over the square of $\Bun_G$, the moduli stack of
$G$-bundles on $X$. Hence the right hand side of the trace formula
\eqref{tr} may be written as the trace of the Frobenius (of the Galois
group of the finite field $k$, over which our curve $X$ is defined) on
the \'etale cohomology of the restriction of ${\mc K}_{d,\rho}$ to the
diagonal in $\Bun_G \times \Bun_G$. We show that the latter is closely
related to the Hitchin moduli stacks of Higgs bundles on $X$. The idea
is then to use the geometry of these moduli stacks to prove the
desired identities of orbital integrals for $G$ and $H$ by
establishing isomorphisms between the corresponding cohomologies. In
\secref{orbital} we formulate a precise conjecture relating the
cohomologies in the case that $G={\rm SL}_2$ and $H$ is a twisted
torus.

A prototype for this is the proof of the fundamental lemma (in the
setting of Lie algebras) given by one of the authors
\cite{Ngo:FL}. However, our situation is different. The argument of
\cite{Ngo:FL} used decomposition of the cohomology of the fibers of
the Hitchin map under the action of finite groups, whereas in our case
the decomposition of the cohomology we are looking for does not seem
to be due to an action of a group.

\medskip

One advantage of the geometric approach is that the moduli stacks and
the sheaves on them that appear in this picture have natural analogues
when the curve $X$ is defined over $k=\C$, rather than a finite
field. Some of our questions may be reformulated over $\C$ as well,
and we can use methods of complex algebraic geometry (some of which
have no obvious analogues over a finite field) to tackle them.

\medskip

This leads to a natural question: what is the analogue of the trace
formula \eqref{tr} over $\C$? We discuss this question in the second,
more speculative, part of the paper and formulate two conjectures in
this direction.

\medskip

The geometric analogue of the right hand (orbital) side of \eqref{tr}
is the cohomology of a sheaf on a moduli stack, as we explain in the
first part. We would like to find a similar interpretation of the left
hand (spectral) side of \eqref{tr}. Again, to simplify our discussion,
consider the unramified case, but with a more general kernel $K$
instead of $K_{d,\rho}$ which is bi-invariant with respect to a fixed
maximal compact subgroup $G(\OO)$ of $G(\AD)$. Then the left hand side
of \eqref{tr} may be rewritten as follows (we ignore the continuous
spectrum for the moment):
\begin{equation}    \label{s1}
\sum_\pi m_\pi N_{\pi},
\end{equation}
where the sum is over irreducible representations $\pi$ of $G(\AD)$,
unramified with respect to $G(\OO)$, $m_\pi$ is the multiplicity of
$\pi$ in the space of automorphic functions, and $N_{\pi}$ is the
eigenvalue of $K_{d,\rho}$ on a spherical vector in $\pi$ fixed by
this subgroup. We wish to interpret this sum as the Lefschetz trace
formula for the trace of the Frobenius on the \'etale cohomology of an
$\ell$-adic sheaf. It is not obvious how to do this, because the set
of the $\pi$'s appearing in \eqref{s1} is not the set of points of a
moduli space (or stack) in any obvious way.

\medskip

However, recall that according to the Langlands correspondence
\cite{L}, the $L$-packets of irreducible
(unramified) tempered automorphic representations of $G(\AD)$ are
supposed to be parametrized by (unramified) homomorphisms
$$
\sigma: W_F \to {}^L G,
$$
where $W_F$ is the Weil group of the function field $F$ and $^L G$ is
the Langlands dual group to $G$. This has been proved by V. Drinfeld
for $G=GL_2$ \cite{Dr1,Dr2} and L. Lafforgue for $G=GL_n, n>2$
\cite{LLaf}, but is unknown in general. Nevertheless, assume that this
is true. Then the eigenvalues $N_{\pi}$ are determined by the Hecke
eigenvalues of $\pi$, and in particular they are the same for all
$\pi$ in the $L$-packet $L_\sigma$ corresponding to any given $\sigma:
W_F \to {}^L G$. We will denote them by $N_{\sigma}$. Further, we
expect that if $\sigma$ is unramified, then there is a unique
representation $\pi$ in the $L$-packet $L_\sigma$ which is unramified
with respect to the subgroup $G(\OO)$.  Hence we will write $m_\sigma$
for $m_\pi$.  Then, assuming the Langlands correspondence, we obtain
that \eqref{s1} is equal to
\begin{equation}    \label{s2}
\sum_\sigma m_\sigma N_{\sigma}.
\end{equation}

Ideally, we would like to describe \eqref{s2} as the trace of the
Frobenius on the \'etale cohomology of a sheaf on a moduli stack,
whose set of $k$-points is the set of equivalence classes of
homomorphisms $\sigma$. Such a stack does not exist if the curve $X$
is defined over a finite field $k$. But when $X$ is defined over $\C$,
we have an algebraic stack $\Loc_{^L G}$ of de Rham local systems on
$X$. Hence we can pose the following question: define a sheaf on this
stack such that its cohomology (representing the left hand side of
\eqref{tr} in the complex case) is isomorphic to the cohomology
representing the right hand side of \eqref{tr}. It would then be
natural to call this isomorphism a {\em geometrization of the trace
  formula} \eqref{tr}.

\bigskip

Our main observation is that the answer may be obtained in the
framework of a {\em categorical form of the geometric Langlands
  correspondence}, which is a conjectural equivalence between derived
categories of coherent sheaves on the moduli stack $\Loc_{^L G}$ of
$^L G$-local systems on a complex curve $X$ and ${\mc D}$-modules on
the moduli stack $\Bun_G$ of $G$-bundles on $X$. Such an equivalence
has been proved in the abelian case by L. Laumon \cite{Laumon} and
M. Rothstein \cite{Roth}, and in the non-abelian case it has been
suggested as a conjectural guiding principle by A. Beilinson and
V. Drinfeld (see, e.g., \cite{F:rev,Laf,LL} for an
exposition).\footnote{This categorical version of the geometric
  Langlands correspondence also appears naturally in the $S$-duality
  picture developed in \cite{KW} (see \cite{F:bourbaki} for an
  exposition).}

So far, this correspondence has been mostly studied at the level of
objects. For example, the skyscraper sheaf supported at a given $^L
G$-local system $E$ on $X$ should go to the Hecke eigensheaf on
$\Bun_G$ with "eigenvalue" $E$. But if we have an equivalence of
categories, then we also obtain non-trivial information about
morphisms; namely, the Hom's between the objects corresponding to each
other on the two sides should be isomorphic. The main point of the
second part of this paper is that for suitable objects the isomorphism
of Hom's yields the sought-after geometric trace formula. (More
precisely, these objects are sheaves on the squares of $\Bun_G$ and
$\Loc_{^L G}$, which may be interpreted as ``kernels'' of functors
acting on the above two categories.)

We propose a conjectural geometrization of the trace formula
\eqref{tr} in this framework. This is still a tentative answer,
because several important issues remain unresolved, as we explain in
\secref{discussion}. Nevertheless, we believe that it contains
interesting features and hence even in this rough form it might
provide a useful framework for a better geometric understanding of the
trace formula as well as the geometric Langlands correspondence.

\medskip

In deriving this formula, we apply the categorical form of the
geometric Langlands correspondence to the sheaf ${\mc K}$, which is
the geometric incarnation of the kernel of the integral operator ${\mb
  K}$. This sheaf is defined on an algebraic stack over the square of
$\Bun_G$. It is natural to ask whether we can obtain meaningful
analogues of the trace formula by applying the categorical Langlands
correspondence to sheaves on $\Bun_G$ itself.

In \secref{rel trace} we show that there is in fact an analogous
statement which may be interpreted as a geometric analogue of the
``relative trace formula'' (also known as Kuznetsov trace formula,
see, e.g., \cite{Jacquet}). This formula has some interesting
features. First, since it involves Whittaker functionals, only generic
automorphic representations appear on the left hand side and hence,
conjecturally, the non-tempered representations should not appear at
all. Second, conjecturally, the factor $m_\sigma$ in the sum
\eqref{s2} should disappear because each generic $L$-packet is
expected to contain a single irreducible representation with a
non-zero Fourier coefficient.

The price we pay for this is that in the sum \eqref{s2} appears a
weighting factor, $L(\sigma,\on{ad},1)^{-1}$, the reciprocal of the
value of the $L$-function of $\sigma$ in the adjoint representation at
$s=1$ (here we assume that $G$ is simple, but our analysis may be
extended to more general reductive groups). Thus (not counting
automorphisms of $\sigma$), we obtain the sum
\begin{equation}    \label{s3}
\sum_\sigma N_{\sigma} \cdot L(\sigma,\on{ad},1)^{-1}.
\end{equation}
Conjecturally, for non-generic $L$-packets the $L$-function
$L(\sigma,\on{ad},s)$ has a pole at $s=1$ (see \cite{GP,Ich1,Ich2})
and hence these representations should disappear from the sum
\eqref{s3}.

The insertion of this factor in the context of the trace formulas
discussed in \cite{L:BE} (which are very close to the trace formulas
considered here) was originally suggested by P. Sarnak \cite{Sarnak}
and further studied by A. Venkatesh \cite{Ven}, for the group $GL_2$
in the number field context. From our point of view, this factor has a
natural geometric interpretation, as one coming from the {\em
  Atiyah--Bott--Lefschetz fixed point formula} (see Sections
\ref{lfpf}, \ref{abl fixed}, and \ref{geom} for more details).

\medskip

On the right hand side of the relative trace formula we obtain certain
analogues of the Kloosterman sums, which are represented geometrically
by cohomologies of moduli stacks that are similar to the ones
appearing in the geometrization of the ordinary trace formula and have
been previously studied in \cite{Ngo} and \cite{FGV:w}. Understanding
relations between these Kloosterman sums, and the corresponding
cohomologies, for $G$ and the groups $H$ discussed above may give us
another way to approach functoriality.

\medskip

The paper is organized as follows. In \secref{kernel} we introduce the
main ingredient of the geometric trace formula, the sheaf ${\mc
  K}_{d,\rho}$ on the square of $\Bun_G$, and its analogues. This
sheaf, viewed as a kernel, gives rise to a functor acting on the
derived category of sheaves on $\Bun_G$. In \secref{eigenvalues} we
consider the action of this functor on Hecke eigensheaves. We prove
that for large $d$ this functor annihilates those eigensheaves which
do not come from smaller groups $H$ by functoriality. Next, we look at
the right hand (orbital) side of the trace formula in
\secref{orbital}. We show that it corresponds to the cohomology of the
restriction of our sheaf to the diagonal. We study this restriction in
more detail and connect it to Hitchin type moduli stacks of Higgs
bundles. We then formulate a precise conjecture linking the cohomology
in the case that $G={\rm SL}_2$ and $H$ is a twisted torus.

In the second part of the paper, we begin by describing a
conjectural geometrization of the left hand (spectral) side of
\eqref{tr} \secref{leap}. Then we explain in \secref{clc} how a
geometrization of the trace formula appears in the context of the
categorical Langlands correspondence. We state the geometric trace
formula in \conjref{gtf conj}. Finally, we develop the needed
formalism and state the relative geometric trace formula (\conjref{rel
  conj}) in \secref{rel trace}. As a byproduct of this discussion, we
formulate a conjectural generalization of the Atiyah--Bott--Lefschetz
fixed point formula for algebraic stacks in \secref{abl fixed} (see
\conjref{conj abl formula}).

\bigskip

\noindent {\bf Acknowledgments.} This paper was conceived as part of a
joint project with Robert Langlands initiated in \cite{FLN}. We are
deeply grateful to him for his collaboration and for generously
sharing with us his ideas and insights on this subject.

We thank Peter Sarnak for drawing our attention to the relative trace
formula, Atsushi Ichino for communicating his and T. Ikeda's
unpublished conjecture on the norm of the Whittaker functional, and
David Nadler for helpful discussions of the conjectural generalization
of the Atiyah--Bott--Lefschetz fixed point formula for stacks
presented in \secref{abl fixed}. E.F. is grateful to David Ben-Zvi,
David Eisenbud, Alexander Givental, Dick Gross, Michael Harris,
Atsushi Ichino, David Kazhdan, Vincent Lafforgue, Sergey Lysenko, and
Martin Olsson for valuable discussions.

\vspace*{10mm}

\begin{center}

{\bf Part I}

\end{center}

\section{The right hand side of the trace formula}    \label{kernel}

Let $X$ be a smooth projective geometrically connected curve $X$ over
a field $k$ which will be a finite field $\Fq$ or the field of complex
numbers. Denote by $F$ the field of rational functions on $X$ and
$\AD$ the ring of ad\'eles of $F$. For each closed point $x \in |X|$
we denote by $F_x$ the completion of $F$ at $x$ and by $\OO_x$ its
ring of integers.

Let $G$ be a reductive group scheme over $X$ which is a quasi-split
form of a split connected reductive algebraic group $\mathbb G$ over
the field $k$. This means that $G$ becomes split after pull-back to an
\'etale cover $X'$ of $X$.

We wish to find a geometric incarnation of the trace formula
\begin{equation}    \label{general}
\on{Tr} {\mb K} = \underset{G(F)\bs G({\mathbb A)}}\int K(x,x) dx,
\end{equation}
where ${\mb K}$ is the operator on the space of automorphic functions
acting by the formula
$$
({\mathbf K} \cdot f)(x) = \underset{G({\mathbb A)}}\int f(xy) H(y)
dy,
$$
where $H$ is a smooth compactly supported function on $G(\AD)$. We
choose a Haar measure on $G({\mathbb A})$ normalized so that the
volume of the fixed maximal compact subgroup
$$
G(\OO) = \prod_{x \in |X|} G(\OO_x)
$$
is equal to $1$.

We can rewrite ${\mb K}$ as an integral operator corresponding to the
kernel
\begin{equation}    \label{K and H}
K(x,y) = \sum_{a \in G(F)} H(x^{-1} a y),
\end{equation}
acting as follows:
$$
({\mb K} \cdot f)(x) = \underset{G(F)\bs G({\mathbb A)}}\int K(x,y)
f(y) dy.
$$

In what follows we will restrict ourselves to unramified automorphic
representations and we will assume that $K$ is in the (restricted)
tensor product of the spherical Hecke algebras (with respect to
$G(\OO_x)$) over all closed points $x \in |X|$.

Let $\Bun_G$ be the algebraic moduli stack of $G$-bundles on $X$.
Denote by $\on{ker}^1(F,G)$ the kernel of the map
$$
H^1(F,G) \to \prod_{x\in |X|} H^1(F_x,G).
$$
For simplicity, we will assume that $\on{ker}^1(F,G) = 0$. Then we
have
$$\Bun_G(k) = G(F) \bs G(\AD)/G(\OO),$$
where $\Bun_G(k)$ is the set (more precisely, groupoid) of $k$-points
of the moduli stack $\Bun_G$ of $G$-bundles on $X$. Hence the kernel
$K$ may be viewed as a function on $\Bun_G(k) \times \Bun_G(k)$.

Then the right hand side of \eqref{tr} may be rewritten as
\begin{equation}    \label{sum b}
\sum_{P \in \Bun_G(k)} \frac{1}{|\on{Aut}(V)|} K(P,P)
\end{equation}
(see, e.g., \cite{BeDh}). For a class of kernels $K$ described below,
we will give a geometric interpretation of this sum as the trace of
the Frobenius on the cohomology of an $\ell$-adic sheaf on an
algebraic stack over the diagonal in $\Bun_G \times \Bun_G$.

\subsection{The dual group}

Let $^L G$ be the Langlands dual group to $G$,
$$^L G=\check\GG \rtimes \Gamma.$$ Here $\check\GG$ is the dual group
of $\GG$ taken here to be defined over the field of complex numbers if
$k=\C$ but to be over $\Ql$ (where $\ell$ is relatively prime
to $q$) if $k=\Fq$, and $\Gamma$ is the (\'etale) fundamental group
$\pi_1(X,x)$ of $X$ for some geometric point $x$ of $X$ that acts on
$\check \GG$ through a finite quotient. This action determines the
quasi-split form $G$. In particular, if this action is trivial, then
$^L G=\GG\times \Gamma$ and $G=X\times\GG$.

We will fix a splitting of $\GG$ and $\check \GG$ so that in
particular we can talk about dominant (co)weights. The fundamental
group $\Gamma$ acts on the abelian group $\Xb_*$ of all weights (or
characters) of the maximal torus of $\check\GG$, preserving the cone
$\Xb_+$ consisting of the dominant ones. Recall that $\Xb_+$ is
equipped with a partial order: $\mu\geq \mu'$ if and only if
$\mu-\mu'$ is a sum of simple roots of $\check\GG$. Since the action
of $\Gamma$ preserves the set of simple roots, it also preserves the
partial order on $\Xb_+$. It also induces a partial order on the set
$[\Xb_+/\Gamma]$ of orbits of $\Gamma$ on $\Xb_+$.

\subsection{The Hecke stack}    \label{aff Gr}

Recall that if we have two $\GG$-bundles $E,E'$ on the (formal) disc
$\on{Spec} k[[t]]$ which are identified over the punctured disc
$\on{Spec} k\ppart$, we obtain a point in the double quotient
$$\GG[[t]]\bs \GG\ppart/\GG[[t]],$$ or, equivalently, a $\GG[[t]]$-orbit
in the affine Grassmannian $$\Gr = \GG\ppart/\GG[[t]]$$ which is an
ind-scheme over $k$ \cite{BD,MV}. These orbits are called Schubert
cells, and they are labeled by elements $\mu$ the set $\Xb_+$ of
dominant wights of the maximal torus in the dual group $\check\GG$. We
denote the orbit corresponding to $\mu$ by $\Gr_\mu$. We will write
$\inv(E,E')=\mu$ if the pair $(E,E')$ belongs to $\Gr_\mu$. Recall
that $\Gr_{\mu'}$ is contained in the closure $\ovl\Gr_\mu$ of
$\Gr_\mu$ if and only if $\mu\geq\mu'$.

Recall the {\em Hecke stack} ${\mc H}^{\BD} = \Hc^\BD(X,\GG)$ (see
\cite{BD}) that classifies quadruples
$$(x,E,E',\phi),$$ where $x\in X$, $E,E'\in \Bun_{\GG}$ and $\phi$ is
an isomorphism
$$
E|_{X{-}\{x\}}\simeq E'|_{X{-}\{x\}}.
$$
We have two natural morphisms $p, p': {\mc H}^{\BD} \to \Bun_{\GG}$
sending such a quadruple to $E$ or $E'$ and the morphism $s: {\mc H}
\to X$. Since $\Bun_{\GG}$ is an algebraic stack, so is
$\Hc^\BD(X,\GG)$. However, if we fix $E'$, then we obtain an
ind-scheme over $X$, which is called the {\em Beilinson-Drinfeld
  Grassmannian} (see \cite{BD,MV}).

Let ${\mc H}^{\BD}{}'(X,\GG)$ be the stack classifying the quadruples
$$(x,E,E',\phi),$$ where $x\in X$, $E \in \Bun_{\GG}$, $E'$ is a
$\GG$-bundle on the disc $D_x$ around the point $x$, and $\phi$ is
an isomorphism
$$
E|_{D_x^\bullet}\simeq E'|_{D_x^\bullet},
$$
where $D_x^\bullet$ is the punctured disc around $x$. We have a
natural morphism $${\mc H}^{\BD}(X,\GG) \to {\mc H}^{\BD}{}'(X,\GG)$$
(restricting $E'$ to $D_x$ and $\phi$ to $D_x^\bullet$), which is in
fact an isomorphism, according to a strong version of a theorem of
Beauville--Laszlo \cite{BL} given in \cite{BD}, Sect. 2.3.7. Therefore
we obtain that the morphism $$s \times p: {\mc H}^{\BD}(X,\GG) \to X
\times \Bun_{\GG}$$ sending the above quadruple to $(x,E)$ is a
locally trivial fibration with fibers isomorphic to the affine
Grassmannian $\Gr=\GG\ppart/\GG[[t]]$.

For every orbit $\mu \in \Xb_+$ we define the closed substack
$\Hc_\mu$ of ${\mc H}^{\BD}(X,\GG)$ by imposing the inequality
\begin{equation} \label{mu}
\inv_x(E,E')\leq \mu.
\end{equation}
It is a scheme over $X\times \Bun_{\GG}$ with fibers isomorphic
to $\ovl\Gr_\mu$.

\medskip

Suppose now that $G$ is a quasi-split form of $\GG$ whose dual group
$^L G$ is a semi-direct product $\check \GG \rtimes \Gamma$ given by a
finite action of $\Gamma$ on $\check \GG$, fixing a given
splitting. Let $\Gamma'$ be the kernel of this action which is a
normal subgroup of finite index in $\Gamma$.  This is equivalent to a
finite \'etale covering $\pi: X'\to X$ with a geometric point $x'\in
X'$ over $x\in X$ such that the pull-back of $G$ to $X'$ is split. The
group scheme $G$ over $X$ is isomorphic to
$$
X' \underset{\Ga/\Ga'}\times \GG,
$$
and a $G$-torsor on $X$ is nothing but a $\GG$-torsor $\wt{E}$ on $X'$
which is $\Gamma/\Gamma'$-equivariant in the sense that for each
$\gamma \in \Gamma/\Gamma'$ we are given an isomorphism $i_\gamma:
\gamma^*(\wt{E}) \overset\sim\longrightarrow \wt{E}$ (where $\gamma^*$
is induced by both the action of $\gamma$ on $X'$ and on $\GG$), and
these isomorphisms satisfy the compatibility conditions $i_{\gamma_1}
i_{\gamma_2} = i_{\gamma_1 \gamma_2}$ for all $\gamma_1,\gamma_2 \in
\Gamma/\Gamma'$.

We define the Hecke stack ${\mc H}^\BD(G,X)$ as
the algebraic stack which classifies quadruples
$$(x,E,E',\phi),$$ where $x\in X$, $E,E'\in \Bun_{G}$ and $\phi$ is
an isomorphism
$$
E|_{X{-}\{x\}}\simeq E'|_{X{-}\{x\}}.
$$
Equivalently, it classifies $(x,\wt{E},\wt{E}',\phi)$ where $x\in X$,
$\wt{E},\wt{E}'$ are $\Gamma/\Gamma'$-equivariant $\GG$-bundles on
$X'$ (in the above sense), and $\phi$ is an isomorphism
$$
\wt{E}|_{X'{-}\{\pi^{-1}(x)\}}\simeq \wt{E}'|_{X'{-}\{\pi^{-1}(x)\}},
$$
which is invariant under the $\Gamma/\Gamma'$-action on $\wt{E}$ and
$\wt{E}'$.

For every geometric point $x$ of $X$, the restriction of $G$ to the
formal disc $D_x$ around $x$ is isomorphic to $\GG$, and there is an
isomorphism between $G|_{D_x}$ and the constant group scheme over
$D_x$ with the fiber $\GG$, well-defined up to the finite action of
$\Gamma$ on $\GG$. It follows that the morphism ${\mc H}^{\BD}(X,\GG)
\to X \times \Bun_G$ sending the quadruple $(x,E,E',\phi)$ to $(x,E)$
is a fibration over $X \times \Bun_G$ with fibers over geometric
points isomorphic to the affine Grassmannian $\Gr=\GG\ppart/\GG[[t]]$.

Furthermore, if $E$ and $E'$ are $G$-torsors over $D_x$ equipped with
an isomorphism over the punctured disc $D_x^\bullet$, then
${\inv}_x(E,E')$ is well-defined as an orbit $[\mu]$ of $\Gamma$ on
$\Xb_+$. Since we have a partial order on the set $[\Xb_+/\Gamma]$ of
orbits of $\Gamma$ on $\Xb_+$, it makes sense to write
$\inv_x(E,E')\leq [\mu]$. Hence, for every orbit $[\mu]$ of $\Gamma$
acting on $\Xb_+$, we define the closed substack $\Hc_{[\mu]}$ of by
imposing the inequality
\begin{equation} \label{mu1}
\inv_x(E,E')\leq [\mu].
\end{equation}
It is again a fibration over $X\times \Bun_G$ with fibers isomorphic
to the union $\ovl\Gr_{[\mu]}$ of $\ovl\Gr_\mu \subset \Gr$ with $\mu$
in the $\Gamma$-orbit $[\mu]$.

\medskip

The following description of ${\mc H}^\BD(X,G)$ will be useful in
\secref{perv sheaf}. Let $\wt{\mc H}^\BD$ be the stack which
classifies quadruples $(x',\wt{E},\wt{E}',\phi)$ where $x'\in X'$,
$\wt{E}$ is a $\Gamma/\Gamma'$-equivariant $\GG$-bundle on $X'$,
$\wt{E}'$ is a $\GG$-bundle on $D_{x'}$, and $\phi$ is an isomorphism
$$
\wt{E}|_{D_{x'}^\bullet} \simeq \wt{E}'|_{D_{x'}^\bullet}.
$$
The morphism
$$
\wt{\mc H} \to X' \times \Bun_G
$$
sending $(x',\wt{E},\wt{E}',\phi)$ to $(x',\wt{E})$ makes it into an
ind-scheme over $X' \times \Bun_G$ whose fibers are isomorphic to the
affine Grassmannian $\Gr$ of $\GG$. We also have a natural action of
$\Gamma/\Gamma'$ on $\wt{\mc H}^\BD$ lifting the action on $X'$, and
\begin{equation}    \label{quot}
{\mc H}^\BD(X,G) \simeq \wt{\mc H}^\BD/(\Gamma/\Gamma').
\end{equation}

\subsection{Symmetric powers}    \label{sym H}

For every positive integer $d$, we introduce a sort of $d$th symmetric
power of $\Hc_{[\mu]}$. Let $X_d$ be the $d$th {\em symmetric power}
of $X$ which is defined as the scheme-theoretic (also known as GIT)
quotient of $X^d$ by the action of the symmetric group $S_d$. Thus,
$X_d$ classifies effective divisors of degree $d$ on $X$. This is a
smooth algebraic variety (see, e.g., \cite{Shaf}, p.95).

The algebraic stack ${\mc H}_{d,[\mu]}$ over $k$ classifies the data
\begin{equation}    \label{objects}
(D,E,E',\phi),
\end{equation}
where
\begin{equation}    \label{D}
D = \sum_{i=1}^r n_i [x_i]
\end{equation}
is an effective divisor on our curve $X$ of degree $d$ (equivalently,
a point of $X_d$), $E$ and $E'$ are two principal $G$-bundles on $X$,
and $\phi$ is an isomorphism between them over $X-\on{supp}(D)$, such that for
each $i$ the inequality
\begin{equation}\label{n mu}
\inv_{x_i}(E,E')\leq [n_i \mu]
\end{equation}
is satisfied. Consider the morphism
$$
{\mc H}_{d,[\mu]} \to X_d\times \Bun_G
$$
sending the quadruple \eqref{objects} to $(D,E)$. Its fiber over a
fixed $D=\sum_i n_i [x_i]$ and $E\in\Bun_G$ isomorphic to the product
\begin{equation}    \label{product}
\prod_{i=1}^r \ovl\Gr_{[n_i \mu_i]}.
\end{equation}

The algebraic stack ${\mc H}_{d,[\mu]}$ is a closed substack of the
$d$th symmetric power ${\mc H}^{\BD}_d$ of the Hecke stack ${\mc
  H}^{\BD}$ which is defined as the classifying stack of the
quadruples \eqref{objects} as above, but with the condition
$\inv_{x_i}(E,E')\leq [n_i \mu]$ removed.  To verify that these
conditions define a closed substack, it is convenient to introduce a
chain version of ${\mc H}_{[\mu]}$. Let $\Hc_{[\mu]}^d$ be the
algebraic stack classifying
\begin{equation}    \label{chains}
(x_\bullet,E_\bullet,\phi_\bullet),
\end{equation}
where
\begin{itemize}
\item $x_\bullet=(x_1,\ldots,x_d)\in X^d$,
\item $E_\bullet=(E_0,\ldots,E_d)$ with $E_i\in \Bun_G$,
\item $\phi_\bullet=(\phi_1,\ldots,\phi_d)$ where $\phi_i$ is an
  isomorphism $$E_{i-1}|_{X{-}\{x_i\}}\simeq E_i|_{X{-}\{x_i\}}$$
\end{itemize}
such that $\inv_{x_i}(E_{i-1},E_i) \leq [\mu]$. We have a morphism
\begin{equation} \label{compose modif}
\pi^d_\Hc: \Hc_{[\mu]}^d \to \Hc^{\BD}_d
\end{equation}
sending $(x_\bullet,E_\bullet,\phi_\bullet)$ to $(D,E,E',\phi)$, where
$D=\sum_{i=1}^d [x_i]$, $E=E_0$, $E'=E_d$, and $\phi$ is obtained by
composing $\phi_1,\ldots,\phi_d$ on their common domain of definition
$X-\on{supp}(D)$. The morphism (\ref{compose modif}) is proper and
$\Hc_{d,[\mu]}$ is its image; its geometric points are characterized
by the inequality (\ref{n mu}).

We also have a chain version $\Hc^d$ of the entire Hecke stack $\Hc$,
which classifies the same objects as above, but without the condition
that $\inv_{x_i}(E_{i-1},E_i) \leq [\mu]$. The morphism \eqref{compose
modif} extends to
$$
\pi^d_\Hc: \Hc^d \to \Hc^{\BD}_d.
$$

\subsection{Examples} \label{example 1}

Let us discuss examples of the above construction. Let $G$ be the
multiplicative group $\Gm$ over $X$. Recall that $\Gm$-principal
bundles on $X$ are equivalent to line bundles, or invertible
$\Oc_X$-modules. Their classifying stack $\Pic$ has a well-known
structure. The degree of a line bundle induces a bijective map from
the set of connected components of $\Pic$ to $\Z$. The kernel of this
map, denoted by $\Pic_0$, classifying line bundles of degree $0$, is
isomorphic to the quotient of an abelian variety -- namely, the
Jacobian $\Jac$ of $X$ -- by the trivial action of $\Gm$ (which is the
group of automorphisms of line bundles)

In this case, the reduced part of the affine Grassmannian is the
discrete set $\Z$ (however, the non-reduced structure at each point is
highly non-trivial\footnote{Indeed, for a general $k$-algebra $R$ the
quotient $R\ppart/R[[t]]$ is isomorphic to the product of $\Z$ and
the set of polynomials of the form $1+r_{-1} t^{-1}+r_{-2}
t^{-2}+\ldots$, where the $r_i$ are nilpotent elements of $R$.}). The
Hecke stack $\Hc^{\BD}$ classifies quadruples
$$(x,L,L',\phi),$$ where $(x,L,L')\in X\times \Pic \times \Pic$ and
$\phi$ is an isomorphism
$$L|_{X{-}\{x\}} \simeq L'|_{X{-}\{x\}}.$$ At the level of $k$-points
we then have $L'=L(\mu x)$ for some integer $\mu \in
\Z$.\footnote{However, this is not true at the level of $R$-points,
  where $R$ is a ring with nilpotents, as explained in the previous
  footnote.} Hence the {\em reduced} part of $\Hc^{\BD}$ is $$\Pic
\times X\times \Z$$ in this case. For every $\mu\in \Z$, which is the
set $\Xb_+$ in this case, we have the inclusion (of the reduced parts,
which we denote by the same symbols by abuse of notation)
$$\Hc_\mu=\Pic\times X\times \{\mu\} \subset \Hc^{\BD}.$$ 
The other projection onto $\Pic$ is $(L,x,\mu) \mapsto L'=L(\mu x)$.

Since $\Gm$ is commutative, $\Hc^{\BD}_d$ can be identified with
$$\Hc^{\BD}_d=\Pic\times (X \times \Z)_d$$ where $(X \times \Z)_d$
denotes the set of orbits of the symmetric group $S_d$ on
$(X\times\Z)^d$. We have two maps $(X \times \Z)_d \to \Z_d$
and $(X \times \Z)_d \to X_d$. The $d$th symmetric power of $\Z$
consists of elements of the form $((d_1,\mu_1),\ldots,(d_r,\mu_r))$
where $d_1,\ldots,d_r$ are positive integers satisfying
$d_1+\cdots+d_r=d$ and $\mu_1 < \cdots < \mu_r$ is a strictly
increasing sequence of integers. The fiber of $(X \times \Z)_d$ over
$((d_1,\mu_1),\ldots,(d_r, \mu_r))$ is $X^{(d_1)} \times\cdots\times
X^{(d_r)}$ projecting onto $X_d$ by the obvious map of adding the
divisors. In particular, for the element $(d,\mu)\in \Z_d$, we have
the component $$\Hc_{d,\mu}=\Pic\times X_d.$$

\medskip

Our second example is a non-split one-dimensional torus $H$ over
$X$. Such a torus is given by a $\mu_2$-torsor (equivalently, \'etale
double cover) $\pi:X' \to X$. Let us denote by $\tau$ the non-trivial
involution on $X'$. The group $\Gamma' \subset \Gamma$ has index 2 in
this case and $\Gamma/\Gamma' = \{ 1,\tau \}$. We assume that $X'$ is
geometrically connected so that $H$ remains non-split over $X\otimes
\bar k$. An $H$-torsor on $X$ is a line bundle $L$ over $X'$ whose
norm down to $X$ is equipped with a trivialization
$\beta:N_{X'/X}(L) \overset{\sim}\longrightarrow \Oc_X$.

The group of connected component of $\Bun_H$ is $\Z/2\Z$.
The neutral component of $\Bun_H$ is a quotient of an abelian variety --
namely, the Prym variety of $X'$ -- by the trivial action of $\Z/2\Z$
which is the group of automorphism of any pair $(L,\beta)$ as
above.

The set $\Xb_+=\Z$ is acted on by the involution $\tau(\mu)=-\mu$. The
set of orbits $[\Xb_+/\Gamma]$ may therefore be indexed by
non-negative integers $m\in \mathbb Z_{\geq 0}$. The moduli stack
$\Hc^\BD$ (again, the reduced part) decomposes into connected
components
\begin{equation}\label{Grassmannian H}
\Hc^\BD=\bigsqcup_{m\in \mathbb Z_{\geq 0}} \Hc_m
\end{equation}
with $\Hc_0=\Bun_H \times X$ and $\Hc_m=\Bun_H \times X'$ for $m\geq
1$. The second projection onto $\Bun_H$ may be described as follows:
for $m=0$ we map $(L,x)\mapsto L$, and for $m\geq 1$ we map $(L,x')
\mapsto L(m(x'-\tau(x')))$.

For every $d\geq 1$ we have the following description of $\Hc_d$:
\begin{equation}\label{Grassmannian H d}
\Hc_{d}=\Bun_H \times (X\sqcup (X'\times \mathbb Z_{>0}))_d.
\end{equation}
The connected components of $(X\sqcup (X'\times \mathbb Z_{>0}))_d$
are the fibers over elements of $(\mathbb Z_{>0})_d$, which are
$r$-tuples $((d_0,m_0),(d_1,m_1),\ldots,(d_r,m_r))$ with
$0=m_0<m_1<\cdots<m_r$, where $d_0 \geq 0$ and $d_1,\ldots,d_r > 0$
are integers such that $d_0+d_1+\cdots+d_r=d$. The corresponding
component of ${\mc H}_d$ is
$$\Bun_H \times X_{d_0}\times X'_{d_1}\times \cdots \times X'_{d_r}.$$

\subsection{The perverse sheaf}    \label{perv sheaf}

Let $\rho$ be a finite-dimensional representation of $^L G$. For
simplicity, we will restrict ourselves to irreducible
representations. If $G$ is split, for example, then such $\rho$ must
be of the form $\rho=\rho_\mu \otimes \rho_{\Lc}$, where $\rho_\mu$ is
the irreducible representation of $\check\GG$ with highest weight
$\mu$ and $\rho_{\Lc}$ is the representation of $\Gamma$ attached to
an irreducible local system $\Lc$ on $X$.

For any semi-simple finite-dimensional representation $\rho$ of $^L
G$, we will define perverse sheaves $\wF_{\rho}$ on $\Hc^\BD$ and
${\mc K}_{d,\rho}$ on $\Hc_d$. They may be constructed more
conceptually using an equivalence of categories suggested in the
Appendix to \cite{FGV:lc}, but for our purposes the explicit
construction given below will suffice.

First, we assume that $G$ is split and $\rho$ is an irreducible
representation of $^L G=\check\GG \times \Gamma$ which is trivial on
$\Gamma$. Recall that we have defined a closed substack $\Hc_\mu$ of
$\Hc^\BD$ such that the projection $\Hc_\mu\to X\times \Bun_G$ is a
locally trivial fibration with fibers isomorphic to the Schubert
variety $\ol\Gr_\mu$.

Using the {\em geometric Satake correspondence} \cite{MV}, we
associate to each irreducible representation $\rho = \rho_\mu$ of $^L
G$ of highest weight $\mu \in \Xb_+$ the sheaf
$$
{\mc K}_\rho={\mc K}_\mu = \IC(\Hc_\mu)[-\dim(X\times
\Bun_G)](-\dim(X\times \Bun_G)/2).
$$
With the above cohomological shift and Tate twist, locally over
$X\times \Bun_G$, ${\mc K}_\rho$ is isomorphic to
$$\underline\Ql \boxtimes \IC(\ol\Gr_\mu),$$ where
$\IC(\ol\Gr_\mu)$ denotes the intersection cohomology sheaf of
$\ol\Gr_\mu$ and $\underline\Ql$ is the constant sheaf on $X\times
\Bun_G$. Hence
${\mc K}_\rho$ is perverse and pure of weight 0 along the fibers of the
projection to $X\times \Bun_G$. By abuse of notation, we will use the
same symbol for $i_* {\mc K}_\rho,$ where $i:\ol\Hc_\mu \to \Hc$,
so that ${\mc K}_\rho$ will also be viewed as a sheaf on $\Hc$.

We define the {\em Hecke functor} ${\mathbb H}_\rho={\mathbb H}_\mu$
as the integral transform corresponding to the kernel ${\mc K}_\rho$
(see \cite{BD}):
$$
{\mathbb H}_\rho({\mc F}) = (s \times p)_*(p'{}^*({\mc F}) \otimes
{\mc K}_\rho).
$$
For $x \in |X|$, let $\Hc_x$ be the fiber of $\Hc$ over $x$, and
$p_x,p'_x: \Hc_x \to \Bun_G$ the corresponding morphisms. Denote by
${\mc K}_{\rho,x}$ the restriction of ${\mc
K}_\rho$ to $\Hc_x$. Define the functor ${\mathbb H}_{\rho,x}$ by the
formula
$$
{\mathbb H}_{\rho,x}({\mc F}) = p_{x*}(p'_x{}^*({\mc F}) \otimes {\mc
  K}_{\rho,x}).
$$
The function corresponding to the sheaf ${\mc K}_{\rho,x}$ is
the kernel $K_{\rho,x}$ of the Hecke operator ${\mb
  H}_{\rho,x}$ and so the functor ${\mathbb H}_{\rho,x}$ on the
derived category $D(\Bun_G)$ of $\ell$-adic sheaves on $\Bun_G$ is a
geometric analogue of the Hecke operator ${\mb H}_{\rho,x}$.

\medskip

Next, suppose that $G$ is split and $\rho$ is a more general
irreducible representation of $^L G=\check\GG \times \Gamma$ of the
form $\rho=\rho_\mu \otimes \rho_\Lc$, where $\rho_\mu$ is the
irreducible representation of $\check\GG$ of highest weight $\mu \in
\Xb_+$ and $\rho_\Lc$ is the representation of $\Gamma$ attached to an
irreducible local system $\Lc$ on $X$. Now we put
\begin{equation}\label{sheaf}
\wF_\rho=\wF_\mu \otimes s^*(\Lc),
\end{equation}
where $s$ is the obvious projection onto $X$. By abuse of notation, we
will use the same symbol for $i_* \wF_\rho$ where $i:\ol\Hc_\mu \to
\Hc^\BD$, so that $\wF_\rho$ could also be viewed as a sheaf on
$\Hc^\BD$ (again, it is perverse on the entire $\Hc^\BD$ only up to a
cohomological shift, but perverse along the fibers of the projection
to $X\times\Bun_G$).

\medskip

Finally, consider the non-split case, following the notation of
\secref{aff Gr}. We have an exact sequence
$$1\to \check\GG \times \Gamma' \to \check\GG\rtimes \Gamma \to
\Gamma/\Gamma' \to 1.$$ Let $\rho$ be an irreducible representation of
$\check \GG \rtimes \Gamma$. Suppose that its restriction to
$\check\GG \times \Gamma'$ is semi-simple, so that
\begin{equation}    \label{decomp}
\rho|_{\check\GG \times \Gamma'} = \bigoplus_{\mu,\Lc}
(\rho_\mu\otimes \rho_\Lc)^{\oplus n_{\mu,\Lc}}
\end{equation}
over a finite collection of pairs $(\mu,\Lc)$, where $\mu\in\Xb_+$ and
$\Lc$ is an irreducible local system on $X'$. Here $\rho_\mu$ is an
irreducible representation of $\check\GG$ of highest weight,
$\rho_\Lc$ is the representation of $\Gamma'$ attached to $\Lc$, and
the integer $n_{\mu,\Lc}$ is the multiplicity.

We attach to this $\rho$ a sheaf $\wF_\rho$ on $\Hc^\BD(X,G)$
as follows. First, we construct, in the same way as above, the sheaf
$\wF'_\rho$ on the stack $\wt\Hc^\BD$ introduced in \secref{aff
Gr}, using the fact that this stack is a locally trivial fibration
over $X' \times \Bun_G$ with fibers isomorphic to the affine
Grassmannian $\Gr$ of $\GG$. Since $\rho$ is a representation of
$\check \GG \rtimes \Gamma$ (and not just $\check\GG \times \Gamma'$),
we obtain that $\wF'_\rho$ is a $\Gamma/\Gamma'$-equivariant
sheaf on $\wt\Hc^\BD$. Hence it descends to a sheaf on the quotient of
$\wt\Hc^\BD$ by $\Gamma/\Gamma'$, which is nothing but $\Hc^\BD(X,G)$,
by formula \eqref{quot}. This gives us the desired sheaf ${\mc
K}_\rho$ on $\Hc^\BD(X,G)$.

We note that this construction yields a generalization of the
geometric Satake correspondence to the quasi-split case.

\subsection{Symmetric power}    \label{symmetric power}

Next, we will construct the symmetric power $\wF_{d,\mu}$ of
$\wF_{\mu}$ for any positive integer $d$.

First, we recall the construction of the symmetric power of a local
system $\Lc$ on $X$. From our point of view, this is a particular case
of the general construction corresponding to the trivial $G$ (or
general $G$ but $\mu=0$).

Let $\Lc$ be a rank $n$ {\em local system} on a curve $X$ over a field
$k$ (which, as before, is either $\C$ or $\Fq$); that is, a rank $n$
locally constant sheaf ($\ell$-adic, if $k=\Fq$) on $X$.  Recall from
\secref{sym H} that $X_d = X^d/S_d$ is the $d$th symmetric power of
$X$. This is a smooth algebraic variety defined over $k$, whose
$k$-points are effective divisors on $X$ of degree $d$.

We define a sheaf on $X_d$, denoted by $\Lc_d$ and called the {\em
$d$th symmetric power of $\Lc$}, as follows:
\begin{equation}    \label{stalk}
{\Lc}_d = \left(\pi^d_*(\Lc^{\boxtimes d})\right)^{S_d},
\end{equation}
where $\pi^d: X^d \to X_d$ is the natural projection. The stalks of
  $\Lc_d$ are easy to describe: they are tensor products of symmetric
  powers of the stalks of $\Lc$. The stalk ${\mc L}_{d,D}$ at a
  divisor $D = \sum_i n_i [x_i]$ is
$$
\Lc_{d,D} = \bigotimes_i S^{n_i}({\mc \Lc}_{x_i}),
$$
where $S^{n_i}({\mc \Lc}_{x_i})$ is the $n_i$-th symmetric power of
the vector space ${\mc \Lc}_{x_i}$.  In particular, the dimensions of
the stalks are not the same, unless $n=1$. In the case when $n=1$ the
sheaf $\Lc_d$ is in fact a rank 1 local system on $X_d$. For all $n$,
$\Lc_d$ is actually a perverse sheaf on $X_d$ (up to cohomological
shift), which is irreducible if and only if ${\mc L}$ is irreducible.

Let denote by $X^{d,\circ}$ the open subscheme of $X^d$ defined by the
equation $x_i\not= x_j$ for all indices $i\not=j$. The action of the
symmetric group $S_d$ on $X^d$ preserves $X^{d,\circ}$ and the
quotient of $X^{d,\circ}$ by $S_d$ is an open subscheme $X_{d,\circ}$
of $X_d$ which classifies multiplicity free effective divisors of
degree $d$ on $X$.

\begin{lem} 
$\Lc_d[d]$ is a perverse sheaf on $X_d$ which is the intermediate
extension from a local system on $X_{d,\circ}$.
\end{lem}

\begin{proof} The natural projection $\pi^d: X^d \to X_d$ is a finite
  flat map which is \'etale over the open subscheme $X_{d,\circ}$ of
  $X_d$. As a finite map, $\pi^d$ is in particular small in the sense
  of Goresky and MacPherson. If $p: Y \to Z$ is a small map and ${\mc
    F}$ a local system on $Y$, then $p_*({\mc F})[\dim Y]$ is a
  perverse sheaf on $Z$ isomorphic to the intermediate extension of
  its restriction to any open dense subset of $Z$. Hence we obtain
  that this is true for the sheaf $\pi^d_*(\Lc^{\boxtimes d})$ on
  $X_d$. Taking the $S_d$-invariants, we obtain the statement of the
  lemma.
\end{proof}

This lemma suggests the following general construction of ${\mc
F}_{d,\rho}$. We note that a similar construction was proposed by
Laumon in \cite{Laumon:Duke} in a slightly different context (for
$GL_n$). Recall that we have a map
$$\pi^d_\Hc:\Hc^{d}_{[\mu]} \to \Hc_{d,[\mu]}$$
defined in (\ref{compose modif}). We have in fact a commutative diagram
\begin{equation}\label{diagram} 
	\xymatrix{\ar @{} [dr] |{}
\Hc^{d}_{[\mu]} \ar[d] \ar[r]^{\pi^d_\Hc} & \Hc_{d,[\mu]} \ar[d] \\ 
X^d \ar[r]_{\pi^d} & X_d }
\end{equation}
which is Cartesian over the open subscheme $X_{d,\circ}$. Recall that
over the same open subscheme, the chain version $\Hc^{d}_{[\mu]}$
defined in \secref{sym H} is canonically isomorphic to the $d$th
Cartesian power of the map
$$\Hc_{[\mu]} \to \Bun_G.$$ By restriction, $(\wF_\rho)^{\boxtimes
d}$ defines a perverse sheaf on
$$\Hc^{d}_{[\mu]} \us{X^d}\times {X^{d,\circ}}$$ which is equipped
with an action of $S_d$. By finite \'etale Galois descent, we get a
perverse sheaf over
$$\Hc_{[\mu]} \us{X^d}\times X_{d,\circ}.$$ Now, we can apply the
functor of intermediate extension to obtain a perverse sheaf
\begin{equation} \label{sym power sheaf}
\wF_{d,\rho} \mbox{ on } \Hc_{d,[\mu]}.
\end{equation}
In the above discussion, we have neglected a shift, since $\wF_\rho$
defined as in (\ref{sheaf}) is only a perverse sheaf after a
shift. Recall that the effect of this shift makes the restrictions of
$\wF_\rho$ to the fibers of the map $\Hc_{[\mu]} \to X\times \Bun_G$
perverse sheaves. With the above definition, the restrictions of
$\wF_{d,\rho}$ to the fibers of the maps
$$\Hc_{d,[\mu]} \to X_d \times \Bun_G$$ over the open subset
$X_d^{\circ}\times \Bun_G$ are also perverse sheaves. We will see that
this is also true for all other fibers.

The above definition of $\wF_{d,\rho}$ is probably the quickest one,
but it is impractical for many purposes as the intermediate extension
is not very explicit. The following lemma, essentially proved by
Mirkovi\'c and Vilonen in \cite{MV}, uses the whole diagram
(\ref{diagram}) to make the definition of $\wF_{d,\rho}$ more
explicit.

\begin{lem} \label{small}
The morphism $\pi^d_\Hc$ of the diagram {\rm (\ref{diagram})} is small
in the stratified sense of \cite{MV}. Moreover, for all $x_\bullet\in
X^d$ mapping to $D\in X_d$, let $\Hc^{d}_{x_\bullet}$ be the fiber of
$\Hc^{d}$ over $x_\bullet$ and $\Hc^\BD_{d,D}$ the fiber of
$\Hc^\BD_{d}$ over $D$. Then the restriction of $\pi^d_\Hc$ to the
fibers of $$\Hc^{d}_{x_\bullet} \to \Hc^\BD_{d,D}$$
is semi-small in the stratified sense. 
\end{lem}

Assume now till the end of this subsection that $G$ is split, $^L
G=\check\GG\times \Gamma$ and $\rho=\rho_\mu\otimes \rho_\Lc$. In this
case the construction of $\wF_{d,\rho}$ can be made more explicit
along the lines of Springer's construction of representations of the
symmetric group. Consider the closed substack $\Hc_\mu^d$ of $\Hc^{d}$
defined in \secref{sym H} and the perverse sheaf
\begin{equation}\label{before push-forward}
\IC(\Hc_\mu^d)\otimes {\rm pr}_{X^d}^*(\Lc^{\boxtimes d})
\end{equation}
on $\Hc^d_\mu$. Next, consider the push-forward
\begin{equation}\label{push-forward}
(\pi^d_\Hc)_*( \IC(\Hc_\mu^d)\otimes {\rm pr}_{X^d}^*(\Lc^{\boxtimes d}))
\end{equation}
to $\Hc_{d,\mu}$ which is contained in $\Hc_d^\BD$. It follows from
the first assertion of Lemma \ref{small} that (\ref{push-forward}) is
a perverse sheaf and is the intermediate extension of its restriction
to any non-trivial open subset; in particular, to $\Hc_{d,\mu}
\us{X^d}\times{X^{d,\circ}}$. This defines an action of $S_d$ on
(\ref{push-forward}). By taking the invariant part of this action and
making cohomological shift by $-\dim(X^d\times\Bun_G)$, we obtain
another definition of $\wF_{d,\rho}$:
\begin{equation}\label{invariant}
\wF_{d,\rho}=\left((\pi^d_\Hc)_*(\IC(\Hc_\mu^d)\otimes {\rm
pr}_{X^d}^* \Lc^{\boxtimes d})\right)^{S_d}[-(d+\dim \Bun_G)](-(d+\dim
\Bun_G)/2).
\end{equation}
The first assertion of Lemma \ref{small} guarantees the equivalence
of the two definitions.

The shift by $-\dim(X^d\times\Bun_G)$ in (\ref{before push-forward})
makes our sheaf perverse along the fibers of the map $\Hc^d_\mu \to
X^d\times\Bun_G$. According to the second assertion of Lemma
\ref{small}, the restriction of $\wF_{d,\rho}$ to the fibers of the
map $\Hc_{d,\mu}\to X_d\times\Bun_G$ is also perverse.

Now it is possible to describe the fibers of $\wF_{d,\rho}$ using the
geometric Satake equivalence \cite{MV}. Fix an effective divisor of
degree $d$
$$
D=\sum_{i=1}^r n_i [x_i]\in X_d, \qquad x_i\not= x_j, \quad
\sum_{i=1}^r n_i=d.
$$
Fix a principal $G$-bundle $E\in \Bun_G$ and choose trivializations of
$E$ on the formal discs $D_{x_i}$. Then the fiber of $\Hc^\BD_d$ over
$(D,E)$ is the product of copies of the affine Grassmannian attached
to the points $x_i$,
\begin{equation}\label{factor}
\Hc^\BD_d(D,E)=\prod_{i=1}^r \Gr_{x_i}.
\end{equation}
The fiber of $\Hc_{d,\mu}$  over $(D,E)$ is
\begin{equation}\label{factor mu}
\Hc_{d,\mu}(D,E)=\prod_{i=1}^r \ol\Gr_{x_i, n_i\mu}.
\end{equation}
The restriction $\wF_{d,\rho}$ to this fiber is the external tensor
product of spherical perverse sheaves
\begin{equation}\label{factor F}
\bigotimes_{i=1}^r \pr_i^* S^{n_i}(\IC(\ovl \Gr_{x_i,\mu})\otimes
\Lc_{x_i}),
\end{equation}
where $S^{n_i}$ is the $n_i$-th symmetric power based on the
symmetrical monoidal structure of the category of spherical perverse
sheaves on the affine Grassmannian. In other words, $S^{n_i}(\IC(\ovl
\Gr_{x_i,\mu})$ is the perverse sheaf on the affine Grassmannian which
corresponds to the $n_i$th symmetric power of $\rho$ under the
geometric Satake equivalence \cite{MV} ($\pr_i$ denotes the projection
of \eqref{factor mu} onto the $i$th factor). This formula follows from
the second assertion of Lemma \ref{small} and from the fact that the
definition of Drinfeld's commutativity constraint used in \cite{MV}
was also based on this lemma.

Let $K_{d,\rho}$ be the function on $\Bun_G(k) \times \Bun_G(k)$
corresponding to the sheaf $\wF_{d,\rho}$. Then it is given by the
formula
\begin{equation}    \label{kernel of kdrho}
K_{d,\rho} = \sum_{D=\sum_i n_i[x_i] \in
 X_d(\Fq)} \quad \prod_i K_{\on{Sym}^{n_i}(\rho),x_i},
\end{equation}
where $K_{\rho,x}$ is the kernel of the Hecke operator ${\mb
  H}_{\rho,x}$.

\subsection{Example} \label{example 2}

Consider the case of $G=X \times {\rm SL}_2$. The dual group is $^L
G=\check\GG \times \Gamma$ where $\check\GG={\rm PGL}_2$. Here
$\check\GG$ is defined over $\C$ if $k=\C$ and over $\Ql$ if $k$ is a
finite field. We fix an irreducible three-dimensional representation
$$\rho=\rho_\mu \otimes \rho_\Ec$$ of $^L G$ as follows. Let $X'\to X$
be an \'etale double cover of $X$ which is geometrically
connected. Choose a geometric point $x'$ of $X'$ over the given point
$x\in X$. This gives rise to a non-trivial character $\rho_\Ec$ of
order two of $\Gamma$.  The first component $\rho_\mu$ is the
irreducible three-dimensional representation of ${\rm PGL}_2$. Let $V$
be the tautological two-dimensional representation of ${\rm
  GL}_2$. Then we have the following formula for the lifting of
$\rho_\mu$ to a representation of ${\rm PGL}_2$:
\begin{equation}\label{three dimension}
\rho_\mu=S^2 V\otimes {\det}(V)^{-1}.
\end{equation}
If we chose the splitting of ${\rm SL}_2$ with the diagonal torus 
$$T=\{{\rm diag}(t,t^{-1})\}$$ and the subgroup of upper triangular
matrices as Borel subgroup, then $\mu$ is simply the cocharacter of
$T$ given by $t\mapsto {\rm diag}(t,t^{-1})$.

For every $d>0$, $\Hc_{d,\mu}$ classifies the data
$$(D,E,E',\phi)$$ where $D\in X_d$ is an effective divisor of degree
$d$, $E,E'$ are rank two vector bundles on $X$ with trivialized
determinant and $\phi$ is an isomorphism of the restrictions of the
vector bundles $E$ and $E'$ to $X-\on{supp}(D)$ such that
$\det(\phi)=1$, and $\phi$ can be extended as a morphism of vector
bundles $E\to E'(D)$.

The representation $\rho_\mu$ corresponds to the intersection
cohomology complex with the appropriate shift
$$\wF_{d,\rho_\mu}\IC(\Hc_{d,\mu})[-(d+\dim(\Bun_G))](-(d+\dim(\Bun_G))/2).$$
The sheaf $\wF_{d,\rho}$ is
$$\wF_{d,\rho}=\IC(\Hc_{d,\mu}) \otimes
\pr_{X_d}^*(\Ec_d)[-(d+\dim(\Bun_G))](-(d+\dim(\Bun_G))/2),$$ where
$\pr_{X_d}$ is the projection on $X_d$ and $\Ec_d$ is the $d$th
symmetric power of $\Ec$, which is a rank one local system on $X_d$ of
order two corresponding to $\rho_{\Ec}$.

\medskip

Next, we consider a closely related example for a one-dimensional
twisted torus $H$ over $X$ as in \secref{example 1}. The dual group
$^L H$ is $\Gm \rtimes \Gamma$, where $\Gamma$ acts through a quotient
$n_H:\Gamma\to \Z/2\Z$, and $\Z/2\Z$ acts by the formula
$\tau(t)=t^{-1}$ for the non-trivial element $\tau \in \Z/2\Z$ and
$t\in\Gm$. Now, we have a homomorphism $^L H \to ^L G$
\begin{equation}    \label{hom L}
\Gm \rtimes \Gamma \to {\rm PGL}_2 \times \Gamma
\end{equation}
mapping $\Gm$ to the diagonal torus of ${\rm PGL}_2$ and an element
$\sigma\in \Gamma$ to the element $$(w_0^{n_H(\sigma)},\sigma)\in {\rm
PGL}_2 \times \Gamma$$ where $w_0$ is the permutation matrix and
$n_H(\sigma)$ denotes the image of $\sigma$ in $\Z/2\Z = \{ 0,1 \}$.

Consider the representation $\rho=\rho_\mu\otimes \rho_\Ec$ of $^L G$,
where $\rho_\mu$ is the irreducible three-dimensional representation
of ${\rm PGL}_2$ and $\rho_\Ec$ is the one-dimensional representation
of $\Gamma$ of order two as above. If we restrict $\rho$ to $^L H$, it
decomposes into a direct sum
$$\rho_H=\rho|_{^L H}=\rho_0 \oplus \rho_1,$$
where $\rho_1$ is a two-dimensional representation of $^L H$
and $\rho_0$ is the one-dimensional representation which is trivial on
$\Gm$, and on the factor $\Gamma$ is given by the formula
$$\rho_0=\rho_\Ec\otimes \rho_{\Ec_H},$$
$\rho_{\Ec_H}$ being given by the character
$n_H:\Gamma\to\Z/2\Z$. In particular, $\rho_0$ is the trivial
representation if $\Ec=\Ec_H$.

Assume from now on that $\Ec=\Ec_H$. Then the representation $\rho_H$
of $^L H$ corresponds to the direct sum
$$\wF_{\rho_H}=\wF_0 \oplus \wF_1,$$ where $\wF_1=\und\Ql|\Hc_0$ is
the constant sheaf on $\Hc_0=X\times\Bun_H$ and $\wF_2=\und\Ql|\Hc_1$
is the constant sheaf on $\Hc_1=X'\times \Bun_H$, in the notation of
(\ref{Grassmannian H}). For every positive integer $d$, we have
$$\wF_{d,\rho_H}=\bigoplus_{d_0+d_1=d} \und\Ql|(\Bun_H\times
X_{d_0}\times X'_{d_1}),$$
where the sum is over all non-negative integers $d_0$ and $d_1$ such
that $d_0+d_1=d$. Here $\Bun_H\times X_{d_0}\times X'_{d_1}$ is a
connected component of $\Hc_{H,d}$ described in (\ref{Grassmannian
H d}).

\subsection{Integral transforms}    \label{int tr}

We have defined a morphism
$$\Hc^\BD_d\to X_d\times \Bun_G\times\Bun_G$$ and a perverse sheaf
$\wF_{d,\rho}$ on $\Hc^\BD_d$ attached to any semi-simple
finite-dimensional representation $\rho$ of $^L G$. Let us denote by
$p_d$ and $p'_d$ the two projections $\Hc^\BD_d \to \Bun_G$ mapping
the quadruple \eqref{objects} to $E$ and $E'$. We use the sheaf
$\wF_{d,\rho}$ to define an integral transform functor ${\mathbb
  K}_{d,\rho}$ on the derived category $D(\Bun_G)$ of $\ell$-adic
sheaves on $\Bun_G$ by the formula
\begin{equation}
{\mathbb K}_{d,\rho}({\mc F}) =p_{d!}(p'_d{}^*({\mc
    F}) \otimes \wF_{d,\rho}).
\end{equation}
Functors of this kind were first introduced in \cite{FGV:lc} in the
case of $GL_n$ under the name ``averaging functors.''  An example is
the functor $\on{Av}^d_E$ of \cite{FGV:lc} which plays the role of a
``projector'' on the Hecke eigensheaf corresponding to an irreducible
rank $n$ local system $E$ on $X$, in the category of ${\mc
  D}$-modules, or $\ell$-adic sheaves, on $\Bun_G$. The functors $\K$
corresponding to these kernels are combinations of the Hecke functors,
in which the positions of the points at which these functors are
applied are allowed to vary over the $d$th symmetric power $X_d$
of $X$.

We generalize this construction to a larger family of integral
transforms on $D(\Bun_G)$ as follows. For any scheme over the
symmetric power $$S\to X_d$$
we form the Cartesian product
$$\Hc^\BD_{d,S}=S \us{X_d}\times \Hc^\BD_d.$$ Denote by $\wF_{d,\rho,S}$
the restriction of $\wF_{d,\rho}$ to $\Hc^\BD_{d,S}$. Let $p_{S}$
and $p'_{S}$ the two projections from $\Hc^\BD_{d,S}$ to $\Bun_G$. Then
we have an integral transform
\begin{equation}
{\mathbb K}_{d,\rho,S}({\mc F}) =p_{S!}(p'_{S}{}^*({\mc F}) \otimes
\wF_{d,\rho,S}).
\end{equation}
If $S$ is a point of $X_d$, this is the usual Hecke functor. If
$S=U_d$, where $U$ is an open subset of $X$, we have an averaging
operator over $U$ that will be useful in the study of the ramified
case.

\subsection{Vector space}    \label{v sp}

Let us form the Cartesian square
\begin{equation}    \label{Cartesian}
\begin{CD}
{\mc M}_{d} @>{{\Delta}_\Hc}>> {\mc H}_{d} \\
@VV{{\mb p}_\Delta}V @VV{{\mb p}}V \\
\Bun_G @>{\Delta}>> \Bun_G \times
\Bun_G
\end{CD}
\end{equation}
where $\Delta$ is the diagonal morphism. Thus, ${\mc
  M}_{d}$ is the fiber product of $X_d \times \Bun_G$ and ${\mc
  H}_{d}$ with respect to the two morphisms to $X_d \times \Bun_G
\times \Bun_G$.

Let ${\mc K} = {\mc K}_{d,\rho,S}$ be a sheaf of the type introduced
in the previous section and ${\mathbb K}$ the corresponding integral
transform functor on the derived category $D(\Bun_G)$ of sheaves on
$\Bun_G$. Let
$$
\ol{\mc K} = {\mb p}_*({\mc K})
$$
(note that ${\mathbf p}$ is proper over the support of ${\mc
  K}$). This is a sheaf on $\Bun_G \times \Bun_G$ which is the kernel
of the functor $\K$.

Note that the sheaf $\Delta_!(\underline{\Ql})$ on $\Bun_G
\times \Bun_G$ is the kernel of the identity functor $\on{Id}$ on the
category $D(\Bun_G)$. Hence the vector space
\begin{equation} \label{rhs-new1}
\on{RHom}(\Delta_!(\underline{\Ql}),\ol{\mc K}),
\end{equation}
may be viewed as a concrete realization of the ``categorical trace''
$$
\on{RHom}(\on{Id},\K)
$$
of the functor $\K$.

By adjunction and base change, the space \eqref{rhs-new1} is
isomorphic to
\begin{equation}    \label{adjunction}
H^\bullet(\Bun_G,\Delta^! {\mb p}_*({\mc K})) =
H^\bullet(\Bun_G,{\mb p}_{\Delta *} \Delta_{\Hc}^!({\mc K})) =
H^\bullet({\mc M}_d,\Delta_{\Hc}^!({\mc K})).
\end{equation}

Let ${\mathbb D}$ be the Verdier duality on ${\mc
H}_{d,\mu}$. It follows from the construction and the fact that
${\mathbb D}(\on{IC}({\mc H}_\mu)) \simeq \on{IC}({\mc H}_\mu)$ that
\begin{equation} \label{verdier}
{\mathbb D}({\mc K}_{d,\rho_\mu
    \otimes \rho_{\Lc}}) \simeq {\mc K}_{d,\rho_\mu \otimes
    \rho^*_{\Lc}}[2(d+\dim \Bun_G)] (d+\dim \Bun_G).
\end{equation}
Therefore, up to a shift and Tate twist (and replacing $\Lc$ by
$\Lc^*$), the last space in \eqref{adjunction} is isomorphic to
\begin{equation}    \label{comp supp}
H^\bullet({\mc M}_d,{\mathbb D}(\Delta_{\Hc}^*({\mathcal K}))) \simeq
H^\bullet_c({\mc M}_d, \Delta_{\Hc}^*({\mathcal K}))^*,
\end{equation}
where $H^\bullet_c(Z,{\mc F})$ is understood as $f_!({\mc F})$, where
$f: Z \to \on{pt}$. Here we use the results of Y. Laszlo and M. Olsson
\cite{LO} on the six operations on $\ell$-adic sheaves on algebraic
stacks and formula \eqref{verdier}.

Using the Lefschetz formula for algebraic stacks developed by
K. Behrend \cite{Behrend}, we find that, formally, the trace of the
(arithmetic) Frobenius on
\begin{equation}    \label{rhs-new} 
H^\bullet(\Bun_G,\Delta^!(\ol{\mc K})) = H^\bullet({\mc
  M}_d,\Delta_{\Hc}^!({\mc K}))
\end{equation}
is equal (up to a power of $q$) to the sum \eqref{sum
  b}, which is the right hand side of the trace formula \eqref{tr}.

\smallskip

For this reason, we will view the vector space \eqref{rhs-new} as a
geometric incarnation of the right hand side of the trace formula
\eqref{tr}.

\smallskip

There are some technical issues that have to be dealt with in order to
make this precise, because ${\mc M}_{d}$ is not an algebraic stack of
finite type. In general, in order to make the trace of the Frobenius
on \eqref{rhs-new} finite, it is necessary to introduce some stability
condition on the objects that ${\mc M}_{d}$ classifies (see
\secref{def stack}). Imposing this stability condition should be
viewed as the geometric counterpart of Arthur's truncation process of
the trace formula (see \cite{CL}, where this is explained in the case
of the trace formula for Lie algebras), but we will not discuss this
issue in the present paper. Here we will consider the moduli stack
${\mc M}_{d}$ without any stability conditions.

\smallskip

In \secref{orbital} we will give a description of the stack ${\mc
  M}_d$ that is reminiscent and closely related to the Hitchin moduli
stack of Higgs bundles on the curve $X$ \cite{Hit1}. We will also
conjecture that $\Delta_{\Hc}^!({\mc K}_{d,\rho})$ is a pure perverse
sheaf on ${\mc M}_d$.

\medskip

If $X$ is a curve over $\C$, then the vector space \eqref{rhs-new}
still makes sense if we consider ${\mc K}$ as an object of either the
derived category of constructible sheaves or of ${\mc D}$-modules on
$\Bun_G$.

\section{Eigenvalues of $\K_{d,\rho}$ on Hecke eigensheaves}
\label{eigenvalues}

In this section we apply the functor $\K_{d,\rho}$ to Hecke
eigensheaves on $\Bun_G$.

\subsection{Cohomology of symmetric powers}

Recall the definition of the symmetric power of a local system ${\mc
L}$ on $X$ given in \secref{symmetric power}. We generalize this
definition slightly, by allowing ${\mc L}$ to be a {\em complex} of
local systems on $X$.\footnote{The reason for this will become clear
in \secref{constant sheaf}.} The resulting object ${\mc L}_d$ will
then be a complex of sheaves on $X_d$.

\medskip

Let us now compute the cohomology of $X_d$ with coefficients
in ${\mc L}_d$. By K\"unneth formula, we have
$$
H^\bullet(X_d,{\mc L}_d) = \left( H^\bullet(X^{d},{\mc L}^{\boxtimes
  d}) \right)^{S_d} = \left( H^\bullet(X,{\mc L})^{\otimes d}
\right)^{S_d},
$$
where the action of $S_d$ on the cohomology is as follows: it acts by
the ordinary transpositions on the even cohomology and by signed
transpositions on the odd cohomology. Thus, we find that
\begin{equation}    \label{sym power}
H^\bullet(X_d,{\mc L}_d) =
\bigoplus_{d_0+d_1+d_2=d}
S^{d_0}(H^0(X,{\mc L})) \otimes \Lambda^{d_1}(H^1(X,{\mc L})) \otimes
S^{d_2}(H^2(X,{\mc L})).
\end{equation}
This is true for any complex of local systems. The cohomological
grading is computed according to the rule that $d_0$ does not
contribute to cohomological degree, $d_1$ contributes $d_1$, and $d_2$
contributes $2d_2$. In addition, we have to take into account the
cohomological grading on ${\mc L}$.

\subsection{Connection to $L$-functions}    \label{L-fun}

Suppose that $k=\Fq$. Let us form a generating
function
\begin{equation}    \label{gen fun}
\sum_{d \geq 0} \on{Tr}(\on{Fr},H^\bullet(X_d,{\mc
  L}_d)) t^d,
\end{equation}
where $\on{Tr}(\on{Fr},H^\bullet(X_d,{\mc L}_d))$ is the trace of the
Frobenius on $H^\bullet(X_d,{\mc L}_d)$, thus the alternating sum of
traces on the groups $H^j(X_d,{\mc L}_d)$.

There are two ways to compute it.

\medskip

{\bf I.} By the Lefschetz formula,
$$
\on{Tr}(\on{Fr},H^\bullet(X_d,{\mc
  L}_d)) = \sum_{D \in X_d(\Fq)} \on{Tr}(\on{Fr}_D,{\mc
  L}_{D}),
$$
where ${\mc L}_D$ is the stalk of ${\mc L}_d$ at $D = \sum_i n_i
[x_i]$, where $\sum_i n_i \deg(x_i) = d$. Since
$$
{\mc L}_D = \bigotimes_i S^{n_i} {\mc L}_{x_i},
$$
we find that the generating function \eqref{gen fun} is equal to
\begin{equation}    \label{L-function}
\prod_{x \in |X|} \on{det}(1-t^{\deg(x)} \on{Fr}_x,{\mc L}_x)^{-1}.
\end{equation}
where $|X|$ is the set of closed points of $X$.

If we set $t=q^{-s}$ in formula \eqref{L-function}, we obtain the
$L$-function $L({\mc L},s)$ of the Galois representation associated to
${\mc L}$ (and defining representation $\on{def}$ of $GL_n$). Thus we
find that $\on{Tr}(\on{Fr},H^\bullet(X_d,{\mc L}_d))$ is the
coefficient of $q^{-ds}$ in $L({\mc L},s)$.

\medskip

{\bf II.}  Formula \eqref{sym power} gives us the following expression
for the generating function \eqref{gen fun}:
$$
\frac{\on{det}(1-t \on{Fr},H^1(X,{\mc L}))}{\on{det}(1-t
  \on{Fr},H^0(X,{\mc L})) \; \on{det}(1-t\on{Fr},H^2(X,{\mc L}))}.
$$
If we substitute $t=q^{-s}$, we obtain the Grothendieck--Lefschetz
formula for the $L$-function $L({\mc L},s)$.

\subsection{Hecke eigensheaves}    \label{h eig}

Recall that $X$ is a geometrically connected smooth proper curve over
a finite field $k=\Fq$, and $G$ a reductive group scheme over $X$
which is a quasi-split form of a constant group $\GG$. The Langlands
dual group of $G$, defined over $\Ql$, is $^L G=\check\GG \rtimes
\Gamma$ where $\Gamma$ is the (\'etale) fundamental group
$\Gamma=\pi_1(X,x)$ that acts on $\check\GG$ by the data defining the
quasi-split form $G$.

Let $W_F$ be the Weil group of the function field $F$ of the curve
$X$. An {\em Arthur parameter} is an equivalence class of
homomorphisms
\begin{equation}\label{Arthur parameter}
\sigma: \SL_2 \times W_F \to {}^L G
\end{equation}
that induces the canonical map $W_F \to \on{Gal}(\ol{F}/F) \to
\Gamma$. According to the conjectures of Langlands and Arthur,
equivalence classes of irreducible automorphic representations of
$G(\AD)$ may be parametrized by Arthur parameters. Such a
representation $\pi$ is called tempered (or Ramanujan) if
$\sigma|_{\SL_2}$ is trivial. Otherwise, it is called non-tempered (or
non-Ramanujan).

If an automorphic representation $$\pi = \bigotimes_{x \in |X|}{}'
\; \pi_x$$ of the ad\'elic group $G(\AD)$ has the Arthur parameter
$\sigma$, then for all closed points $x \in |X|$ where $\pi_x$ is
unramified the restriction $\sigma|_{W_F}$ is also unramified and the
Satake parameter of $\pi_x$ is equal to the conjugacy class
$$
\sigma\left(\begin{pmatrix} q^{1/2} & 0 \\ 0 & q^{-1/2} \end{pmatrix}
\times \on{Fr}_x \right) \in {}^L G.
$$

\medskip

Now we discuss a geometric analogue of this -- in the unramified case,
in order to simplify the discussion. The geometric analogue of the
notion of automorphic representations (or, more precisely, the
corresponding spherical automorphic function) is the notion of Hecke
eigensheaf which we now recall.

\medskip

Let $\rho$ be a continuous representation of $\LG$ on a
finite-dimensional vector space $V$ and let $\sigma$ be an Arthur
parameter. Then we obtain a representation $\rho \circ \sigma$ of
$\SL_2 \times W_F$ on $V$. The standard torus of $\SL_2 \subset \SL_2
\times W_F$ defines a $\Z$-grading on ${\rho\circ\sigma}$:
$$
{\rho\circ\sigma} = \bigoplus_{i\in \Z} (\rho\circ\sigma)_i,
$$
where each $\rho\circ\sigma_i$ is a continuous representation of
$W_F$. Assume that each of them is unramified. Then it gives rise to
an $\ell$-adic local system $\Lc_{(\rho\circ\sigma)_i}$ on $X$. Now we
define $\Lc_{\rho\circ\sigma}$ to be a complex of local systems on $X$
with the trivial differential
$$\Lc_{\rho\circ\sigma} = \bigoplus_{i\in \Z}
\Lc_{(\rho\circ\sigma)_i}[-i].$$

Recall the Hecke functors $\H_\rho$ and ${\mathbb H}_{\rho,x}$
introduced in \secref{perv sheaf} (geometric analogues of the Hecke
operators corresponding to $\rho$). The following definition is a
slight generalization (to the case when $\sigma|_{\SL_2}$ is
non-trivial) of the definition given in \cite{BD}.

\begin{definition}
A sheaf ${\mc F}$ on $\Bun_G$ is called a {\em Hecke eigensheaf with
the eigenvalue $\sigma$} \eqref{Arthur parameter} if for any
representation $\rho$ of $^L G$ we have an isomorphism
\begin{equation}    \label{iso}
{\mathbb H}_\rho({\mc F}) \simeq \Lc_{\rho \circ \sigma} \boxtimes
{\mc F},
\end{equation}
and these isomorphisms are compatible for different $\rho$ with
respect to the structures of tensor categories.
\end{definition}

\medskip

It follows from the above identity that for every $x\in X$, we have an
isomorphism
\begin{equation}    \label{isox}
{\mathbb H}_{\rho,x}({\mc F}) \simeq \Lc_{\rho \circ \sigma,x} \otimes
{\mc F},
\end{equation}
where $\Lc_{\rho \circ \sigma,x}$ is the stalk of $\rho \circ \sigma$
at $x$. Formula \eqref{isox} follows from (\ref{iso}) by base
change for the cohomology with compact support.

\medskip

In $X$ is a smooth projective connected curve over $\C$, then the
above definition also makes sense if we take as $\sigma$ a
homomorphism $\SL_2 \times \pi_1(X,x) \to {}^L G$ and as ${\mc
  L}_{(\rho \circ \sigma)_i}$ the corresponding local system on $X$.

\bigskip

Now we compute the ``eigenvalues'' of the integral transform ${\mathbb
  K}_{d,\rho}$ on ${\mc F}_\sigma$. The following result is Lemma 2.6
of \cite{FLN}.

\begin{lem} \label{Lemma 1 bis}
If ${\mc F} = {\mc F}_\sigma$ is a Hecke eigensheaf with eigenvalue
$\sigma$, then for every representation $\rho$ of $\LG$ and every
positive integer $d$ we have
\begin{equation} \label{d-th coeff}
{\mathbb K}_{d,\rho}({\mc
  F}_\sigma)=H^\bullet(X_d,\Lc_{\rho\circ\sigma,d}) \otimes {\mc
  F}_\sigma.
\end{equation}
\end{lem}

\begin{proof} We will restrict ourselves to the split case. The general
case is obtained from the split case by the descent method. From now
on, the group scheme $G$ is of the form $\GG\times X$, $\LG$ is
$\check\GG\times\Gamma$, the representation $\rho$ is assumed to be of
the form $\rho=\rho_\mu\otimes\rho_\Lc$ where $\rho_\mu$ is a
irreducible representation of $\check\GG$ of highest weight $\mu$ and
$\rho_\Lc$ is the continuous representation of $\Gamma$ attached to an
irreducible local system $\Lc$ over $X$.

Consider the chain version of the Hecke stack from
\secref{sym H} with truncation parameter $\mu$,
$$\Hc^d_\mu\to X^d \times\Bun_G\times\Bun_G,$$
with the (shifted) perverse sheaf
$$\wF^d_{\rho}=\IC(\Hc^d_\mu)\otimes \pr_{X^d}^*(\Lc^{\boxtimes d}).$$
As in \secref{int tr}, let us denote by $p_d, p'_d$ the two
projections to $\Bun_G$ and $q_d=\pr_{X^d} \times p_d$, where
$\pr_{X^d}$ is the projection to $X^d$. Consider the integral
transform
$$D(\Bun_G) \to D(X^d\times \Bun_G)$$
given by
$${\mathbb H}_\rho^d: {\mc F} \mapsto q_{d!}(p'_d{}^* {\mc F}\otimes
\wF^d_\rho).$$ By using repeatedly (\ref{iso}), we obtain
$${\mathbb H}_\rho^d({\mc F}_\sigma)=\Lc_{\rho\circ\sigma}^{\boxtimes
  d} \boxtimes {\mc F}_\sigma.$$
By pushing forward along $X^d$, we get
\begin{equation}    \label{pr bun}
\pr_{\Bun_G,!}{\mathbb H}_\rho^d({\mc
F}_\sigma)=H^\bullet(X,\Lc_{\rho\circ\sigma})^{\otimes d} \otimes {\mc
F}_\sigma.
\end{equation}
The symmetric group $S_d$ acts on the right hand side in the obvious
way and the invariant part is
$$H^\bullet(X_d,\Lc_{\rho\circ\sigma,d}) \otimes {\mc F}_\sigma.$$

On the other hand, recall the morphism
$$\pi_\Hc^d:\Hc^d\to\Hc_d.$$ By definition, the perverse sheaf
$\wF_{d,\rho}$ is the invariant part of $\pi_{\Hc!}^d \wF^d_\rho$
under the action of $S_d$. Hence the the $S_d$-invariant part of the
left hand side of \eqref{pr bun} is ${\mathbb K}_{d,\rho}({\mc
  F}_\sigma)$, and we obtain the desired formula (\ref{d-th coeff}).
\end{proof}

\begin{cor}    \label{van}
For all $d>(2g-2)\dim\rho$, we have ${\mathbb K}_{d,\rho}({\mc
F}_\sigma) \equiv 0$, unless $\rho \circ \sigma$ or $(\rho \circ
\sigma)^*$ has non-zero invariants under the action of $W_F$ (or
$\pi_1(X,x)$, if $X$ is defined over $\C$).
\end{cor}

\begin{proof} Indeed, if $\rho \circ \sigma$ and
$(\rho \circ \sigma)^*$ have zero spaces of invariants, then
  $$H^0(X,\rho \circ \sigma) = H^2(X,\rho \circ \sigma) = 0$$ and $\dim
  H^1(X,\rho \circ \sigma) = (2g-2)\dim\rho$ (since it is then equal
  to the Euler characteristic of the constant local system of rank
  $\dim\rho$ on $X$). Hence we obtain that
$$
H^\bullet(X_d,\Lc_{\rho\circ\sigma,d}) \simeq
\Lambda^d(H^1(X,\rho \circ \sigma)),
$$
and the statement of the corollary follows.
\end{proof}

We can translate the above lemma and corollary to the level of
functions: the eigenvalue of the operator ${\mb K}_{d,\rho}$ (see
\secref{int tr}) on an unramified Hecke eigenfunction $f_\sigma$ is
equal to the $q^{-ds}$-coefficient of the automorphic $L$-function
$L(\rho,\sigma,s)$. This can indeed be shown directly at the level of
functions as follows (see also the proof of Lemma 2.6 in \cite{FLN}).

When ${\mathbb K}_{d,\rho}$ acts on ${\mc F}$, for a fixed effective
divisor $D = \sum_i n_i [x_i]$ of degree $d$ we apply modifications at
the closed points $x_i$. These modifications are described by the
Schubert varieties $\ol{\on{Gr}}_{n_i \mu}$ with the ``kernel''
$S^{n_i}(\IC(\ol\Gr_\mu))$, which is the perverse sheaf on
$\ol{\on{Gr}}_{n_i \mu}$ corresponding to the $n_i$th symmetric power
of the representation $\rho$.

When we apply the Hecke functor corresponding to
$S^{n_i}(\IC(\ol\Gr_\mu))$ at the point $x_i$ to our sheaf ${\mc F}$,
we simply tensor it by $S^{n_i}(\rho \circ \sigma)_{x_i}$, because
${\mc F}$ is a Hecke eigensheaf with eigenvalue $\sigma$. So the net
result is that for each fixed $D$ we tensor ${\mc F}$ with
$$
\bigotimes_i S^{n_i}(\rho \circ \sigma)_{x_i}
$$
(note that in general this is a complex of vector
spaces). But this is precisely the stalk of $(\rho \circ
\sigma)_d$ at $D$ (see formula \eqref{stalk}). In other words,
because we are applying ${\mathbb K}_{d,\rho}$ to a Hecke
eigensheaf ${\mc F}$, we are effectively replacing each Hecke
operator by the corresponding ``eigenvalue'', that is,
$S^{n_i}(\rho \circ \sigma)_{x_i}$.

Then we need to integrate over all possible $D$. This simply means
taking the cohomology of $X_d$ with coefficients in $(\rho \circ
\sigma)_d$. The result is
$$
{\mathbb K}_{d,\rho}({\mc F}) \simeq
H^\bullet(X_d,(\rho \circ \sigma)_d) \otimes {\mc F},
$$
as in \lemref{Lemma 1 bis}. Now, each step of this calculation makes
sense at the level of functions. Hence we obtain that the eigenvalue
of the operator ${\mb K}_{d,\rho}$ on $f_\sigma$ is equal to the trace of
the Frobenius on $H^\bullet(X_d,(\rho \circ \sigma)_d)$,
which is the $q^{-ds}$-coefficient of the automorphic $L$-function
$L(\rho,\sigma,s)$, as we have seen in \secref{L-fun}.

If $\rho \circ \sigma$ and $(\rho \circ \sigma)^*$ have zero spaces of
invariants, then $L(\sigma,\on{def},s)$ is a polynomial of degree
$(2g-2)\dim\rho$ in $q^{-s}$, and so ${\mb K}_{d,\rho}(f_\sigma) = 0$
in this case for all $d>(2g-2)\dim\rho$. This is function-theoretic
version of Corollary \ref{van}.

\subsection{Example: constant sheaf}    \label{constant sheaf}

Assume that $G$ is split, $\LG=\check\GG\times\Gamma$,
$\rho=\rho_\mu\otimes \rho_0$, where $\rho_\mu$ is the irreducible
representation of highest weight $\mu$, $\rho_0$ is the trivial
representation of $\Gamma$. In this case $\wF_\rho=\wF_\mu$ is the
intersection cohomology complex of $\Hc_\mu$ shifted by
$-\dim(X\times\Bun_G)$.

Let ${\mc F}_0 = \underline\Ql|\Bun_G$, the constant sheaf on $\Bun_G$.
This is the geometric analogue of the trivial representation of
$G({\mathbb A})$.  Let us apply the Hecke operator ${\mathbb
H}_{\mu,x}$ to ${\mc F}_0$. For any $G$-principal bundle $E$, the fiber
of $p^{-1}(E)\cap \Hc_x$ is isomorphic to $\Gr_x$, once we have
chosen a trivialization of $E$ on the formal disc $D_x$. Thus the
fiber of ${\mathbb H}_{\rho,x}({\mc F}_0)$ at $E$ is isomorphic to
$${\mathbb H}_{\rho,x}({\mc
F})_E=H^\bullet(\ol{\on{Gr}}_\mu,\IC(\ol\Gr_\mu)).$$ This isomorphism is
actually canonical in the sense that it does not depend on the choice
of the trivialization of $E$ on $D_x$, so that the above isomorphism
can be put in a family with respect to the parameter $E$. Hence we
obtain that
$${\mathbb H}_{\rho,x}({\mc F}_0) \simeq
H^\bullet(\ol{\on{Gr}}_\mu,\IC(\ol\Gr_\mu)) \otimes {\mc F}_0.$$

By the geometric Satake correspondence \cite{MV},
$$
H^\bullet(\ol{\on{Gr}}_\mu,\IC(\ol\Gr_\mu)) \simeq \rho_\mu^{\on{gr}},
$$
a complex of vector spaces, which is isomorphic to the representation
$\rho_\mu$ with the cohomological grading corresponding to the
principal grading on $\rho_\mu$. In other words, this is $(\rho \circ
\sigma_0)_x$, where $\sigma_0: W_F \times \SL_2 \to {}^L G$ is
trivial on $W_F$ and is the principal embedding on $\SL_2$. We conclude
that the constant sheaf on $\Bun_G$ is a Hecke eigensheaf with the
eigenvalue $\sigma_0$.  This is in agreement with the fact that
$\sigma_0$ is the Arthur parameter of the trivial automorphic
representation of $G({\mathbb A})$.

For example, if $\rho_\mu$ is the defining representation of $GL_n$,
then the corresponding Schubert variety is ${\mathbb P}^{n-1}$, and we
obtain its cohomology shifted by $(n-1)/2$, because the intersection
cohomology sheaf $\IC(\ol\Gr_\mu)$ is the constant sheaf placed in
cohomological degree $-(n-1)$, that is
$$
H^\bullet(\ol{\on{Gr}}_\mu,\IC(\ol\Gr_\mu)) = \theta^{(n-1)/2} \oplus
\theta^{(n-3)/2} \oplus \ldots \oplus \theta^{-(n-1)/2},
$$
where $\theta^{1/2} = \Ql[-1](-1/2)$. This agrees with the fact that
the principal grading takes values $(n-1)/2,\ldots,-(n-1)/2$ on the
defining representation of $GL_n$, and each of the corresponding
homogeneous components is one-dimensional.

Next, we find that
$$
{\mathbb K}_{d,\rho_\mu}(\underline\Ql) \simeq
H^\bullet(X_d,(\rho_\mu \circ \sigma_0)_d) \otimes
\underline\Ql,
$$
where $\rho_\mu \circ \sigma_0$ is the complex described in the
previous section. At the level of functions, we are multiplying the
constant function by
\begin{equation}    \label{prod of zeta1}
\on{Tr}(\on{Fr},H^\bullet(X_d,(\rho_\mu \circ \sigma_0)_d)) =
\prod_{i \in P(\rho_\mu)} \zeta(s-i)^{\dim \rho_{\mu,i}},
\end{equation}
where $P(\rho_\mu)$ is the set of possible values of the principal
grading on $\rho_\mu$ and $\rho_{\mu,i}$ is the corresponding subspace
of $\rho_\mu$.

For example, if $\rho_\mu$ is the defining representation of $GL_n$,
then formula \eqref{prod of zeta1} reads
\begin{equation}    \label{prod of zeta}
\prod_{k=0}^{n-1} \zeta(s+k-(n-1)/2).
\end{equation}

\subsection{Decomposition of the trace}    \label{decomp trace}

Now we wish to use Corollary \ref{van} to decompose the trace of ${\mb
  K}_{d,\rho}$ for $d>(2g-2)\dim \rho$ as a sum over subgroups of $^L
G$.

Let $\sigma$ be an Arthur parameter. We attach to it two subgroups of
$^L G$: $^\la G = {}^\la G_\sigma$ is the centralizer of the image of
$\SL_2$ in $^L G$ under $\sigma$, and $^\la H = {}^\la H_\sigma$ is
the Zariski closure of the image of $W_F$ in $^\la G_\sigma$ under
$\sigma$. The representation $\rho \circ \sigma$ has non-zero
invariants if and only if the restriction of $\rho$ to $^\la H_\sigma$
has non-zero invariants. Thus, assuming that Arthur's conjectures are
true, we obtain that the trace of ${\mb K}_{d,\rho}$ decomposes as a
double sum: first, over different homomorphisms $\phi: \SL_2 \to {}^L
G$, and second, for a given $\phi$, over the subgroups $^\la H$ of the
centralizer $^\la G_\phi$ of $\phi$ having non-zero invariants in
$\rho$:
\begin{equation}    \label{Phi and Psi}
\on{Tr} {\mb K}_{d,\rho} = \sum_\phi \sum_{^\la H \subset {}^\la G_\phi}
\Phi_{\phi,{}^\la H}.
\end{equation}
The summands $\Phi_{\on{triv},{}^\la H}$ correspond to the Ramanujan
representations. Denote their sum by $(\on{Tr} {\mb K}_{d,\rho})_R$.

Note that if $\phi$ is non-trivial, then the rank of $^\la G_\phi$ is
less than that of $^L G$. As explained in \cite{FLN}, we would like to
use induction on the rank of $^L G$ to isolate the Ramanujan part
$(\on{Tr} {\mb K}_{d,\rho})_R$ in $\on{Tr} {\mb K}_{d,\rho}$. In other
words, we need to isolate and remove the contribution of the
non-Ramanujan representations. In \cite{FLN} it was shown how to
isolate the contribution of the trivial representation.

We then wish to decompose $(\on{Tr} {\mb K}_{d,\rho})_R$ over $^\la
H$,
\begin{equation}    \label{sum2}
(\on{Tr} {\mb K}_{d,\rho})_R = \sum_{^\la H \subset {}^L G} (\on{Tr}
  K^H_{d,\rho_H})_R.
\end{equation}
Here the sum should be over all possible $^\lambda H \subset {}^L G$
such that $^\lambda H$ has non-zero invariant vectors in $\rho \circ
\sigma$, and $K^H_{d,\rho_H}$ is the operator corresponding to $\rho_H
= \rho|_{^\la H}$ for the group $H(\AD)$. Proving formula \eqref{sum2}
is the main step in the strategy to prove functoriality outlined in
\cite{FLN}. More specifically, we would like to use the orbital side
of the trace formula \eqref{tr} to establish \eqref{sum2}.

\medskip

Actually, comparisons of trace formulas should always be understood as
comparisons of their {\em stabilized} versions. Therefore the traces
in \eqref{sum2} should be replaced by the corresponding stable traces
(see \cite{FLN} for more details). Since the fundamental lemma has
been proved \cite{Ngo:FL}, the connection between the actual trace
formula and the stabilized trace formula is now well-understood.

\medskip

The spectral side of \eqref{sum2} (again, assuming Arthur's
conjectures) is the sum of the eigenvalues of ${\mb K}_{d,\rho}$ which
are expressed in terms of the coefficients of the $L$-function of the
corresponding Arthur parameter $\sigma: W_F \to {}^\la H$, for
different $^\la H$. A precise formula for these eigenvalues is
complicated in general, but we can compute its asymptotics as $d \to
\infty$. If we divide ${\mb K}_{d,\rho}$ by $q^d$ (geometrically, this
makes sense because of the dimension of $X_d$), then the asymptotics
will be very simple:
\begin{equation}    \label{asymptotics}
q^{-d} (\on{Tr} {\mb K}_{d,\rho})_R \sim \sum_{^\lambda H \subset {}^L G} \quad
\sum_{\sigma': W_F \to {}^\lambda H} N_\sigma \begin{pmatrix}
d+m_\sigma(\rho)-1 \\ m_\sigma(\rho)-1 \end{pmatrix},
\end{equation}
where $N_\sigma$ is the multiplicity of automorphic representations in
the corresponding $L$-packet.

Indeed, the highest power of $q$ comes from the highest cohomology,
which in this case is
$$
H^{2d}(X_d,(\rho \circ \sigma)_d) = \on{Sym}^{d}(H^2(X,\rho
\circ \sigma))
$$
($d_0=0,d_1=0$ and $d_2=d$ in the notation of formula \eqref{sym
power}). We have $\dim H^2(X,\rho \circ \sigma) = m_\sigma(\rho)$, the
multiplicity of the trivial representation in $\rho \circ
\sigma$ (we are assuming here that this trivial representation
splits off as a direct summand in $\rho \circ \sigma$), and
\begin{equation}    \label{dim}
\dim \on{Sym}^{d}(H^2(X,\rho \circ \sigma)) = \begin{pmatrix}
  d+m_\sigma(\rho)-1 \\ m_\sigma(\rho)-1 \end{pmatrix}.
\end{equation}

Thus, as a function of $q^d$, the eigenvalues of $q^{-d} {\mb
  K}_{d,\rho}$ on the Ramanujan representations grow as $O(1)$ when $d
\to \infty$. For the non-Ramanujan representations corresponding to
non-trivial $\phi: \SL_2 \to {}^L G$, they grow as a higher power of
$q^d$. For instance, the eigenvalue corresponding to the trivial
representation of $G(\AD)$ (for which $\phi$ is a principal embedding)
grows as $O(q^{d(\rho,\mu)})$, where $\mu$ is the highest weight of
$\rho$. In general, it grows as $O(q^{da})$, where $2a$ is the maximal
possible highest weight of the image of $\SL_2 \subset {}^L G$ under
$\phi$ acting on $\rho$. Thus, the asymptotics of the non-Ramanujan
representations dominates that of Ramanujan representations. This is
why we wish to isolate the Ramanujan part first and then decompose it
over $^\la H$.

\bigskip

The main question we want to answer is the following: how to observe
the decompositions \eqref{Phi and Psi} and \eqref{sum2} on the orbital
side of the trace formula, using geometry?

\subsection{Example of $\SL_2$ and the twisted torus
  continued}    \label{ex sl2}

Let us go back to the case that $G = X \times \SL_2$ and the Langlands
dual group is $\LG=\PGL_2\times\Gamma$. We consider the representation
$\rho=\rho_\mu\otimes\rho_\Ec$ of $^L G$, where $\rho_\mu$ is the
three-dimensional adjoint representation of $\PGL_2$ and $\Ec$ is a
non-trivial rank one local system over $X$ of order two which is non
trivial. We are now looking at the effect of the operator ${\mathbb
  K}_{d,\rho}$ on a Hecke eigenfunction $f_\sigma$ attached to an
Arthur parameter $\sigma:\SL_2\times \Gamma\to \LG$.

Consider first the case of the Arthur parameter $\sigma=\sigma_0$
whose restriction to the factor $\SL_2$ is non-trivial, and hence an
isomorphism. It then follows that the restriction of $\sigma$ to $W_F$
must be trivial. The automorphic representation of $\SL_2(\AD)$
attached to this parameter is the trivial one-dimensional
representation. In this case ${\mc F}_{\rho\circ\sigma}$ is the graded
local system
$${\mc F}_{\rho\circ\sigma_0} = (\und\Ql[2]\oplus \und\Ql \oplus
\und\Ql[-2]) \otimes \Ec,$$ that is, the direct sum of three copies of
$\Ec$ put in cohomological degrees $-2,0$ and $2$. Since
$H^i(X,\Ec)=0$ for $i\not=1$, we have
$$H^\bullet(X,{\mc F}_{\rho\circ\sigma_0})=H^1(X,\Ec)\otimes (\Ql[1] \oplus
\Ql[-1] \oplus \Ql[-3]),$$
the direct sum of three copies of $H^1(X,\Ec)$ put in degrees $-1,1$
and $3$. In particular $H^\bullet(X,{\mc F}_{\rho\circ\sigma_0})$ vanishes in
even degrees. It follows from an obvious generalization of formula
(\ref{sym power}) that for a large integer $d$, we have
$$H^\bullet(X_d,{\mc F}_{d,\rho\circ\sigma_0})=0.$$ Hence we obtain that
the integral transform $\K_{d,\rho}$ annihilates the contribution of
the trivial representation of $\SL_2(\AD)$, which is the only
non-Ramanujan representation in this case. Thus, only Ramanujan
representations of $\SL_2(\AD)$ contribute to the trace formula for
large $d$ in this case.

\medskip

Now we look at the Arthur's parameters $\sigma$ whose restriction to
$\SL_2$ factor is trivial. The restriction of $\sigma$ to $W_F$ is of
the form $$\sigma(\alpha)=(\sigma_+(\alpha),\varphi(\alpha))\in
\check\GG\times \Gamma, \qquad \alpha\in W_F,
$$
where $\varphi: W_F \to \Gamma$ is the canonical homomorphism, for
some homomorphism $\sigma_+: W_F \to \check\GG$. The local system
$\Lc_{\rho\circ\sigma}$ is of the form
$$\Lc_{\rho\circ\sigma}=\Lc_{\rho_\mu\circ\sigma_+} \otimes \Ec,$$
where $\Lc_{\rho_\mu\circ\sigma_+}$ is the rank three local system
given by the representation of $W_F$ given by the composition of
$\sigma_+:W_F\to \PGL_2$ and the three-dimensional representation
$\rho_\mu$ of $\PGL_2$.

If either $H^0(X,\Lc_{\rho\circ\sigma})$ or
$H^2(X,\Lc_{\rho\circ\sigma})$ is non-zero, then the three-dimensional
representation $ \rho\circ\sigma$ must have either non-zero space of
invariants or coinvariants. Since $\Lc_{\rho_\mu\circ\sigma}$ is
semi-simple because of the purity, $\rho\circ\sigma$ has non-zero
invariants. Let $^\la H$ be the closure of the image of
$\sigma_+$. We want to understand when this subgroup has an invariant
vector $v$ in the adjoint representation $V$ of $\PGL_2$.

We identify $\PGL_2$ with the group ${\rm SO}_3$ preserving the
Killing form on $V$. There are two possibilities for the invariant
vector: be isotropic or anisotropic.

If the invariant vector is $v$ is isotropic, $v$ belongs to the
two-dimensional vector space $\langle v\rangle^\bot$. Since we have
assumed $\sigma_+$ to be semi-simple, there exists another vector
$v_1\notin \langle v\rangle^\bot$ such that $\sigma_+(W_F)$ preserves
the line $\langle v_1 \rangle$. Let $v_{2}$ be a vector generating the
one-dimensional vector space $\langle v\rangle^\bot \cap \langle
v_1\rangle^\bot$.  We have a decomposition of $V$ into lines preserved
by $\sigma_+(W_F)$
$$V=\langle v\rangle \oplus \langle v_1\rangle \oplus \langle v_2
\rangle.$$ One can check that $\sigma_+(W_F)$ is contained in the
split torus ${\rm SO}_2 \subset {\rm SO}_3$. This means that $\sigma$
is a parameter of an Eisenstein series, so that the rank three local
system $\Lc_{\rho_\mu\circ\sigma_+}$ is of the form
$$\Lc_{\rho_\mu\circ\sigma_+}=\Ec_1\oplus \Ec_1^{-1} \oplus \und\Ql$$
for some rank one local system $\Ec_1$ on $X$.

If $\Lc_{\rho_\mu\circ\sigma}=\Lc_{\rho_\mu\circ\sigma_+} \otimes \Ec$
contains the trivial local system, the only possibility is
$\Ec=\Ec_1$, so that
$$\Lc_{\rho_\mu\circ\sigma}=\und\Ql \oplus \und\Ql \oplus\Ec.$$

If the invariant vector $v$ is anisotropic, then $v\notin \langle
v\rangle^\bot$ and we have an orthogonal decomposition of $V$
$$V=\langle v\rangle \oplus \langle v\rangle^\bot.$$ In this case,
$\sigma_+(W_F)$ is contained in the subgroup ${\rm O}_2$ of ${\rm
SO}_3$,
$${\rm O}_2=({\rm O}(1)\times {\rm O}_2) \cap {\rm SO}_3.$$ If
$\sigma_+(W_F)$ is contained in ${\rm SO}_2$, we end up in the
previous case. If it is not contained in ${\rm SO}_2$, then $\sigma_+$
induces a surjective homomorphism $\rho_{\Ec_1}:W_F \to \Z/2\Z$
corresponding to a certain non-trivial rank one local system $\Ec_1$
of order two. In this case $\Lc_{\rho_\mu\circ\sigma_+}$ has the form
$$\Lc_{\rho_\mu\circ\sigma_+}=\Ec_1\oplus {\mc V}_1,$$ where $V_1$ is
a rank two local system of determinant $\det({\mc V}_1)=\Ec_1$. If the
image of $\sigma_+$ contains only two elements, then we again end up
in the same situation as above, i.e., ${\mc V}_1=\und\Ql \oplus \Ec_1$
and $\sigma$ is the Arthur parameter of an Eisenstein series that we
have encountered above.

Finally, if the image of $\sigma_+$ contains more than two elements,
then ${\mc V}_1$ is an irreducible rank two local system. Now if
$$\Lc_{\rho_\mu\circ\sigma}=(\Ec_1\otimes \Ec) \oplus ({\mc
V}_1\otimes\Ec)$$ contains the trivial local system, the only
possibility is $\Ec_1=\Ec$, because ${\mc V}_1\otimes \Ec$ is an
irreducible rank two local system. This means that $\sigma$ factors
through the dual group $^L H$ of the twisted one-dimensional torus
$H=H_{\Ec}$ attached to $\Ec$.

To summarize, apart from the Eisenstein series, only (cuspidal)
automorphic representations with the Arthur parameter factorizing
through the dual group $^L H$ of the twisted one-dimensional torus
$H_{\Lc}$ contribute to the trace of the operator ${\mb
  K}_{d,\rho}$. Actually, the parameter of the Eisenstein series which
contributes to the trace also factors through $^L H$ as follows:
$$\sigma: W_F \to \Gm \rtimes W_F \to {\rm PGL}_2\times \Gamma,$$
where $W_F \to \Gm \rtimes \Gamma$ is given by $\alpha\to
\rho_{\Ec}(\alpha)\rtimes \varphi(\alpha)$. The contribution of these
representations must be exactly the same as the trace of the operator
${\mathbf K}_{d,\rho_H}$ for the twisted torus $H$ corresponding to
the restriction $\rho_H$ of $\rho$ from $\LG$ to $^L
H$.

\section{The moduli stack of $G$-pairs}    \label{orbital}

As we have seen in \secref{v sp}, the right hand side of the trace
formula \eqref{tr} is equal (up to a power of $q$) to the trace of the
(arithmetic) Frobenius on the vector space
\begin{equation}    \label{restr}
H^\bullet({\mc M}_{d},\Delta_{\Hc}^!({\mc K}_{d,\rho})).
\end{equation}
Recall from \secref{symmetric power} that the sheaf ${\mc K}_{d,\rho}$
is supported on the substack $\Hc_{d,[\mu]}$ of $\Hc$. In this section
we will mostly consider the split case, so in order to simplify our
notation we will write $\Hc_{d,\mu}$ for $\Hc_{d,[\mu]}$. Let ${\mc
M}_{d,\mu}$ be the fiber product of $X_d \times \Bun_G$ and ${\mc
H}_{d,\mu}$ with respect to the two morphisms to $X_d \times
\Bun_G \times \Bun_G$. In other words, we replace $\Hc_d$ by
$\Hc_{d,\mu}$ in the upper right corner of the diagram
\eqref{Cartesian}. The sheaf $\Delta_{\Hc}^!({\mc K}_{d,\rho})$ is
supported on ${\mc M}_{d,\mu} \subset {\mc M}_d$, and hence the vector
space \eqref{restr} is equal to
\begin{equation}    \label{restr1}
H^\bullet({\mc M}_{d,\mu},\Delta_{\Hc}^!({\mc K}_{d,\rho})).
\end{equation}

In this section we show that the stack ${\mc M}_{d,\mu}$ has a
different interpretation as a moduli space of objects that are closely
related to {\em Higgs bundles}. More precisely, we need ``group-like''
versions of Higgs bundles (we will call them ``$G$-pairs''). The
moduli spaces of (stable) Higgs bundles has been introduced by Hitchin
\cite{Hit1} (in characteristic $0$) and the corresponding stack (in
characteristic $p$) has been used in \cite{Ngo:FL} in the proof of the
fundamental lemma. In addition, there is an analogue ${\mc A}_{d,\mu}$
of the Hitchin base and a morphism $h_{d,\mu}: {\mc M}_{d,\mu} \to
{\mc A}_{d,\mu}$ analogous to the Hitchin map. We would like to use
this Higgs bundle-like realization of ${\mc M}_{d,\mu}$ and the
morphism $h_{d,\mu}$ in order to derive the decompositions \eqref{Phi
  and Psi} and \eqref{sum2}.

\subsection{Definition of the moduli stack}    \label{def stack}

Let us assume that $G$ is split over $X$ and $\mu$ is a fixed dominant
coweight.  The groupoid ${\mc M}_{d,\mu}(k)$ classifies the triples
$$
(D,E,\varphi),
$$
where $D = \sum_i n_i [x_i]\in X_d$ is an effective divisor of degree $d$,
$E$ is a principal $G$-bundle on a curve $X$, and $\varphi$ is a section
of the adjoint group bundle
$$
\on{Ad}(E) = E \underset{G}\times G
$$
(with $G$ acting on the right $G$ by the adjoint action) on $X {-}
\on{supp}(D)$, which satisfies the local conditions
\begin{equation}\label{invariant 4}
{\rm inv}_{x_i}(\varphi) \leq n_i\mu
\end{equation}
at $D$. Since we have defined $\Hc_{d,\mu}$ as the image of
$\Hc^d_\mu$ in $\Hc^\BD_d$, it is not immediately clear how to make
sense of these local conditions over an arbitrary base (instead of
$\on{Spec}(k)$). There is in fact a functorial description of
$\Hc_{d,\mu}$ and of ${\mc M}_{d,\mu}$ that we will now explain.

We will assume that $G$ is semi-simple and simply-connected. The
general case is not much more difficult.  Let
$\omega_1,\ldots,\omega_r$ denote the fundamental weights of $G$ and
$$\rho_{\omega_i}:G\to \GL(V_{\omega_i})$$ the Weyl modules of highest
weight $\omega_i$. Using the natural action of $G$ on ${\rm End}
(V_{\omega_i})$, we can attach to any $G$-principal bundle $E$ on $X$
the vector bundle
\begin{equation}
\on{End}_{\omega_i}(E)=E \underset{G}\times {\rm End}(V_{\omega_i}).
\end{equation}
The section $\varphi$ of $\on{Ad}(E)$ on $X{-}\on{supp}(D)$ induces a
section $\on{End}_{\omega_i}(\varphi)$ of the vector bundle
$\on{End}_{\omega_i}(E)$ on $X{-}\on{supp}(D)$. The local conditions
(\ref{invariant 4}) are equivalent to the property that for all $i$,
$\on{End}_{\omega_i}(\varphi)$ may be extended to a section
$$\varphi_i\in \on{End}_{\omega_i}(E)\otimes_{\Oc_X} \Oc_X(\langle
\mu,\omega_i\rangle D).$$ Though
the $\varphi_i$ determine $\varphi$, we will keep $\varphi$ in the
notation for convenience.

Thus, we obtain a provisional functorial description of ${\mc
M}_{d,\mu}$ as the stack classifying the data
\begin{equation}\label{phi phi}
(D,E,\varphi,\varphi_i)
\end{equation}
with $D\in X_d$, $E\in \Bun_G$, $\varphi$ is a section of $\on{Ad}(E)$
on $X{-}\on{supp}(D)$, $\varphi_i$ are sections of
$\on{End}_{\omega_i}(E)\otimes_{\Oc_X} \Oc_X(\langle
\mu,\omega_i\rangle D)$ over $X$ such that
$$\varphi_i|_{X{-}\on{supp}(D)}=\on{End}_{\omega_i}(\varphi).$$

Sometimes it will be more convenient to package the data
$(\varphi,\varphi_i)$ as a single object $\wt\varphi$ which has values
in the closure of $$(t_i \rho_{\omega_i}(g))_{i=1}^r \subset
\prod_{i=1}^r \on{End}(V_{\omega_i}).$$ where $g\in G$ and
$t_1,\ldots,t_r \in \Gm$ are invertible scalars. This way Vinberg's
semi-group \cite{Vinberg} makes its appearance naturally in the
description of ${\mc M}_{d,\mu}$. We will discuss this in more detail
elsewhere, and for the time being will stick to the more concrete
description (\ref{phi phi}) of the moduli space.

\subsection{Comparison with the Hitchin fibration} 

It is instructive to note that the stack ${\mc M}_{d,\mu}$ is very
similar to the moduli stacks of Higgs bundles (defined originally by
Hitchin \cite{Hit1} and considered, in particular, in
\cite{Ngo:Endo,Ngo:FL}). The latter stack ${\mc N}_D$ also depends on
the choice of an effective divisor
$$
D = \sum_i n_i [x_i]
$$
on $X$ and classifies pairs $(E,\phi)$, where $E$ is again a
$G$-principal bundle on $X$ and $\phi$ is a section of the adjoint
{\em vector bundle}
$$
\on{ad}(E) = E \underset{G}\times \g
$$
(here $\g = \on{Lie}(G)$) defined on $X{-}\on{supp}(D)$, which is
allowed to have a pole of order at most $n_i$ at $x_i$. In other
words,
$$
\phi \in H^0(X,\on{ad}(E) \otimes {\mc O}_X(D)).
$$
This $\phi$ is usually referred to as a {\em Higgs field}.

In both cases, we have a section which is regular almost everywhere,
but at some (fixed, for now) points of the curve these sections are
allowed to have singularities which are controlled by a divisor. In
the first case we have a section $\varphi$ of adjoint group bundle
$\on{Ad}(E)$, and the divisor is $D \cdot \mu$, considered as an
effective divisor with values in the lattice of integral weights of
$^L G$. In the second case we have a section $\phi$ of the adjoint Lie
algebra bundle $\on{ad}(E)$, and the divisor is just the ordinary
effective divisor.

An important tool in the study of the moduli stack ${\mc N}_D$ is the
{\em Hitchin map} \cite{Hi2} from ${\mc N}_D$ to an affine space
$$
{\mc A}_D  \simeq \bigoplus_i
H^0(X,{\mc O}_X((m_i+1)D),
$$
where the $m_i$'s are the exponents of $G$. It is obtained by, roughly
speaking, picking the coefficients of the characteristic polynomial of
the Higgs field $\phi$ (this is exactly so in the case of $GL_n$; but
one constructs an obvious analogue of this morphism for a general
reductive group $G$, using invariant polynomials on its Lie
algebra). A point $a \in A_E$ then records a stable conjugacy class in
$\g(F)$, where $F$ is the function field, and the number of points in
the fiber over $a$ is related to the corresponding orbital integrals
in the Lie algebra setting (see \cite{Ngo:Endo,Ngo:FL}).

More precisely, ${\mc A}_D$ is the space of section of the bundle
$${\mathfrak t}/W \underset{\Gm}\times {\mc O}_X(D)^\times$$ obtained by
twisting ${\mathfrak t}/W = \on{Spec}(k[{\mathfrak t}]^W)$, equipped with
the $\Gm$-action inherited from ${\mathfrak t}$, by the $\Gm$-torsor
${\mc O}_X(D)^\times$ on $X$ attached to the line bundle ${\mc
O}_X(D)$. Recall that $k[\mathfrak t/W]$ is a polynomial algebra with
homogeneous generators of degrees $d_1+1,\ldots,d_r+1$.

In our present setting, we will have to replace ${\mathfrak t}/W$ by
$T/W$.  Recall that we are under the assumption that $G$ is
semi-simple and simply-connected. First, recall the isomorphism of
algebras
$$k[G]^G=k[T]^W=k[T/W].$$ It then follows from \cite{Bourbaki},
Th. VI.3.1 and Ex. 1, that $k[G]^G$ is a polynomial
algebra generated by the functions
$$g\mapsto \on{tr}(\rho_{\omega_i}(g)),$$ where $\omega_1,\ldots,
\omega_r$ are the fundamental weights of $G$.

For a fixed divisor $D$, the analogue of the Hitchin map for ${\mc
M}_{d,\mu}(D)$ (the fiber of ${\mc M}_{d,\mu}$ over $D$) is the
following map:
\begin{equation}
{\mc M}_{d,\mu}(D) \to \bigoplus_{i=1}^r H^0(X,\Oc_X(\langle
\mu,\omega_i \rangle D ))
\end{equation}
defined by attaching to $(D,E,\varphi,\varphi_i)$ the collection of traces
$${\rm tr}(\varphi_i)\in H^0(X,\Oc_X(\langle \mu,\omega_i \rangle D )).$$
By letting $D$ vary in $X_d$, we obtain a fibration
$$h_{d,\mu}:{\mc M}_{d,\mu} \to {\mc A}_{d,\mu}$$ where ${\mc
A}_{d,\mu}$ is a vector bundle over $X_d$ with the fiber
$\bigoplus_{i=1}^r H^0(X,\Oc_X(\langle \mu,\omega_i \rangle D ))$ over
an effective divisor $D\in X_d$.

The morphism $h_{d,\mu}:\M_{d,\mu} \to \Ac_{d,\mu}$ is very similar to
the Hitchin fibration. We are now going to outline the geometric
properties of the Hitchin fibration that can be carried over to our
new situation.  We will discuss in more detail the case of $\SL_2$ and
postpone the general case for another occasion.

There exists a Picard stack
${\mc P}_{d,\mu} \to {\mc A}_{d,\mu}$ that plays the same role as 
${\mc P} \to {\mc A}$ constructed in \cite{Ngo:Endo}. In particular, 
for every $a\in {\mc A}_{d,\mu}$, the Picard stack ${\mc P}_{a}$ acts
on the fiber ${\mc M}_{a}$ of $h_{d,\mu}$ over $a$ containing a 
principal homogeneous space as an open substack. For generic $a$,
${\mc P}_{a}$ acts simply transitively on ${\mc M}_{a}$. In general,
there exists a product formula that expresses the quotient $[{\mc M}_{a}/ 
{\mc P}_{a}]$ in terms of local data as in \cite{Ngo:Endo,Ngo:FL}.

Assume from now on that $G$ is semisimple. As in \cite{CL}, there
exists a open substack $\M_{d,\mu}^{\rm st}$ of $\M_{d,\mu}$ that is
proper over $\Ac_{d,\mu}$.  This open substack depends on the choice
of a stability condition. However, its cohomology should be
independent of this choice. Moreover, there exists an open subset
$\Ac_{d,\mu}^{\rm ani}$ of $\Ac_{d,\mu}$ whose $\bar k$-points are the
pairs $(D,b)$ such that as an element of $(T/W)(F\otimes_{k} \bar k)$,
$b$ corresponds to a regular semisimple and anisotropic conjugacy
class in $G(F\otimes_{k} \bar k)$. The preimage $\M_{d,\mu}^{\rm ani}$
of $\Ac_{d,\mu}^{\rm ani}$ is contained in $\M_{d,\mu}^{\rm st}$ for
all stability conditions. In particular, the morphism $\M_{d,\mu}^{\rm
  ani}\to \Ac_{d,\mu}^{\rm ani}$ is proper.

Let ${\mc P}_{d,\mu}^{\rm tor}$ be the open substack of ${\mc
  P}_{d,\mu}$ such that for every $a\in \Ac_{d,\mu}$, the fiber ${\mc
  P}_{a}^{\rm tor}$ is an open subgroup of ${\mc P}_{a}$ and the
component group $\pi_{0}({\mc P}_{a}^{\rm tor})$ is the torsion part
of $\pi_{0}(\Pc_{a})$. Whatever stability condition we choose,
${\mc P}_{d,\mu}^{\rm tor}$ acts on $\M_{d,\mu}^{\rm st}$.

Our goal is to understand the cohomology \eqref{restr1}. According to
formula \eqref{verdier}, up to a shift and Tate twist (and replacing
the local system $\Lc$ by $\Lc^*$), it is isomorphic to the dual of
the cohomology with compact support of $\Delta_{\Hc}^*({\mc
  K}_{d,\rho})$.

\begin{conj}    \label{perverse}
The restriction to the diagonal $\Delta_{\Hc}^*({\mc K}_{d,\rho})$ is
a pure perverse sheaf.
\end{conj}

To compute the cohomology with compact support of
$\Delta_{\Hc}^*({\mc K}_{d,\rho})$, we consider the sheaf $(h_{d,\mu}^{\rm
  st})_{!}\Delta_{\Hc}^*({\mc K}_{d,\rho}) = (h_{d,\mu}^{\rm
  st})_{*}\Delta_{\Hc}^*({\mc K}_{d,\rho})$ on ${\mc A}_{d,\mu}^{\rm st}$
(recall that $h_{d,\mu}^{\rm st}$ is proper). By Deligne's purity
theorem, \conjref{perverse} implies that $(h_{d,\mu}^{\rm
  st})_{*}\Delta_{\Hc}^*({\mc K}_{d,\rho})$ is a pure complex. Hence,
geometrically, it is isomorphic to a direct sum of shifted simple
perverse sheaves.

\medskip

As we explained in \secref{decomp trace}, we wish to compare the {\em
  stable} trace formulas for the given group $G$ and the groups $H$
(depending on $\rho$ and $d$). Hence we need to isolate geometrically
a part in the cohomology of $(h_{d,\mu}^{\rm st})_{*}\Delta_{\Hc}^*({\mc
  K}_{d,\rho})$ which corresponds to the stable trace. For this, we
follow the strategy of \cite{Ngo:Endo,Ngo:FL}.

\medskip

Let $\pi_{0}({\mc P}_{d,\mu}^{\rm tor})$ denote the sheaf of connected
components of ${\mc P}_{d,\mu}^{\rm tor}$.  This is a sheaf of finite
abelian groups for the \'etale topology of $\Ac_{d,\mu}$. As in
\cite{Ngo:FL}, we define the {\em stable part} $$(h_{d,\mu}^{\rm
  st})_{*}\Delta_{\Hc}^*({\mc K}_{d,\rho})_{\rm st}$$ as the largest
direct factor of $(h_{d,\mu}^{\rm st})_{*}\Delta_{\Hc}^*({\mc K}_{d,\rho})$
on which $\pi_{0}^{\rm tor}({\rm P}_{d,\mu})$ acts unipotently.

Note that in the above formula the subscript ``st'' and the
superscript ``st'' have completely different meaning: the subscript
refers to the stable trace formula (see \secref{adelic} below), and
the superscript refers to imposing a stability condition in the sense
of geometric invariant theory.

By analogy with the main theorem of \cite{Ngo:FL, CL}, we state the
following conjecture.

\begin{conj}\label{support}
  The support of any simple perverse constituent of $(h_{d,\mu}^{\rm
    st})_{*}\Delta_{\Hc}^*({\mc K}_{d,\rho})_{\rm st}$ is
  $\Ac_{d,\mu}$. In particular, for every integer $n$, the perverse
  sheaf $$^p { H}^n ((h_{d,\mu}^{\rm st})_{*}\Delta_{\Hc}^*({\mc
    K}_{d,\rho})_{\rm st})$$ is the intermediate extension of its
  restriction to any non-empty open subscheme of $\Ac_{d,\mu}$.
\end{conj}

We will now explain the arithmetic meaning of this conjecture. 

\subsection{Ad\'elic description}    \label{adelic}

Let ${\mc M}_{d,\mu}(D)$ be the fiber of ${\mc M}_{d,\mu}$ over
$D=\sum_{i=1}^r n_i[x_i]$. If the curve $X$ is defined over a finite
field $k$, we can give an ad\'elic expression for the groupoid
$M_{d,\mu}(D)$ of $k$-points of ${\mc M}_{d,\mu}(D)$. Namely, let
\begin{align*}
\wt{M}_{d,\mu}(D) &= \{ \gamma,(g_x)_{x\in |X|}|\gamma \in G(F), (g_x)
\in G({\mathbb A})/ G(\Oc_{\mathbb A}), g_x^{-1} \gamma g_x \in
{G({\mathcal O}_x) \la_x(t_x) G({\mathcal O}_x)} \} \\ &
\subset G(F) \times G({\mathbb A})/ G(\Oc_{\mathbb A}),
\end{align*}
where $\la_x \leq d_x\mu$ if $x$ belongs to the support of $D$ and
$\lambda_x=0$, otherwise. The group $G(F)$ acts on $\wt{M}_{d,\mu}(D)$
by the formula
$$h\left(\gamma,(g_x)\right)=\left(h\gamma h^{-1}, (hg_x)_{x\in
  X}\right),$$
and
$$
M_{d,\mu}(D) = G(F) \backslash \wt{M}_{d,\mu}(D).
$$
Indeed, $(g_x)_{x\in |X|}$ determines a $G$-bundle $E$ on $X$ and
$\gamma$ determines a section $\varphi$ of $\on{Ad}(E)$ on $X -
\on{supp}(D)$ satisfying conditions (\ref{invariant 4}).

The map $h_{d,\mu}$ sends a point $(D,\gamma,(g_{x}))\in \M_{d,\mu}(k)$ 
to $(D,b)\in {\mc A}_{d,\mu}(k)$
where $b=(b_{1},\ldots,b_{r})$ with $b_{i}=\tr(\gamma)\in 
H^{0}(X,{\mc O}_{X}(\langle \mu,\omega_i \rangle D))$. 

Let $(D,b)\in\Ac_{d,\mu}^{\rm ani}(k)$. According to the Lefschetz
trace formula (see \secref{v sp}), the trace of the Frobenius
acting on the stalk of $(h_{d,\mu}^{\rm st})_{*}\Delta_{\Hc}^*({\mc
  K}_{d,\rho})$ over $(D,b)$ may be expressed as
\begin{equation}    \label{sum over gamma}
\sum_{\gamma \in G(F)/\on{conj.}} O_\gamma(K_{d,\rho,D}),
\end{equation}
where $\gamma$ runs over the set of conjugacy classes of $G(F)$ that
map to $b$.  The global orbital integral $O_\gamma(K_{d,\rho,D})$ is
defined by the formula
$$
O_\gamma(K_{d,\rho,D}) = \int_{G_\gamma(F)\bs G(\AD)}
H_{d,\rho,D}(g^{-1}\gamma g) dg
$$
where the $G_{\gamma}(F)$ stabilizer of $\gamma$ in $G(F)$ is a
discrete subgroup of $G(\AD)$ and $dg$ is the Haar measure on $G(\AD)$
such that $G(\Oc_{\AD})$ has volume one. The function under the
integration sign is
$$
H_{d,\rho,D} = \bigotimes_{x \in |X|} H_{\on{Sym}^{n_x(D)}(\rho)}
$$
on $G(\AD)$, so that we have
$$
K_{d,\rho,D}(x,y) = \sum_{a \in G(F)} H_{d,\rho,D}(x^{-1} a y)
$$
(see formula \eqref{K and H}). Here we write
$$
D = \sum_{x \in |X|} n_x(D) [x],
$$
and $H_{\on{Sym}^n(\rho)}$ is the element of the spherical Hecke
algebra of $G(\OO_x)$ bi-invariant functions on $G(F_x)$ corresponding
to the representation $\on{Sym}^{n}(\rho)$ under the Satake
isomorphism. In particular, for all but finitely many $x \in |X|$,
$H_{\on{Sym}^{n_x(D)}(\rho)}$ is the characteristic function of
$G(\OO_x)$, the unit element of the Hecke algebra.

As in \cite{Ngo:Endo}, for every $(D,b)\in \Ac_{d,\mu}^{\rm ani}(k)$
the trace of the Frobenius operator on the fiber of $(h_{d,\mu}^{\rm
  st})_{*}\Delta_{\Hc}^*({\mc K}_{d,\rho})_{\rm st}$ over $(D,b)$ can
be expressed in terms of stable orbital integrals
\begin{equation}
\prod_{x\in |X|} SO_{b}(H_{\on{Sym}^{n_x(D)}(\rho)}).
\end{equation}

If $(D,b)\notin \Ac_{d,\mu}^{\rm ani}(k)$ the expression is more
complicated. Nevertheless, following \cite{CL}, it is reasonable to
expect that it is the correct contribution of $b$ to the Arthur stable
trace formula.

\subsection{The case of $\SL_2$}

Let us consider again our example with $G=\SL_2$ and $\mu$ the highest
weight of the adjoint representation of ${\rm PGL}_2$. As cocharacter
of the maximal torus of $\SL_2$, $\mu$ is simply $t\mapsto {\rm
diag}(t,t^{-1})$. The moduli stack ${\mc M}_{d,\mu}$ classifies the
data
$$(D,E,\varphi)$$ where $D\in X_d$ is an effective divisor on $X$ of
degree $d$, $E$ is a rank two vector bundle over $X$ equipped with a
trivialization of its determinant and $\varphi$ is an automorphism of
$E|_{X{-}\on{supp}(D)}$ of determinant one that can be extended to a
homomorphism of vector bundles
$$\varphi_1:E\to E(D).$$

Let ${\mc A}_{d,\mu}$ be a fibration over $X_d$, whose fiber over
$D\in X_d$ is the vector space $H^0(X,\Oc_X(D))$. We have a
Hitchin-like map
$$h_{d,\mu}:{\mc M}_{d,\mu} \to {\mc A}_{d,\mu},$$ sending
$(D,E,\varphi)$ to $(D,b)$, where
$$b=\tr(\varphi_1) \in H^0(X,\Oc_X(D)).$$

Using Riemann-Roch theorem, it is easy to see that 
$$\dim(\Ac_{d,\mu})=2d-g+1$$
for every $d \geq 2g-1$. This formula is in fact true for all $d\geq
g$ and we have a formula for $d\leq g$ as well.

Recall that $\Ac_{d,\mu}$ classifies all pairs $(D,b)$ with $D\in X_d$
and $b\in H^0(X,\Oc_X(D))$. If $b$ is a non-zero section of
$H^0(X,\Oc_X(D))$, the divisor $D'=\dv(b)+D$ is an effective divisor
of degree $d$ that is linearly equivalent to $D$. Here $\dv(b)$ is the
principal divisor attached to $b$ viewed an a non-zero rational
function on $X$. Such a $D'$ determines $b$ up to a scalar in
$k^\times$. If $\Ac_{d,\mu}^\times$ denote the complement of the zero
section $b=0$, then we have a map
$$\Ac_{d,\mu}^\times \to X_d\underset{\on{Pic}_d}\times X_d$$
given by $(D,b)\mapsto (D,D')$ which is a $\Gm$-torsor. Here $X_d\to
{\rm Pic}_d$ is the Abel-Jacobi map from $X_d$ to the $d$-th component
of Picard's variety (if we replace the Picard variety by the Picard
stack, the above map is an isomorphism). In particular, we have the
dimension formula
$$\dim(\Ac_{d,\mu})=\dim(X_d\underset{\on{Pic}_d}\times X_d)+1.$$
Now, we recall a classical formula from the theory of special divisors
on compact Riemann surfaces.

\begin{lem}[Martens] \label{Martens} We have the following formula for
  the dimension of $X_d\underset{\on{Pic}_d}\times X_d$
\begin{equation}
\dim(X_d\underset{\on{Pic}_d}\times X_d)=
\begin{cases} 
d      & \text{if $1\leq d\leq g-1$}, \\
2d-g & \text{if $g\leq d$}.
\end{cases}
\end{equation}
Here $g$ denotes the genus of $X$.
\end{lem}

\begin{proof}
  Suppose first $1\leq d\leq g-1$.  According to Martens' theorem
  \cite[p.191]{ACGH}, for every integer $d$ such that $1\leq d\leq
  g-1$, the Abel-Jacobi morphism $X_d \to {\rm Pic}_d$ from $X_d$ to
  the $d$-th component of he Picard variety is semi-small over its
  image. In other words,
$$\dim(X_d\underset{\on{Pic}_d}\times X_d)=d.$$ 

Suppose next $g\leq d\leq 2g-2$. In this case, the map $X_d\to {\rm
  Pic}_d$ is surjective and its generic fiber is a projective space of
dimension $d-g$. It follows that the irreducible component of
$X_d\underset{\on{Pic}_d}\times X_d$ that dominates ${\rm Pic}_d$ has
dimension $2d-g$. It is enough to prove that the other components, if
any, have less dimension. Let $r$ be a positive integer and let denote
by $W_d^r$ the locally closed subvariety of $\Pic_d$ consisting in
line bundles $L$ of degree $d$ such that
$$\dim H^0(X,L)=d+1-g+r.$$
By Riemann-Roch theorem, $\dim H^0(X,L')=r$ where $L'=L^{-1}\otimes K$
where $K$ is the canonical sheaf. We have $\deg(L')=d'=-d+2g-2$. After
Martens' theorem the subvariety of line bundles $L'$ of degree $d'$
such that $\dim(H^0(X,L'))=r$ is less or equal to $d'-2r+2$. It
follows that the part of $X_d \underset{\on{Pic}_d}\times X_d$ over
$W_d^r$ is of dimension at most $$2(d-g+r)+d'-2r+2=d,$$ which is less
than $2d-g$.

Finally, if $d\geq 2g-1$, then for all $D\in X_d$, $\dim
H^0(X,\Oc_X(D))=d+1-g$ after Riemann-Roch theorem. It follows that the
map $X_d\to {\rm Pic}_d$ is a projective bundle of rank $d-g$ and the
total dimension of $X_d \underset{\on{Pic}_d}\times X_d$ is $2d-g$.
\end{proof}

We can use Hitchin's device called {\em spectral curve} to describe
the fibers of $h_{d,\mu}$. Recall that for any sections $b_1\in
H^0(X,\Oc_X(D))$ and $b_2\in H^0(\Oc_X(2D))$, we have a curve on the
total space of the line bundle $\Oc_X(D)$ defined by the equation
$$t^2-b_1 t+ b_2=0.$$
For every $b\in H^{0}(X,\Oc_{X(D)}$, let $Y_{D,b}$ be the spectral
curve corresponding to particular parameters $b_1=b$ and $b_2=1_{2D}$
the constant function $1$ on $X$ considered as a global section of the
line bundle $\Oc_X(2D)$. The map $\varphi:E\to E(D)$ defines the
structure of a module over the symmetric $\Oc_X$-algebra of
$\Oc_X(-D)$ on $E$, so that $E$ can be seen as an $\Oc$-module on the
total space of $\Oc_X(D)$. The Cayley--Hamilton theorem implies that
$E$ is supported on the spectral curve $Y_{D,b}$.

For every $(D,b)\in {\mc A}_{d,\mu}$, the fiber 
$${\mc M}_{D,b}=h_{d,\mu}^{-1}(D,b)$$ 
classifies $\Oc_{Y_{D,b}}$-modules ${\mc F}$ such that by pushing along
the finite flat map $p_{D,b}:Y_{D,b}\to X$, we get a rank two locally
free $\Oc_X$-module $p_{D,b*}{\mc F}$, equipped with a trivialization
of the determinant. If $Y_{D,b}$ is reduced, the fact that
$p_{D,b*}{\mc F}$ is a rank two vector bundle implies that ${\mc F}$
is a torsion-free module of generic rank one and vice versa. If
$Y_{D,b}$ is smooth, ${\mc F}$ is in fact an invertible sheaf.

In the present $\SL_2$ case, the group ${\mc P}_{D,b}$ consists of
invertible sheaves ${\mc L}$ on $Y_{D,b}$ with the trivial norm down
to $X$. It acts on ${\mc M}_{D,b}$ by tensor product
$$({\mc L},{\mc F}) \mapsto {\mc L}\otimes_{\Oc_{Y_{D,b}}} {\mc F}$$
because 
$$\det(p_{D,b*}({\mc L}\otimes_{\Oc_{Y_{D,b}}} {\mc F}))={\rm Nm}_{Y_{D,b}/X}
({\mc L}) \otimes_{\Oc_X} \det(p_{D,b*}({\mc F})).$$

For $b\in H^0(X,\Oc_X(D))$ such that the section $b^2-4_{2D}$ of the
line bundle $\Oc_X(2D)$ has a multiplicity free divisor, the spectral
curve $Y_{D,b}$ is smooth. In this case, ${\mc P}_{D,b}$ is isomorphic
to the quotient of an abelian variety by the trivial action of
$\Z/2\Z$. In particular, ${\mc P}_{D,b}$ is proper and connected. It
is also known that ${\mc P}_{D,b}$ acts simply transitively on ${\mc
  M}_{D,b}$ in this case.

For a general parameter $(D,b)$, as in \cite{Ngo:FL}, we need to keep
track of two invariants attached to $\Pc_{D,b}$: the component group
$\pi_{0}(\Pc_{D,b})$ and its invariant $\delta$.

The component group $\pi_{0}(\Pc_{D,b})$ is related to the theory of
endoscopy. In the present situation, it is fairly easy to compute. If
the spectral curve $Y_{D,b}$ has at least one unibranch ramification
point, in particular irreducible, then $\pi_{0}(\Pc_{D,b})=0$. If this
curve has only ramification points with two branches but still
irreducible, then $\pi_{0}(\Pc_{D,b})=\Z/2\Z$. This is the case if and
only if the normalization of $Y_{D,b}$ is an unramified covering of
$X$. Finally, if the spectral curve $Y_{D,b}$ has two irreducible
components, then $\pi_{0}(\Pc_{D,b})=\Z$. See the last section of
\cite{Ngo:Endo} for this calculation.

The $\delta$ invariant $\delta(D,b)$ which is defined as the dimension
of the affine part of the group $\Pc_{D,b}$ is a rough but efficient
measure of singularity of ${\mc M}_{D,b}$. This invariant can be
calculated from the discriminant as follows. Let ${\rm discr}(D,b)$ be
the effective divisor attached to the section $b^2-4_{2D}$ of
$\Oc_X(2D)$. In other words, we have
$$ {\rm discr}(D,b)=\dv(b^2-4_{2D})+2D$$
where $\dv(b^2-4_{2D})$ is the (virtual) divisor attached to the rational
function $b^2-4_{2D}$. Let us write ${\rm discr}(D,b)$ in the form ${\rm
  discr}(D,b)=\sum_{x\in X} n_x[x]$. Then
\begin{equation}\label{delta}
 \delta(D,b)=\sum_{x\in X} \left[\frac{n_x}{2}\right]
\end{equation}
where $\left[\frac{n_x}{2}\right]$ denotes the integer part of
$\frac{n_x}{2}$. A convenient way to express $\delta(D,b)$ from ${\rm
  discr}(D,b)$ is as follows : we write ${\rm discr}(D,b)$ under the
form $D_1+2D_2$ where $D_1,D_2$ are effective divisors and $D_1$ is
multiplicity free. Then
$$\delta(D,b)=\deg(D_2).$$
We will stratify $\Ac_{d,\mu}$ by this invariant $\delta$
$$\Ac_{d,\mu}=\bigsqcup_{\delta} \Ac_{d,\mu}(\delta)$$
where $\Ac_{d,\mu}(\delta)$ is the locally closed subscheme of
$\Ac_{d,\mu}$ consisting of points $(D,b)$ such that
$\delta(D,b)=\delta$. In view of some general results proved in
\cite{Ngo:FL}, the following estimate of dimensions provides some
evidence in favor of Conjecture \ref{support}.

\begin{prop}\label{dimension} 
  Assume $d\geq 2g-1$ and $\delta\leq d-2g+1$. Then the codimension in
  $\Ac_{d,\mu}$ of the stratum $\Ac_{d,\mu}(\delta)$ is equal to
  $\delta$.
\end{prop}

\begin{proof} For every integer $e\leq d$, let denote
  $\Ac_{d,\mu}(\delta,e)$ the locally closed subscheme of
  $\Ac_{d,\mu}(\delta)$ defined by the condition $\delta(D,b)=\delta$
  and the degree of the zero divisor of the rational function $b$
  equals $e$. Let $(D,b)$ be a geometric point of
  $\Ac_{d,\mu}(\delta,e)$. Let us denote $E_+=\dv(b)_+$ the divisor of
  zeros of $b$ and $E_-=\dv(E)_-$ the divisor of poles; we have then
  $\deg(E_+)=\deg(E_-)=e$.  The requirement that the $\dv(b)+D$ be
  effective, implies that $D=E_-+D_0$ for a certain effective divisor
  $D_0$ of degree $\deg(D_0)=d_0=d-e$. Consider the discriminant
  divisor
$${\rm discr}(D,b)=\dv(b^2-4_{2D})+2D$$
where 
$$\dv(b^2-4_{2D})=\dv(b^2-4_{2D})_+-2E_-$$
since $\dv(b^2-4_{2D})_-=\dv(b^2)$ and therefore
$${\rm discr}(D,b)=\dv(b-2_D)_+ + \dv(b+2_D)_+ + 2D_0.$$
We observe that $b$ defines a finite flat morphism $b:X\to \mathbb
P^1$ of degree $e$ and in terms of this map $\dv(b-2_D)_-=b^{-1}(2)$ and
$\dv(b+2)=b^{-1}(-2)$. In particular, these effective divisors are
disjoint.

Now we write $\dv(b-2_D)_+=D'_1+2D'_2$ and $\dv(b+2_D)_+=D''_1+2 D''_2$
where $D'_1$ and $D''_1$ are multiplicity free effective divisors. We
will use the obvious notation $d'_1=\deg(D'_1)$, $d'_2=\deg(D'_2)$
etc. Then we have
$$\delta=d'_2+ d''_2+d_0$$
since $D'_1+D''_1$ is also a multiplicity free divisor. 

We consider the map from $\Ac_{d,\mu}(\delta,e)$ to
$$X_{d'_1}\times X_{d'_2} \times X_{d''_1} \times X_{d''_2} \times X_{d_0}.$$
This is a $\Gm$-torsor over the closed subset of the above product
defined by only one condition: $D'_1+2D'_2$ and $D''_1+2 D''_2$ are
linearly equivalent. Now this subvariety of the quintuple product is
smooth of dimension
$$d'_1+d'_2+d''_1+d''_2+d_0 -g=2d-g-\delta$$
if at least one of the integers $d'_1,d'_2,d''_1,d''_2$ is larger or
equal to $2g-1$. We prove that either $d'_{1}$ or $d''_{1}$ is larger
than equal to $2g-1$.

Recall the equalities
$$d'_1+2d'_2+d_0=d''_1+2d''_2+d_0=d$$
and 
$$\delta=d'_2+d''_2+d_0.$$ 
Under the assumption $d-\delta\geq 2g-1$, we can proceeds as follows:
assume $d'_2 \leq d''_2$, then $2d'_2+d_0 \leq \delta$ so that the
inequality
$$d'_1 \geq d-\delta \geq 2g-1$$
is satisfied. It follows that the dimension of $\Ac_{d,\mu}(\delta,e)$
is $2d-g+1-\delta$.
\end{proof}

\subsection{The case of the twisted torus}

Now we consider a non-split example. Namely, let $H$ be the
one-dimensional torus over $X$ attached to a $\mu_{2}$-torsor
$\Ec_{H}:X'\to X$ as in \secref{example 1}. The dual group of $H$ is
$^{L} H = \Gm\rtimes \Gamma$ where $\Gamma$ acts non-trivially on
$\Gm$ through the quotient corresponding to $\Ec_{H}$. We have defined
a homomorphism $^{L} H \to \PGL_2\times \Gamma$ in \eqref{hom L}. We
will consider the restriction $\rho_{H}$ of the representation
$\rho=\rho_{\mu}\times \rho_{\Ec_{H}}$ of $\PGL_2\times \Gamma$ to
$^{L} H$, where $\rho_{\mu}$ is the adjoint representation of $\PGL_2$
and $\rho_{\Ec_{H}}$ is the representation of $\Gamma$ corresponding
to $\Ec_{H}$. The three-dimensional representation $\rho_H$ breaks
into a direct sum of a one-dimensional and a two-dimensional
representations:
$$
\rho_{H}=\rho_{0}\oplus \rho_{1}.
$$
We have seen earlier that the perverse sheaf $\wF_{d,\rho}$ also
breaks into a direct sum
\begin{equation}
	\wF_{d,\rho_{H}}=\bigoplus_{d_{0}+d_{1}=d} \und\Ql |
	(\Bun_{H}\times X_{d_{0}} \times X'_{d_{1}}).
\end{equation}

The relevant moduli stack ${\mc M}_{d,[\mu]}$, which we will denote
here by $\M_{d,\rho_{H}}$, is a disjoint union
\begin{equation}
\M_{d,\rho_{H}}=\bigsqcup_{d_{0}+d_{1}=d}
\M_{d,\rho_{H}}(d_{0},d_{1}),
\end{equation}
where $\M_{d,\rho_{H}}(d_{0},d_{1})$ classifies the data
$(L',D_{0},D_{1},\theta)$, where
$$(L',D_{0},D_{1}) \in \Bun_{H}\times X_{d_{0}} \times X'_{d_{1}}$$
and $\theta$ is an isomorphism of line bundles
$$\theta:L(D_{1}-\tau(D_{1})) \to L$$
such that $\tau(\theta)\theta=1$. Recall that $\tau$ is the
involution on $X'$ such that $X'/\tau = X$.
We derive from $\theta$ a symmetric isomorphism $\Oc_{X'}(D_1) \to
\tau^*(\Oc_{X'}(D_1))$ or, equivalently a descent datum of the line
bundle $\Oc_{X'}(D_1)$ to $X$. This defines a point in the cartesian
product $X'_{d_1} \times_{\Pic_{d_1} (X')} \Pic_{d_1/2}(X)$.  Let us
denote
$$\Ac_{d,\rho_{H}}(d_{0},d_{1})= X_{d_0} \times (X'_{d_1}
\times_{\Pic_{d_1} (X')} \Pic_{d_1/2}(X)).$$
This space is empty for odd integers $d_1$.
The Hitchin base is then the disjoint union
$$\Ac_{d,\rho_{H}}=\bigsqcup_{d_{0}+d_{1}=d}
\Ac_{d,\rho_{H}}(d_{0},d_{1}).$$ 
The Hitchin map,
restricted to each component $\M_{d,\rho_{H}}(d_{0},d_{1})$, is just
the projection
$$\M_{d,\rho_{H}}(d_{0},d_{1})=\Bun_{H}\times
\Ac_{d,\rho_{H}}(d_{0},d_{1}) \to \Ac_{d,\rho_{H}}(d_{0},d_{1})$$
that has relative dimension $\dim(\Bun_{H})=g-1$.

\subsection{Conjecture}

We recall that our goal is to analyze the (stable parts) of the
cohomologies of the restrictions of the sheaves ${\mc K}_{d,\rho}$ and
${\mc K}_{d,\rho_H}$ to ${\mc M}_{d,\rho}$ and ${\mc M}_{d,\rho_H}$,
respectively, and to use them to establish the desired identities
\eqref{sum2}. We return to our main example (see \secref{ex sl2}):
$G=\SL_2$, $\rho = \rho_\mu \otimes \Lc$, where $\rho_\mu$ is the
adjoint representation of $^L G = \PGL_2$ and $\Lc$ is an order two
local system on $X$, and $H$ is the twisted torus over $X$
corresponding to $\Lc$. We state the following conjecture in this case.

\medskip

Recall that $\pi_{0}(\Bun_{H})=\Z/2\Z$ and let us denote by
$\Bun_{H}^{0}$ the neutral component of $\Bun_{H}$. Let us also denote
by $\tau$ the involution of $\Ac_{d,\rho_{H}}$ given by
$$\tau(L',D_{0},D_{1},\theta)=(L',D_{0},
\tau(D_{1}),\tau(\theta)).$$ This involution acts on the cohomology of
$\Ac_{d,\rho_{H}}$ (with compact support) and we will denote by the
upper script $\tau$ the space of invariants of this action.

\begin{conj}
  For sufficiently large $d$ there exists a quasi-isomorphism between
  complexes of $\ell$-adic vector spaces equipped with Frobenius
  operators
$$R\Gamma_c(\Ac_{d,\mu},(h_{d,\mu}^{\rm st})_{*}\Delta_{\Hc}^*({\mc
  K}_{d,\rho})_{\rm st}) \simeq R\Gamma_{c}(\Bun_{H}^{0},\Ql) \otimes
R\Gamma_c(\Ac_{d,\rho_{H}},\Ql)^{\tau},$$
modulo the contribution of the Eisenstein series mentioned at the end 
of the \secref{ex sl2}.
\end{conj}

We remark that the vector space on the right hand side may be viewed
as a subspace of $R\Gamma_c(\M_{d,\rho_{H}},\Ql)$, since
$\M_{d,\rho_{H}} \simeq \Bun_{H} \times \Ac_{d,\rho_{H}}$. Taking
$\tau$-invariants is due to the fact that the map $\Ac_{d,\rho_{H}}
\to \Ac_{d,\rho}$ is generically two-to-one, and taking the neutral
component $\Bun_{H}^{0}$ is the geometric counterpart of the factor
$1/2$ in the corresponding stable trace formula.

We hope that this conjecture may provide a starting point for further
investigation of trace formulas by geometric methods.

\newpage

\vspace*{10mm}

\begin{center}

{\bf Part II}

\end{center}

\bigskip

In this, more speculative, part of the paper we propose two
conjectural ``geometric trace formulas'' in the case of a complex
algebraic curve $X$.

\section{The geometric trace formula}    \label{leap}

In \secref{v sp} we have interpreted geometrically the right hand side
of the trace formula \eqref{tr}. Now we turn to the left hand
side. From now on we will assume that $G$ is a constant group scheme
$\GG \times X$, where $\GG$ is a split reductive group defined over
the ground field $k$. Henceforth, in order to simplify our notation,
we will not distinguish between $G$ and $\GG$.

\subsection{The left hand side of the trace formula}

As we discussed in \secref{h eig}, the $L$-packets of (unramified)
irreducible automorphic representations should correspond to
(unramified) homomorphisms $W_F \times SL_2 \to {}^L G$. Assuming this
conjecture and ignoring for the moment the contribution of the
continuous spectrum, we may write the left hand side of
\eqref{general} as
\begin{equation}    \label{trace1}
\on{Tr} {\mb K} = \sum_\sigma m_\sigma N_\sigma,
\end{equation}
where $\sigma$ runs over the unramified homomorphisms $W_F \times SL_2
\to {}^L G$, $m_\sigma$ is the multiplicity of the irreducible
automorphic representation unramified with respect to $G(\OO)$ in the
$L$-packet corresponding to $\sigma$, and $N_\sigma$ is the eigenvalue
of the operator ${\mb K}$ on an unramified automorphic function
$f_\sigma$ on $\Bun_G(\Fq)$ corresponding to a spherical vector in
this representation:
$$
{\mb K} \cdot f_\sigma = N_\sigma f_\sigma.
$$

Thus, recalling \eqref{sum b}, the trace formula \eqref{general}
becomes
\begin{equation}    \label{classical}
\sum_{\sigma: W_F \times SL_2 \to {}^L G} m_\sigma N_\sigma =
\sum_{P \in \Bun_G(\Fq)} \frac{1}{|\on{Aut}(V)|} K(P,P).
\end{equation}

Consider, for example, the case of ${\mb K} = {\mb H}_{\rho,x}$, the
Hecke operator corresponding to a representation $\rho$ of $^L G$ and
$x \in |X|$. Then, according to the conjectures of Langlands and
Arthur,
\begin{equation}    \label{Nsigma}
N_\sigma = \on{Tr} \left(\sigma\left(\begin{pmatrix}
        q_x^{1/2} & 0 \\ 0 & q_x^{-1/2} \end{pmatrix} \times
      \on{Fr}_x\right),\rho \right).
\end{equation}
The operators ${\mb K}$ introduced in \secref{int tr} are generated by
the Hecke operators ${\mb H}_{\rho,x}$. Therefore the eigenvalue
$N_\sigma$ for such an operator ${\mb K}$ is expressed in terms of the
traces of $\sigma(\on{Fr}_x)$ on representations of $^L G$.

\subsection{Lefschetz fixed point formula
  interpretation}    \label{lfpf}

It is tempting to try to interpret the left hand side of
\eqref{classical} as coming from the Lefschetz trace formula for the
trace of the Frobenius on the cohomology of an $\ell$-adic sheaf on a
moduli stack, whose set of $k$-points is the set of $\sigma$'s.
Unfortunately, such a stack does not exist if $k$ is a finite field
$\Fq$ (or its algebraic closure).  On the other hand, if $X$ is over
$\C$, then there is an algebraic stack $\Loc_{^L G}$ of (de Rham) $^L
G$-local systems on $X$; that is, $^L G$-bundles on $X$ with flat
connection. But in this case there is no Frobenius acting on the
cohomology whose trace would yield the desired number (the left hand
side of \eqref{classical}).  Nevertheless, we will define a certain
vector space (when $X$ is over $\C$), which we will declare to be a
geometric avatar of the left hand side of \eqref{classical} (we will
give an heuristic explanation for this in \secref{abl fixed}). We will
then conjecture that this space is isomorphic to \eqref{rhs-new} --
this will be the statement of the ``geometric trace formula'' that we
propose in this paper.

\medskip

Let us first consider the simplest case of the Hecke operator ${\mb
  K}={\mb H}_{\rho,x}$. In this case the eigenvalue $N_\sigma$ is
given by formula \eqref{Nsigma}, which is essentially the trace of the
Frobenius of $x$ on the vector space which is the stalk of the local
system on $X$ corresponding to $\sigma$ and $\rho$. These vector
spaces are fibers of a natural vector bundle on $X \times \Loc_{^L G}$
(when $X$ is defined over $\C$).

Indeed, we have a tautological $^L G$-bundle ${\mc T}$ on $X \times
\on{Loc}_{^L G}$, whose restriction to $X \times \sigma$ is the $^L
G$-bundle on $X$ underlying $\sigma \in \on{Loc}_{^L G}$. For a
representation $\rho$ of $^L G$, let ${\mc T}_\rho$ be the associated
vector bundle on $X \times \on{Loc}_{^L G}$. It then has a partial
flat connection along $X$. Further, for each point $x \in |X|$, we
denote by ${\mc T}_x$ and ${\mc T}_{\rho,x}$ the restrictions of ${\mc
  T}$ and ${\mc T}_\rho$, respectively, to $x \times \on{Loc}_{^L G}$.

It is tempting to say that the geometric incarnation of the left hand
side of \eqref{classical} in the case ${\mb K}={\mb H}_{\rho,x}$ is
the cohomology
$$
H^\bullet(\Loc_{^L G},{\mc T}_{\rho,x}).
$$
However, this would only make sense if the vector bundle ${\mc
  T}_{\rho,x}$ carried a flat connection (i.e., a ${\mc D}$-module
structure) and this cohomology was understood as the de Rham
cohomology (which is the analogue of the \'etale cohomology of an
$\ell$-adic sheaf that we now wish to imitate when $X$ is defined over
$\C$).

Unfortunately, ${\mc T}_{\rho,x}$ does not carry any natural
connection. If we have a vector bundle ${\mc V}$ with a flat
connection on a smooth algebraic variety $Y$, then its de Rham
cohomology may be defined as the {\em coherent} cohomology of the
tensor product
$$
{\mc V} \otimes \Lambda^\bullet(TY),
$$
where $TY$ is the tangent sheaf to $Y$. This cohomology is
well-defined even if ${\mc V}$ does not have a flat connection, and
so, by making a leap of faith, we can declare the cohomology
\begin{equation}    \label{geom trace}
H^\bullet({\on{Loc}_{^L G}},{\mc T}_{\rho,x} \otimes
\Lambda^\bullet(T\on{Loc}_{^L G}))
\end{equation}
as the geometric incarnation of the left hand side of
\eqref{classical}.

However, this formula is only the first approximation to the right
formula, because $\on{Loc}_{^L G}$ is not a variety, but an algebraic
stack. We suggest that the correct formula is
$$
\on{RHom}(\Delta_*(\OO),\Delta_*({\mc T}_{\rho,x})) =
\on{RHom}(\Delta^*\Delta_*(\OO),{\mc T}_{\rho,x}),
$$
where $\Delta$ is the diagonal. This will be discussed further in
Sections \ref{appl} and \ref{abl fixed}.

On the other hand, as explained in \secref{v sp}, the geometric avatar
of the right hand side of \eqref{classical} is the (de Rham)
cohomology
$$
H^\bullet(\Bun_G,\Delta^!(\ol{\mc K}_{\rho,x})).
$$

The conjectural geometric trace formula for $\K={\mathbb H}_{\rho,x}$
would then be an isomorphism
\begin{equation}    \label{gtf for hecke}
\on{RHom}(\Delta^*\Delta_*(\OO),{\mc T}_{\rho,x}) \simeq
H^\bullet(\Bun_G,\Delta^!(\ol{\mc K}_{\rho,x})).
\end{equation}

\subsection{The case of $\K=\K_{d,\rho}$}    \label{case K}

Now we generalize this to other integral operators ${\mb K}$, in
particular, the one that is of most interest to us, ${\mb K} = {\mb
  K}_{d,\rho}$ (see formula \eqref{kernel of kdrho}). In this case
$N_\sigma$ on the left hand side of the trace formula
\eqref{classical} is given by the formula
$$
\sum_{D=\sum_i n_i[x_i] \in
 X_d(\Fq)} \quad \prod_i \on{Tr} \left(\sigma\left(\begin{pmatrix}
        q_{x_x}^{1/2} & 0 \\ 0 & q_{x_i}^{-1/2} \end{pmatrix} \times
      \on{Fr}_{x_i}\right), \on{Sym}^{n_i}(\rho) \right).
$$
The corresponding sheaf on $\Loc_{^L G}$ is defined by symmetrization
of the vector bundles ${\mc T}_\rho$, as follows.

First, let us recall from \secref{symmetric power} that if $E$ is a
flat vector bundle on $X$, then we define its $d$th symmetric power
$E_d$ by the formula
\begin{equation}    \label{s p}
E_d = \left( \pi^d_* (E^{\boxtimes d}) \right)^{S_d},
\end{equation}
where $\pi^d: X^d \to X_d$ is the symmetrization map. This is a ${\mc
  D}$-module on $X_d$ (but not a vector bundle if the rank of $E$ is
greater than $1$).  We may perform the same construction over any base
$B$ parametrizing a family of flat vector bundles. In other words,
suppose that we have a vector bundle ${\mc E}$ over $X \times B$
equipped with a flat connection along $X$. Then we obtain a coherent
sheaf ${\mc E}_d$ on $X_d \times B$ which is moreover a ${\mc
  D}$-module along $X_d$. The restriction of ${\mc E}_d$ to $X_d
\times b$ is $E_d$ given by formula \eqref{s p}, where $E = {\mc
  E}|_{X \times b}$.

Let us apply this in the case of $B = \Loc_{^L G}$. If $\rho$ is again
a representation of $^L G$, we have a tautological vector bundle ${\mc
  T}_\rho$ on $X \times \Loc_{^L G}$ equipped with a flat connection
along $X$. Let ${\mc T}_\rho^d$ be the sheaf on $X^d \times \Loc_{^L
  G}$ obtained by taking the Cartesian power of ${\mc T}_\rho$ along
$X$. It carries a flat connection along $X^d$. Denote by $\wt\pi^d$
the symmetrization morphism $X^d \times \Loc_{^L G} \to X_d \times
\Loc_{^L G}$ along the first factor and set
$$
{\mc T}_{d,\rho} = \left( \wt\pi^d_* ({\mc T}_\rho^{d}) \right)^{S_d}.
$$
This is a coherent sheaf on $X_d \times \Loc_{^L G}$ which carries
the structure of a ${\mc D}$-module along $X_d$. For each
effective divisor $D = \sum_i n_i[x_i]$ on $X$, the restriction of
${\mc T}_{d,\rho}$ to $D \times \Loc_{^L G}$ is the vector bundle
$\bigotimes_i {\mc T}^{(n_i)}_{\rho,x_i}$, where ${\mc
  T}^{(n)}_{\rho,x_i}$ is the tautological vector bundle on $\Loc_{^L
  G}$ corresponding to the representation $\on{Sym}^{n}(\rho)$ and
$x_i \in |X|$.

Let $\pi: X_d \times \Loc_{^L G} \to \Loc_{^L G}$ be the
projection onto the second factor. Then
\begin{equation}    \label{Fdrho}
{\mc F}_{d,\rho} = \pi_!({\mc T}_{d,\rho})
\end{equation}
(where the direct image corresponds to taking the fiberwise de Rham
cohomology) is the coherent sheaf on $\Loc_{^L G}$ that captures the
functor $\K=\K_{d,\rho}$ the way the vector bundle ${\mc T}_{\rho,x}$
captures the Hecke functor $\K={\mathbb H}_{\rho,x}$. Hence, following
the same argument as in the case of $\K={\mathbb H}_{\rho,x}$, we may
view the vector space
\begin{equation}    \label{declare}
H^\bullet({\on{Loc}_{^L G}},{\mc F}_{d,\rho} \otimes
\Lambda^\bullet(T\on{Loc}_{^L G}))
\end{equation}
as the first approximation to the geometric avatar of the left hand
side of the trace formula \eqref{classical} for $K=K_{d,\rho}$ (it is
approximate because $\on{Loc}_{^L G}$ is an algebraic stack).

\bigskip

On the other hand, as explained in \secref{v sp}, the geometric avatar
of the right hand side of \eqref{classical} is the (de Rham)
cohomology
$$
H^\bullet(\Bun_G,\Delta^!(\ol{\mc K}_{d,\rho})).
$$
Therefore, in the first approximation, we obtain the following
conjectural geometrization of the trace formula \eqref{classical}:
\begin{equation}    \label{geometric1}
H^\bullet({\on{Loc}_{^L G}},{\mc F}_{d,\rho} \otimes
\Lambda^\bullet(T\on{Loc}_{^L G})) \simeq
H^\bullet(\Bun_G,\Delta^!(\ol{\mc K}_{d,\rho})).
\end{equation}

However, since $\Loc_{^L G}$ is an algebraic stack in general, the
left hand side of \eqref{geometric1} has to be modified. In the next
section we will argue that \eqref{geometric1} should be replaced by
the isomorphism \eqref{geometric2} which we will propose as the
sought-after ``geomeric trace formula''. We will also show that this
isomorphism naturally appears as a corollary of the categorical
Langlands correspondence (also known as ``non-abelian Fourier--Mukai
transform'').

\section{Categorical Langlands correspondence}    \label{clc}

According to the geometric Langlands conjecture \cite{BD} (see, e.g.,
\cite{F:rev} for an exposition), for each (sufficiently generic) $^L
G$-local system $\sigma$ on $X$ there exists a ${\mc D}$-module (or
perverse sheaf) ${\mc F}_\sigma$ on $\Bun_G$ which is a Hecke
eigensheaf with ``eigenvalue'' $\sigma$. This property means that we
have a system of isomorphisms
\begin{equation}    \label{Hecke}
\H_\rho({\mc F}_\sigma) \simeq (\rho \circ \, \sigma) \boxtimes {\mc
  F}_\sigma
\end{equation}
compatible with the tensor product structures, where $\H_\rho$ is the
Hecke functor associated to the representation $\rho$ of $^L G$,
acting from the category of ${\mc D}$-modules on $\Bun_G$ to the
category of ${\mc D}$-modules on $X \times \Bun_G$ (see \secref{aff
  Gr}).

\subsection{Equivalence of categories}

Categorical Langlands correspondence extends this to a kind of
spectral decomposition of the derived category of ${\mc D}$-modules
on $\Bun_G$. More precisely, one hopes that there exists an
equivalence of derived categories\footnote{\label{foot} It is expected
  (see \cite{FW}, Sect.~10) that in general there is not a canonical
  equivalence, but rather a $\Z_2$-gerbe of such equivalences. This
  gerbe is trivial, but not canonically trivialized. One gets a
  particular trivialization of this gerbe, and hence a particular
  equivalence $C$, for each choice of the square root of the canonical
  line bundle $K_X$ on $X$.}

\medskip

\begin{equation}    \label{na fm}
\boxed{\begin{matrix} \text{derived category of} \\
    \OO\text{-modules on } \on{Loc}_{^L G} \end{matrix}} \quad
    \longleftrightarrow \quad \boxed{\begin{matrix}
    \text{derived category of} \\ \D\text{-modules on } \Bun_G
    \end{matrix}}
\end{equation}

\bigskip

\noindent under which the skyscraper coherent sheaf supported at
$\sigma \in \on{Loc}_{^L G}$ goes to ${\mc F}_\sigma$. This has been
proved in the abelian case by Laumon \cite{Laumon} and Rothstein
\cite{Roth}, as a version of the Fourier--Mukai transform. In the
non-abelian case this categorical correspondence (``non-abelian
Fourier--Mukai transform'') has been suggested by Beilinson and
Drinfeld (see, e.g., \cite{F:rev,Laf}). It also follows from the
S-duality/Mirror Symmetry picture of Kapustin and Witten \cite{KW}
(see \cite{FW,F:bourbaki} for an exposition). It has not yet been made
into a precise conjecture in the literature, so we only use it as a
guiding principle.\footnote{D. Arinkin, V. Drinfeld, and D. Gaitsgory
  have announced in recent talks a more precise formulation of this
  conjecture.}

In what follows, we will denote by $D(\on{Loc}_{^L G})$ and
$D(\Bun_G)$ the categories that should appear on the two sides of
\eqref{na fm}.

This conjectural equivalence has been mostly studied at the level of
objects (for instance, in the case when $C$ is applied to the
skyscraper sheaf $\OO_\sigma$ on $\Loc_{^L G}$, see the next
subsection). But it should also yield important information at the
level of morphisms.  In particular, if we denote the equivalence going
from left to right in the above diagram by $C$, then we should have
isomorphisms
\begin{equation}    \label{RHom}
\on{RHom}_{D(\on{Loc}_{^L G})}({\mc F}_1,{\mc F}_2) \simeq 
\on{RHom}_{D(\Bun_G)}(C({\mc F}_1),C({\mc F}_2))
\end{equation}
and
\begin{equation}    \label{RHom1}
\on{RHom}({\mathbb F}_1,{\mathbb F}_2) \simeq 
\on{RHom}(C({\mathbb F}_1),C({\mathbb F}_2)),
\end{equation}
where ${\mathbb F}_1$ and ${\mathbb F}_2$ are arbitrary two functors from the
category $D(\on{Loc}_{^L G})$ to itself. We will use the latter to produce
the geometric trace formula \eqref{geometric2}.

\subsection{Compatibility with Wilson/Hecke functors}

There is also another important property that the equivalence $C$ is
expected to satisfy. On both sides we have natural functors labeled by
representations $\rho$ of $^L G$. On the right hand side these are the
{\em Hecke functors} $\H_\rho$ acting from the category of ${\mc
  D}$-modules on $\Bun_G$ to the category of ${\mc D}$-modules on $X
\times \Bun_G$ (see \secref{aff Gr}).  On the left hand side these are
the functors acting from the category of $\OO$-modules on
$\on{Loc}_{^L G}$ to the category of sheaves on $X \times \on{Loc}_{^L
  G}$, which are ${\mc D}$-modules along $X$ and $\OO$-modules along
$\on{Loc}_{^L G}$. Following \cite{KW}, we will call them {\em Wilson
  functors} and denote them by $\W_\rho$.

By definition,
\begin{equation}    \label{wilson}
\W_\rho({\mc F}) = {\mc T}_\rho \otimes p_2^*({\mc F}),
\end{equation}
where ${\mc T}_\rho$ is the tautological vector bundle on $X \times
\on{Loc}_{^L G}$ defined in \secref{lfpf} and $p_2: X \times
\on{Loc}_{^L G} \to \on{Loc}_{^L G}$ is the natural projection. Note
that, by construction, ${\mc T}_\rho$ carries a connection along $X$
and so the right hand side of \eqref{wilson} is a ${\mc D}$-module
along $X$. The Wilson functor $\W_\rho$ may also be written as the
"integral transform" functor corresponding to the $\OO$-module
$\Delta_*({\mc T}_\rho)$ on $X \times \on{Loc}_{^L G} \times
\on{Loc}_{^L G}$, where $\Delta$ is the diagonal embedding of $X
\times \on{Loc}_{^L G}$:
$$
{\mc F} \mapsto q_*(q'{}^*({\mc F}) \otimes \Delta_*({\mc T}_\rho))),
$$
where $q$ and $q'$ are the projections onto $X \times \on{Loc}_{^L G}$
(the first factor) and $\on{Loc}_{^L G}$ (the second factor),
respectively.

We may also consider more general functors built from the Wilson
functors -- for example, the functor $\W_{d,\rho}$ defined by the
formula
$$
\W_{d,\rho}({\mc F}) = {\mc F}_{d,\rho} \otimes p_2^*({\mc F}),
$$
where ${\mc F}_{d,\rho}$ is an $\OO$-module on $\Loc_{^L G}$ given by
formula \eqref{Fdrho}.

The compatibility of the equivalence $C$ with the Wilson/Hecke
functors is the statement that we should have a family of
isomorphisms, which are compatible with the tensor product structures
on the Wilson and Hecke functors,
\begin{equation}    \label{functorial isom}
C(\W_\rho) \simeq \H_\rho, \qquad \rho \in \on{Rep}({}^L G).
\end{equation}
This implies that we have functorial isomorphisms
\begin{equation}    \label{WH}
C(\W_\rho({\mc F})) \simeq \H_\rho(C({\mc F})), \qquad \rho \in
\on{Rep}({}^L G).
\end{equation}

In particular, observe that the skyscraper sheaf $\OO_\sigma$ at
$\sigma \in \on{Loc}_{^L G}$ is obviously an eigensheaf of the Wilson
functors:
$$
\W_\rho(\OO_\sigma) = (\rho \circ \, \sigma) \boxtimes \OO_\sigma.
$$
Indeed, if we tensor a skyscraper sheaf with a vector
bundle, this is the same as tensoring it with the fiber of this vector
bundle at the point of support of this skyscraper sheaf. Therefore
\eqref{WH} implies that $C(\OO_\sigma)$ must be a Hecke eigensheaf on
$\Bun_G$ with ``eigenvalue'' $\sigma$ (see formula
\eqref{Hecke}). Hence we recover the traditional formulation of the
geometric Langlands correspondence.

\subsection{Application to the geometric trace
  formula}    \label{appl}

Consider now the isomorphism \eqref{RHom1} in the case ${\mathbb F}_1
= \on{Id}$ and ${\mathbb F}_2 = \W_{d,\rho}$. 
It follows from the construction of the functors $\W_{d,\rho}$ and
$\H_{d,\rho}$ and formulas \eqref{functorial isom} and \eqref{WH} that
$$
C(\W_{d,\rho}) \simeq \K_{d,\rho}.
$$
Hence \eqref{RHom1} gives us an isomorphism
$$
\on{RHom}(\on{Id},\W_{d,\rho}) \simeq \on{RHom}(\on{Id},\K_{d,\rho}).
$$

We should also have an isomorphism of the RHoms's of the kernels
defining these functors. On the right hand side this is the RHom
$$
\on{RHom}(\Delta_!(\underline{{\mc C}}),\ol{\mc K}_{d,\rho})
$$
(where $\underline{\mc C}$ is the constant sheaf on $\Bun_G$ and
$\ol{\mc K}_{d,\rho} = {\mathbf p}_*({\mc K}_{d,\rho})$) in the
derived category of ${\mc D}$-modules on $\Bun_G \times \Bun_G$, which
we have proposed in formula \eqref{rhs-new1} as a geometric analogue
of the right hand side of the trace formula.

On the left hand side of the categorical Langlands correspondence, the
kernel corresponding to $\W_{d,\rho}$ is just $\Delta_*({\mc
  F}_{d,\rho})$, where ${\mc F}_{d,\rho}$ is the $\OO$-module on
$\Loc_{^L G}$ defined by \eqref{Fdrho}. In other words, it is
supported on the diagonal in $\on{Loc}_{^L G} \times \on{Loc}_{^L G}$
and is just the $\OO$-module push-forward of ${\mc F}_{d,\rho}$ from
the diagonal.

Thus, we find that the categorical version of the geometric Langlands
correspondence should yield the following isomorphism
\begin{equation}    \label{geometric}
\on{RHom}(\Delta_*(\OO),\Delta_*({\mc F}_{d,\rho})) \simeq
\on{RHom}(\Delta_!(\underline{{\mc C}}),\ol{\mc K}_{d,\rho}).
\end{equation}

\medskip

We recall from \secref{v sp} that the right hand side is isomorphic to
$$
H^\bullet({\mc M}_d,\Delta_{\Hc}^!({\mc K}_{d,\rho})).
$$

We may also rewrite the left hand side of \eqref{geometric}
using adjunction as follows:
$$
\on{RHom}(\Delta^* \Delta_*(\OO),{\mc F}_{d,\rho}).
$$

Hence we obtain the following isomorphism, which we conjecture as a
{\em geometric trace formula}.

\begin{conj}    \label{gtf conj}
There is an isomorphism of vector spaces:
\begin{equation}    \label{geometric2}
\on{RHom}(\Delta^* \Delta_*(\OO),{\mc F}_{d,\rho}) \simeq
H^\bullet(\Bun_G,\Delta^!(\ol{\mc K}_{d,\rho})).
\end{equation}
\end{conj}

\medskip

Thus, starting with the categorical Langlands correspondence, we have
arrived at what we propose as a geometrization of the trace formula.

We have a similar conjecture for more general functors $\K$ built from
the Hecke functors of the type introduced in \secref{int tr}. On the
other side of the categorical Langlands correspondence it corresponds
to a functor $\W$ built in the same way from the Wilson functors, so
that
$$
C(\K) \simeq \W.
$$
We then expect to have an analogous isomorphism
\begin{equation}    \label{more general}
\on{RHom}(\Delta^* \Delta_*(\OO),{\mc F}) \simeq
H^\bullet(\Bun_G,\Delta^!(\ol{\mc K})).
\end{equation}
where $\ol{\mc K}$ is the kernel of $\K$ and ${\mc F}$ is the
$\OO$-module on $\Loc_{^L G}$ such that $\Delta_*({\mc F})$ is the
kernel of $\W$.

\begin{remark}
It is tempting to rewrite the left hand side of \eqref{more general}
as
$$
\on{RHom}(\OO,\Delta^! \Delta_*({\mc F})) = H^\bullet(\Loc_{^L
  G},\Delta^! \Delta_*({\mc F})),
$$
where $\Delta^!$ is the right adjoint functor to $\Delta_*$. Then
\eqref{more general} becomes more symmetrical:
\begin{equation}    \label{more general1}
H^\bullet(\Loc_{^L G},\Delta^! \Delta_*({\mc F})) \simeq
H^\bullet(\Bun_G,\Delta^!(\ol{\mc K})).
\end{equation}
However, the definition of the functor $\Delta^!$ for the derived
category of coherent sheaves on the stack $\Loc_{^L G}$ is tricky.
Defining the left hand side of \eqref{more general1} essentially
amounts to rewriting it as the left hand side of \eqref{more
  general}. That is why we stated the conjecture in the form
\eqref{more general}.\qed
\end{remark}

\medskip

As a special case, let both $\K$ and $\W$ be the identity
functors. Then \eqref{more general} gives rise to an isomorphism
\begin{equation}    \label{trivial case}
\on{RHom}(\Delta^* \Delta_*(\OO),\OO) \simeq
H^\bullet_{\on{dR}}(\Bun_G).
\end{equation}
The right hand side of this formula is known \cite{Tel,NS,HS}, and
hence this special case already gives us a good test for the
conjecture. The isomorphism \eqref{trivial case} holds in the abelian
case $G = {\mathbb G}_m$, as we will see at the end of in \secref{ex}.

\subsection{Connection to the Atiyah--Bott-Lefschetz fixed point
  formula}    \label{abl fixed}

Here we give an heuristic explanation why we should think of the space
\eqref{more general} as a geometric incarnation of the sum appearing
on the left hand side of \eqref{classical}.

Let ${\mc F}$ be a coherent sheaf on $\on{Loc}_{^L G}$ built from the
vector bundles ${\mc T}_{\rho,x}, x \in |X|$ (like ${\mc F}_{d,\rho}$
constructed below). We would like to interpret taking the trace of the
Frobenius on the (coherent!) cohomology $H^\bullet(\on{Loc}_{^L
  G},{\mc F})$ as the sum over points $\sigma$ of $\on{Loc}_{^L G}$,
which we think of as the fixed points of the Frobenius automorphism
acting on a moduli of homomorphisms
$$
\pi_1(X \underset{k}\otimes \ol{k}) \to {}^L G.
$$

Recall the Atiyah--Bott--Lefshetz fixed point formula \cite{AB} (see
also \cite{Illusie}, \S 6). Let $M$ be a smooth proper scheme, $V$ a
vector bundle on $M$, and $\V$ the coherent sheaf of sections of
$V$. Let $u$ be an automorphism acting on $M$ with isolated fixed
points, and suppose that we have an isomorphism $\gamma: u^*(\V) \simeq
\V$. Then
\begin{equation}    \label{l factor}
\on{Tr}(\gamma,H^\bullet(M,\V)) = \sum_{p \in M^u}
\frac{\on{Tr}(\gamma,\V_p)}{\on{det}(1-\gamma,T^*_p M)},
\end{equation}
where
$$
M^u = \Gamma_u \underset{M \times M}\times \Delta
$$
is the set of fixed points of $u$, the fiber product of the graph
$\Gamma_u$ of $u$ and the diagonal $\Delta$ in $M \times M$, which we
assume to be transversal to each other (as always, the left hand side
of \eqref{l factor} stands for the {\em alternating} sum of traces on
the cohomologies).

Now, if we take the cohomology not of $\V$, but of the tensor
product $\V \otimes \Omega^\bullet(M)$, where $\Omega^\bullet(M)
= \Lambda^\bullet(T^*(M))$ is the graded space of differential forms,
then the determinants in the denominators on the right hand side of
formula \eqref{l factor} will get cancelled and we will obtain the
following formula:
\begin{equation}    \label{l factor1}
\on{Tr}(\gamma,H^\bullet(M,\V \otimes \Omega^\bullet(M))) = \sum_{p \in M^u}
\on{Tr}(\gamma,\V_p).
\end{equation}

We would like to apply formula \eqref{l factor1} in our
situation. However, $\on{Loc}_{^L G}$ is not a scheme, but an
algebraic stack (unless $^L G$ is a torus), so we need an analogue of
\eqref{l factor1} for algebraic stacks (and more generally, for
derived algebraic stacks, since $\on{Loc}_{^L G}$ is not smooth as an
ordinary stack).

Let $M$ be as above and $\Delta: M \hookrightarrow M^2$ the diagonal
embedding. Observe that
\begin{equation}    \label{diagonal}
\Delta^* \Delta_*(\V) \simeq \V \otimes \Omega^\bullet(M),
\end{equation}
and hence we can rewrite \eqref{l factor1} as follows:
\begin{equation}    \label{l factor2}
\on{Tr}(\gamma,H^\bullet(M,\Delta^* \Delta_*(\V))) = \sum_{p \in M^u}
\on{Tr}(\gamma,\V_p).
\end{equation}

Now we propose formula \eqref{l factor2} as a {\em conjectural
  generalization of the Atiyah--Bott--Lefschetz fixed point formula}
to the case that $M$ is a smooth algebraic stack, and more generally,
smooth derived algebraic stack (provided that both sides are
well-defined). We may also allow here $\V$ to be a perfect complex (as
in \cite{Illusie}).

Even more generally, we drop the assumption that $u$ has fixed points
(or that the graph $\Gamma_u$ of $u$ and the diagonal $\Delta$ are
transversal) and conjecture the following general fixed point formula
of Atiyah--Bott type for (derived) algebraic stacks.

\begin{conj}    \label{conj abl formula}
\begin{equation}    \label{l factor3}
\on{Tr}(\gamma,H^\bullet(M,\Delta^* \Delta_*(\V))) =
\on{Tr}(\gamma,H^\bullet(M^u,i_u^* \Delta_*(\V))),
\end{equation}
where $M^u = \Gamma_u \underset{M \times M}\times \Delta$
is the fixed locus of $u$ and $i_u: M^u \to M \times M$.
\end{conj}

We want to apply \eqref{l factor2} to the left hand side of
\eqref{more general}, which is
$$
\on{RHom}(\Delta^* \Delta_*(\OO),{\mc F}).
$$
This is not exactly in the form of the left hand side of \eqref{l
  factor2}, but it is very close. Indeed, if $M$ is a smooth scheme,
then
$$
\on{RHom}(\Delta^* \Delta_*(\OO),{\mc F}) \simeq {\mc F} \otimes
\Lambda^\bullet(TM),
$$
so we obtain the exterior algebra of the tangent bundle instead of the
exterior algebra of the cotangent bundle. In our setting, we want
$\gamma$ to be the Frobenius, and so a fixed point is a homomorphism
$\pi_1(X) \to {}^L G$. The tangent space to $\sigma$ (in
the derived sense) should then be identified with the cohomology
$H^\bullet(X,\on{ad} \circ \, \sigma)[1]$, and the cotangent space with
its dual. Hence the trace of the Frobenius on the exterior algebra of
the tangent space at $\sigma$ is the $L$-function
$L(\sigma,\on{ad},s)$ evaluated at $s=0$. Poincar\'e duality implies
that the trace of the Frobenius on the exterior algebra of the
cotangent space at $\sigma$ is
$$
L(\sigma,\on{ad},1) = q^{-d_G} L(\sigma,\on{ad},0),
$$
so the ratio between the traces on the exterioir algebras of the
tangent and cotangent bundles at the fixed points results in an
overall factor which is a power of $q$.\footnote{The exact power of
  $q$ depends on whether we choose the geometric or the arithmetic
  Frobenius here. We recall from \secref{v sp} that the right hand
  side of the trace formula \eqref{classical} is equal to the trace of
  the Frobenius also up to an overall factor which is a power of
  $q$. Since this is an heuristic argument, we ignore these factors.}

Hence, by following this argument and switching from $\C$ to $\Fq$, we
obtain (up to a power of $q$) the trace formula \eqref{classical} from
the isomorphism \eqref{more general}.

\medskip

\begin{remark} The sheaf $\Delta^* \Delta_*({\mc O}_M)$ is the
  structure sheaf of the self-intersection of the diagonal of a
  derived stack $M$. In \cite{BN} it is interpreted as the loop space
  of $M$. Hence the left hand side of our conjectural generalization
  \eqref{l factor3} of the Atiyah--Bott--Lefschetz formula may be
  viewed as the cohomology of the pullback of $\V$ from $M$ to its loop
  space. David Nadler has outlined to us a possible proof of \eqref{l
    factor2} using the methods of \cite{BFN,BN}.\qed
\end{remark}

\subsection{Discussion}    \label{discussion}

The isomorphism \eqref{more general} is still tentative, because there
are some unresolved issues in the definition of the two sides of the
isomorphism \eqref{more general}.

The first issue is the definition of the stack $\Loc_{^L
  G}$. Presumably, it should be defined not as a stack, but as a
derived (or differential graded, or DG) stack. This should correspond
to the inclusion of the non-tempered (or non-Ramanujan) automorphic
representations, which should be parametrized, according to Arthur, by
homomorphisms
$$
\sigma: W_F \times SL_2 \to {}^L G
$$
(for tempered representations, $\sigma|_{SL_2}$ is expected to be
trivial). For instance, the trivial representation of $G(\AD)$
corresponds to $\sigma$ that is trivial on $W_F$ and is a principal
$SL_2$-embedding on the second factor. Recent results of V. Lafforgue
\cite{Laf} show that one might be able to take such $\sigma$ into
account by considering a derived version of the stack $\Loc_{^L G}$. There
are also additional insights pointing in the same directions coming
from the study of $S$-duality in physics (see a brief discussion at
the end of \cite{F:bourbaki} and references therein).

The second issue is the precise definition of $\on{RHom}$'s in the two
categories involved in the Langlands correspondence \eqref{na fm},
which appear on the two sides of the isomorphism \eqref{more
  general}. This is in fact the main question in finding the correct
formulation of the categorical Langlands correspondence.

Finally, as we have discussed in Part I, in order to become an
effective tool, trace formulas have to be stabilized. Hence an
analogue of stabilization needs to be worked out in the geometric
setting as well.

\medskip

But suppose we could solve these problems and really make sense of
\eqref{more general} as a geometric trace formula. What would we be able
to learn from it?

First of all, we would have a good framework for the geometric trace
formula that works for curves over $\Fq$ and over $\C$.
Second, we could try to establish the isomorphism
\eqref{more general} by using fixed-point formulas,
for example. Here the description of ${\mathbb G}_m$-fixed points in
the moduli space of stable flat bundles due to Hitchin
\cite{Hit1,Hit2} and Simpson \cite{Simp} could be useful.

\medskip

The main application we have in mind is to use the isomorphism
\eqref{more general} to establish functoriality. Namely, we wish to
express each side as a direct sum of two vectors spaces corresponding
to the tempered (Ramanujan) and non-tempered (non-Ramanujan)
representations, and then to decompose further the former as a direct
sum over the groups $H$ (see \secref{decomp trace}). In
\secref{orbital} we made the first steps toward this goal for the
right hand (orbital) side of \eqref{geom trace}, using an analogue of
the Hitchin fibration. The fact that we now have the left hand
(spectral) side of the isomorphism \eqref{more general} may help us
understand better how to make this decomposition. For instance, the
locus of ${\mathbb G}_m$-fixed points in the moduli space of stable
$^L G$-local systems decomposes into the union of the $\phi=0$ locus
(where $\phi$ is the Higgs field) and other loci with $\phi \neq
0$. It is possible that the former corresponds to the
tempered part and the latter to the non-tempered part.

\section{Relative geometric trace formula}    \label{rel trace}

In the previous section we discussed a geometrization of the trace
formula \eqref{tr}. We have seen in the previous section that a
geometric analogue of this formula may be constructed within the
framework of the categorical form of the Langlands correspondence, as
the statement that the RHom's of kernels of certain natural functors
are isomorphic. These kernels are sheaves on algebraic stacks over the
squares $\Bun_G \times \Bun_G$ and $\on{Loc}_{^L G} \times
\on{Loc}_{^L G}$.

It is natural to ask what kind of statement we may obtain if we
consider instead the RHom's of sheaves on the stacks $\Bun_G$ and
$\on{Loc}_{^L G}$ themselves.

In this section we will show that this way we obtain what may be
viewed as a geometric analogue of the so-called {\em relative trace
  formula}. On the spectral side of this formula we also have a sum
like \eqref{trace1}, but with one important modification; namely, the
eigenvalues $N_{\sigma}$ are weighted with the factor
$L(\sigma,\on{ad},1)^{-1}$.  The insertion of this factor was
originally suggested by Sarnak in \cite{Sarnak} and further studied by
Venkatesh \cite{Ven}, for the group $GL_2$ in the number field
context. Another difference is that the summation is restricted to the
generic representations, which means (at least, conjecturally) that we
remove the contribution of the non-Ramanujan representations as well
as the multiplicity factors $m_\sigma$.

We will also analyze the corresponding orbital side and express is as
the cohomology of a sheaf on an algebraic stack. This will give us
what we may call a {\em relative geometric trace formula}. It will be
the statement about an isomorphism of two cohomology spaces, one on
$\Bun_G$ and the other on $\Loc_{^L G}$. We will see that this
isomorphism arises naturally from the categorical form of the
geometric Langlands correspondence \eqref{na fm}.

\subsection{Relative trace formula}

First, let us recall the setup of the relative trace formula.

Let $G$ be a split simple algebraic group over $k=\Fq$.
In order to state the relative trace formula, we need
to choose a non-degenerate character of $N(F)\bs N(\AD)$, where $N$ is
a maximal unipotent subgroup of $G$. A convenient way to define it is
to consider a twist of the group $G$. Let us pick a maximal torus $T$
such that $B=TN$ is a Borel subgroup. If the maximal torus $T$ admits
the cocharacter $\crho: \Gm \to T$ equal to half-sum of all positive
roots (corresponding to $B$), then let $K_X^{\crho}$ be the $T$-bundle
on our curve $X$ which is the pushout of the $\Gm$-bundle $K_X^\times$
(the canonical line bundle on $X$ without the zero section) under
$\crho$. If $\crho$ is not a cocharacter of $T$, then its square is,
and hence this $T$-bundle is well-defined for each choice of the
square root $K_X^{1/2}$ of $K_X$. We will make that choice once and for
all.\footnote{This choice is related to the ambiguity of the
equivalence \eqref{na fm}, see the footnote on page \pageref{foot}.} 

Now set
$$
G^K = K_X^{\crho} \underset{T}\times G, \qquad N^K =
K_X^{\crho} \underset{T}\times N,
$$
where $T$ acts via the adjoint action. For instance, if $G=GL_n$, then
$GL_n^K$ is the group scheme of automorphism of the rank $n$
bundle ${\mc O} \oplus K_X \oplus \ldots \oplus
K_X^{\otimes(n-1)}$ on $X$ (rather than the trivial bundle ${\mc
O}^{\oplus n}$).

We have
\begin{equation}    \label{comm}
N^K/[N^K,N^K] = K_X^{\bigoplus \on{rank}(G)}.
\end{equation}
Now let $\psi: \Fq \to \C^\times$ be an additive
character, and define a character $\Psi$ of $N^K(\AD)$ as follows
$$
\Psi((u_x)_{x \in |X|}) = \prod_{x \in |X|} \prod_{i=1}^{\on{rank} G}
\psi(\on{Tr}_{k_x/k} \on{Res}_x(u_{x,i})),
$$
where $u_{x,i} \in K_X(F_x)$ is the $i$th projection of $u_x \in
N^K(F_x)$ onto $K_X(F_x)$ via the isomorphism \eqref{comm}. We
denote by $k_x$ the residue field of $x$, which is a finite extension
of the ground field $k=\Fq$.

By the residue formula, $\Psi$ is trivial on the subgroup $N^K(F)$
(this was the reason why we introduced the twist). It is also trivial
on $N^K({\mc O})$.

In what follows, in order to simplify notation, we will denote
$G^K$ and $N^K$ simply by $G$ and $N$.

Given an automorphic representation $\pi$ of $G(\AD)$, we have the
Whittaker functional $W: \pi \to \C$,
$$
W(f) = \int_{N(F)\bs N(\AD)} f(u) \Psi^{-1}(u) du,
$$
where $du$ is the Haar measure on $N(\AD)$ normalized so that the
volume of $N(\OO)$ is equal to $1$.

We choose, for each unramified automorphic representation $\pi$, a
non-zero $G(\OO)$-invariant function $f_\pi \in \pi$ on $G(F)\bs
G(\AD)$.

Let $K$ again be a kernel on the square of $\Bun_G(k) =
G(F) \bs G(\AD)/G(\OO)$ and ${\mb K}$ the
corresponding integral operator acting on unramified automorphic
functions. The simplest unramified version of the relative trace
formula reads (here we restrict the summation to cuspidal automorphic
representations $\pi$)
\begin{multline}    \label{simplest}
\sum_{\pi} \ol{W_{\Psi}(f_{\pi})} \; W_\Psi(K \cdot f_\pi) ||f_{\pi}||^{-2} = \\
\underset{N(F) \bs N(\AD)}\int \quad \underset{N(F) \bs N(\AD)}\int
K(u_1,u_2) \Psi^{-1}(u_1) \Psi(u_2) du_1 du_2,
\end{multline}
(see, e.g., \cite{Jacquet}), where
$$
||f||^2 = \underset{G(F)\bs G(\AD)}\int |f(g)|^2 dg,
$$
and $dg$ denotes the invariant Haar measure normalized so that the
volume of $G(\OO)$ is equal to $1$. Note that
$$
||f||^2 = q^{d_G} L(G) ||f||^2_T,
$$
where $||f||^2_T$ is the norm corresponding to the Tamagawa measure,
$d_G = (g-1) \dim G$,
$$
L(G) = \prod_{i=1}^\ell \zeta(m_i+1),
$$
where the $m_i$ are the exponents of $G$.

The following conjecture was communicated to us by B. Gross and
A. Ichino. In the case of $G=SL_n$ or $PGL_n$, formula \eqref{rs1}
follows from the Rankin--Selberg convolution formulas, as we explain
in \secref{rs sec} below. Other cases have been considered in
\cite{GP,Ich1,Ich2}. For a general semi-simple group $G$ of adjoint
type formula \eqref{rs1} has been conjectured by A. Ichino and
T. Ikeda assuming that $\pi$ is square-integrable. Note that if a
square-integrable representation is tempered, then it is expected to
be cuspidal, and that is why formula \eqref{rs1} is stated only for
cuspidal representations.

Recall that an $L$-packet of automorphic representations is called
{\em generic} if each irreducible representation $\pi = \bigotimes'
\pi_x$ from this $L$-packet has the property that the local $L$-packet
of $\pi_x$ contains a generic representation (with respect to a
particular choice of non-degenerate character of $N(F_x)$).

\begin{conj}    \label{gen}
Suppose that the $L$-packet corresponding to an unramified $\sigma:
W_F \to {}^L G$ is generic. Then it contains a unique, up to an
isomorphism, irreducible representation $\pi$ such that $W_\Psi(f_\pi)
\neq 0$, with multiplicity $m_\pi = 1$. Moreover, if this $\pi$ is
in addition cuspidal, then the following formula holds:
\begin{equation}    \label{rs1}
|W_\Psi(f_\pi)|^2 ||f_\pi||^{-2} = q^{d_N-d_G} L(\sigma,\on{ad},1)^{-1}
|S_\sigma|^{-1},
\end{equation}
where $S_\sigma$ is the (finite) centralizer of the image of $\sigma$
in $^L G$, $d_N = -(g-1)(4 \langle \rho,\crho \rangle - \dim N)$ (see
formula \eqref{dim BunN}), and $d_G = (g-1) \dim G$.
\end{conj}

Furthermore, we expect that if the $L$-packet corresponding to
$\sigma: W_F \times SL_2 \to {}^L G$ is non-generic, then
\begin{equation}    \label{non-gen}
L(\sigma,\on{ad},1)^{-1} = 0.
\end{equation}
Ichino has explained to us that according to Arthur's conjectures,
square-integrable non-tempered representations are non-generic.

Recall that we have ${\mb K} \cdot f_\pi = N_\sigma f_\pi$. Therefore
\conjref{gen} and formula \eqref{simplest} give us the following:
\begin{multline}    \label{relative}
q^{-d_G} \sum_{\sigma: W_F \to {}^L G} N_\sigma \cdot
L(\sigma,\on{ad},1)^{-1} |S_\sigma|^{-1} =
\\ q^{-d_N} \int \int K(u_1,u_2) \Psi^{-1}(u_1) \Psi(u_2) du_1 du_2.
\end{multline}
On the left hand side we sum only over unramified $\sigma$, and only
those of them contribute for which the corresponding $L$-packet of
automorphic representations $\pi$ is generic.

We expect that the left hand side of formula \eqref{relative} has the
following properties:

\begin{itemize}

\item[(1)] It does not include homomorphisms $\sigma: SL_2 \times W_F
  \to \LG$ which are non-trivial on the Artur's $SL_2$.

\item[(2)] The multiplicity factor $m_\sigma$ of formula
  \eqref{classical} disappears, because only one irreducible
  representation from the $L$-packet corresponding to $\sigma$ shows
  up (with multiplicity one).

\item[(3)] Since $S_\sigma$ is the group of automorphisms of $\sigma$,
  the factor $|S_\sigma|^{-1}$ makes the sum on the left hand side
  \eqref{relative} look like the Lefschetz fixed point formula for
  stacks.

\end{itemize}

\subsection{The case of $G=PGL_n$}    \label{rs sec}

Here we derive formula \eqref{rs1} for $G=PGL_n$ using the
Rankin-Selberg convolution method, following \cite{Lys}. Denote by
$\Bun_n$ the moduli stack of rank $n$ bundles on $X$ and $\Bun^m_n$
the connected component corresponding to vector bundles of degree $m$.
Let $\sigma$ be an irreducible $n$-dimensional $\ell$-adic
representation of $W_F$ and ${\mc F}_\sigma$ the pure perverse sheaf
of weight 0 on $\Bun_{n}$ (irreducible on each $\Bun^m_n$) which is a
Hecke eigensheaf with respect to $\sigma$ (see \cite{FGV:lc}).  Let
$f_\sigma$ be the corresponding function on $\Bun_{GL_n}(\Fq)$. The
following formula is derived in \cite{Lys}:
$$
\sum_{L \in \Bun^m_n(\Fq)} \frac{1}{|\on{Aut}(L)|} f_\sigma(L)
f_{\sigma^*}(L) = \on{Res}_{s=1} L(\sigma,\on{ad},s)^{-1}
$$
(here $\on{ad}$ denotes the adjoint representation of $GL_n$). In
fact, S. Lysenko \cite{Lys} has given a geometric interpretation of
this formula, which is compatible with the categorical Langlands
correspondence.

In addition, for irreducible $\sigma$ we have
$$
W_\Psi(f_\sigma) = q^{d_N-d_G},
$$
where
$$
d_N = -\frac{(g-1)n(n-1)(2n-1)}{6}, \qquad d_G = (g-1)n^2,
$$
see \cite{LL}.

Using this formula, it is straightforward to obtain for $G=\PGL_n$
that
$$
|W_\Psi(f_\pi)|^2 ||f_\pi||^{-2} = \frac{1}{n} q^{d_N-d_G}
L(\sigma,\on{ad},1)^{-1},
$$
where $\on{ad}$ denotes the adjoint representation of $^L G =
\SL_n$. This agrees with formula \eqref{rs1}, because $S_\sigma$ is
the center of $^L G$ in this case (here we are again under the
assumption that $\sigma$ is irreducible).

\subsection{Geometric meaning: right hand side}

Now we discuss the geometric meaning of formula \eqref{relative},
starting with the right hand side. Let $\Bun^{\F_T}_N$ be the moduli
stack of $B=B^K$ bundles on $X$ such that the corresponding $T$-bundle
is $\F_T = K_X^{\crho}$. Note that
\begin{equation}    \label{dim BunN}
\dim \Bun^{\F_T}_N = d_N = -(g-1)(4 \langle \rho,\crho \rangle - \dim N).
\end{equation}

Let $\on{ev}: \Bun^{\F_T}_N \to {\mathbb G}_a$ be the map
constructed in \cite{FGV:w}.

For instance, if $G=GL_2$, then $\Bun^{\F_T}_N$ classifies
rank two vector bundles $\V$ on $X$ which fit in the exact sequence
\begin{equation}    \label{exten}
0 \to K_X^{1/2} \to \V \to K_X^{-1/2} \to 0,
\end{equation}
where $K_X^{1/2}$ is a square root of $K_X$ which we have fixed.
The map $\on{ev}$ assigns to such $\V$ its extension class in
$\on{Ext}(\OO_X,K_X) = H^1(X,K_X) \simeq {\mathbb G}_a$. For other
groups the construction is similar (see \cite{FGV:w}).

On ${\mathbb G}_a$ we have the Artin-Schreier sheaf ${\mc L}_\psi$
associated to the additive character $\psi$. We define the sheaf
$$
\wt\Psi = \on{ev}^*({\mc L}|_\psi)
$$
on $\Bun^{\F_T}_N$. Next, let $p: \Bun^{\F_T}_N \to \Bun_G$ be the
natural morphism. Let
$$
\Psi = p_!(\wt\Psi)[d_N-d_G]((d_N-d_G)/2).
$$

Then the right hand side of \eqref{relative} is equal to the trace of
the Frobenius on the vector space
\begin{equation}    \label{lhs1}
\on{RHom}(\Psi,{\mathbb K}(\Psi)).
\end{equation}
Here we use the fact that ${\mathbb D} \circ \K \simeq \K \circ
{\mathbb D}$ and ${\mathbb D}(\wt\Psi) \simeq \on{ev}^*({\mc
  L}|_{\psi^{-1}})[2d_N](d_N)$.

\subsection{Geometric meaning: left hand side}    \label{geom}

As discussed above, we don't have an algebraic stack parametrizing
homomorphisms $\sigma: W_F \to {}^L G$ if our curve $X$ is defined
over a finite field $\Fq$. But such a stack exists when $X$ is over
$\C$, though in this case there is no Frobenius operator on the
cohomology whose trace would yield the desired number (the left
hand side of \eqref{relative}). In this subsection we will define a
certain vector space (when $X$ is over $\C$) and conjecture that it is
isomorphic to the vector space \eqref{lhs} which is the geometric
incarnation of the right hand side of \eqref{relative} (and which is
well-defined for $X$ over both $\Fq$ and $\C$). This will be our
``relative geometric trace formula''. In the next subsection we will
show that this isomorphism is a corollary of the categorical version
of the geometric Langlands correspondence.

In order to define this vector space, we will use the coherent sheaf
${\mc F}_{d,\rho}$ on $\on{Loc}_{^L G}$ introduced in formula
\eqref{Fdrho}.  We propose that the geometric incarnation of the left
hand side of \eqref{relative} in the case when $\K = \K_{d,\rho}$ is
the cohomology
\begin{equation}   \label{rhs}
H^\bullet(\on{Loc}_{^L G},{\mc F}_{d,\rho}).
\end{equation}

The heuristic explanation for this proceeds along the lines of the
explanation given in the case of the ordinary trace formula in
\secref{abl fixed}, using the Atiyah--Bott--Lefschetz fixed point
formula.\footnote{We note that applications of the
  Atiyah--Bott--Lefschetz fixed point formula in the context of Galois
  representations have been previously considered by M. Kontsevich in
  \cite{Ko}.} We wish to apply it to the cohomology \eqref{rhs}. If
$\Loc_{^L G}$ were a smooth scheme, then we would have to multiply the
number $N_\sigma$ which corresponds to the stalk of ${\mc F}_{d,\rho}$
at $\sigma$, by the factor
\begin{equation}    \label{weighting factor}
\on{det}(1-\on{Fr},T^*_\sigma \Loc_{^L G})^{-1}.
\end{equation}
Recall that the tangent space to $\sigma$ (in the derived sense) may
be identified with the cohomology $H^\bullet(X,\on{ad} \circ \,
\sigma)[1]$. Using the Poincar\'e duality, we find that the factor
\eqref{weighting factor} is equal to
$$
L(\sigma,\on{ad},1)^{-1}.
$$

Therefore, if we could apply the Lefschetz fixed point formula to the
cohomology \eqref{rhs} and write it as a sum over all $\sigma: W_F \to
\LG$, then the result would be the left hand side of \eqref{relative}
(up to a factor that is a power of $q$). (Note however that since
$\Loc_{^L G}$ is not a scheme, but an algebraic stack, the weighting
factor should be more complicated for those $\sigma$ which admit
non-trivial automorphisms, see the conjectural fixed point formula
\eqref{l factor3} in \secref{abl fixed}.)

This leads us to the following {\em relative geometric trace formula}
(in the case of the functor $\K_{d,\rho}$).

\begin{conj}    \label{rel conj}
We have the following isomorphism of vector spaces:
\begin{equation}    \label{geom relative}
H^\bullet(\on{Loc}_{^L G},{\mc F}_{d,\rho}) \simeq
\on{RHom}_{\on{Bun_G}}(\Psi,{\mathbb K}_{d,\rho}(\Psi)).
\end{equation}
\end{conj}

It is clear how to generalize this to other functors $\K$ of the form
described in \secref{int tr}.

Now we explain how the isomorphism \eqref{geom relative}
fits in the framework of a
categorical version of the geometric Langlands correspondence.

\subsection{Interpretation from the point of view of the categorical
  Langlands correspondence}

We start by asking what is the ${\mc D}$-module on $\Bun_G$
corresponding to the structure sheaf $\OO$ on $\on{Loc}_{^L G}$ under
the functor $C$ discussed in \secref{clc}. The
following answer was suggested by Drinfeld (see \cite{Laf}): it is the
sheaf $\Psi$ that we have discussed above.

The rationale for this proposal is the following: we have
$$
\on{RHom}_{\on{Loc}_{^L G}}(\OO,\OO_\sigma) = \C, \qquad \forall
\sigma,
$$
where $\OO_\sigma$ is again the skyscraper sheaf supported at
$\sigma$. Therefore, since $C(\OO_\sigma) = {\mc F}_\sigma$, we should
have, according to \eqref{RHom},
$$
\on{RHom}_{\Bun_G}(C(\OO),{\mc F}_\sigma) = \C, \qquad \forall
\sigma.
$$
According to the conjecture of \cite{LL}, the sheaf $\Psi$ has just
this property:
$$
\on{RHom}_{\Bun_G}(\Psi,{\mc F}_\sigma)
$$
is the one-dimensional vector space in cohomological degree $0$ (if we
use appropriate normalization for ${\mc F}_\sigma$).

This vector space should be viewed as a geometric incarnation of the Fourier coefficient of the
automorphic function corresponding to ${\mc F}_\sigma$.

This provides
some justification for the assertion that\footnote{As explained in the
footnote on page \pageref{foot}, $C$ corresponds to a particular
choice of $K_X^{1/2}$. Given such a choice, $C(\OO)$ should be the
character sheaf $\Psi$ associated to that $K_X^{1/2}$.}
\begin{equation}    \label{CO}
C(\OO) = \Psi.
\end{equation}
Let us take this for granted.

Now let us look at the left hand side of the isomorphism \eqref{geom
  relative}. We can rewrite it as the de Rham cohomology of a sheaf on
  $X_d$, whose fiber at $D = \sum n_i [x_i] \in X_d$ is
\begin{equation}    \label{rhs1}
\on{RHom}_{\on{Loc}_{^L G}}\left(\OO,\bigotimes_i {\mc
T}^{(n_i)}_{\rho,x_i} \right) = \on{RHom}\left(\OO,\prod_i
\W^{(n_i)}_{\rho,x_i}(\OO) \right),
\end{equation}
where $\W^{(n_i)}_{\rho,x_i}$ is the Wilson operator corresponding to
the representation $\on{Sym}^{n_i} \rho$, specialized to the point
$x_i \in X$ (see formula \eqref{wilson}).

Using the compatibility \eqref{WH} with the Wilson/Hecke
operators and formulas \eqref{RHom} and \eqref{CO},
we obtain that \eqref{rhs1} should be isomorphic to
\begin{equation}    \label{lhs}
\on{RHom}_{\Bun_G}\left(\Psi,\prod_i \H^{(n_i)}_{\rho,x_i}(\Psi) \right),
\end{equation}
where $\H^{(n_i)}_{\rho,x_i}$ is the Hecke operator corresponding to the
representation $\on{Sym}^{n_i} \rho$, specialized to the point
$x_i \in X$. Hence the right hand side of \eqref{geom relative} is
isomorphic to the de Rham cohomology of a sheaf on $X_d$ whose
fiber at $D = \sum n_i [x_i] \in X_d$ is \eqref{lhs}. But varying
$D$ over $X_d$ gives us precisely our averaging functor ${\mathbb
K}_{d,\rho}$.  Therefore the result is
\begin{equation}    \label{lhs2}
\on{RHom}_{\Bun_G}\left(\Psi,{\mathbb K}_{d,\rho}(\Psi) \right),
\end{equation}
which is the right hand side of \eqref{geom relative}.

Thus, we obtain that the relative geometric trace formula \eqref{geom
  relative} follows from the categorical version of the geometric
  Langlands correspondence.

\medskip

Let's recap. We started out with the classical relative trace formula
\eqref{relative}. The right hand side of \eqref{relative} had a simple
geometric interpretation as the trace of the Frobenius on the vector
space \eqref{lhs}. To interpret the left hand side of \eqref{relative}
we made a ``leap of faith'' and replaced it by the cohomology
\eqref{rhs}. However, we have just shown that the resulting geometric
relative trace formula \eqref{geom relative} is a meaningful statement
from the point of view of the categorical geometric Langlands
correspondence. This is an indication that \eqref{geom relative} is
actually a reasonable conjecture (and hence so is our ``leap of
faith'').

\subsection{Applications to functoriality}

As explained in \cite{FLN} (following \cite{L:BE,L:PR,L:Shaw}), our
goal is to express the trace of the averaging operator ${\mb
  K}_{d,\rho}$ for the group $G$ as the sum of traces for the groups
$H$ dual to subgroups $^\la H \subset {}^L G$ for which representation
$\rho$ contains non-zero invariant vectors. This has important
applications to functoriality.

Likewise, in the setting of the relative trace formula, we are
interested in decomposing the left hand side as a sum over the groups $H$.
It is instructive to see why and how this decomposition should come
about, first for an arbitrary $G$ and then in the example of
$G=GL_2$ and $\rho$ the second symmetric power of its defining representation.

\bigskip

Let us consider the general case first. Look at the left hand side of
formula \eqref{relative}. According to Corollary \ref{van} (see
also \cite{FLN}, Lemma 2.6), $N_\sigma = 0$ for $d>(2g-2)\dim \rho$,
unless $(\rho \circ \, \sigma)$ contains non-zero invariant
vectors. From now on we will assume that $d>(2g-2)\dim \rho$.  Then we
obtain that the left hand side of formula \eqref{relative} decomposes
into a sum
\begin{equation}    \label{sum}
\sum_{^\lambda H \subset {}^L G} \quad \sum_{\sigma': W_F \to
  {}^\lambda H} N_{\sigma'} \cdot L(\sigma',\on{ad}_{^L
  G},q^{-1})^{-1} |S_\sigma|^{-1}
\end{equation}
over all possible $^\lambda H \subset {}^L G$ such that
$^\lambda H$ has a trivial subrepresentation in $\rho \circ
\sigma$. Then each $\sigma': W_F \to {}^\lambda H$ gives rise to
$\sigma: W_F \to {}^L G$ and to $\rho \circ \, \sigma': W_F \to
\on{Aut} \rho$ (we view $\rho$ as a representation of $^\la H$
obtained by restriction from $^L G$). We have
\begin{align}
N_{\sigma'} &= \on{Tr}\left(\on{Fr},\sum_i (-1)^i
H^i(X_d,(\rho \circ \, \sigma')_d)\right) \\ \notag
&=q^{-ds}\on{-coeff.} \on{of} L(\sigma',\rho,s).
\end{align}

Observe that $^L \g = {}^\la \h \oplus R_{H}$, where ${}^\la \h$
is the Lie algebra of $^\la H$ and $R_{H}$ is a certain
representation of $^\la H$. We have
$$
L(\sigma,\on{ad}_{^L G},1) = L(\sigma',\on{ad}_{^\la
  H},1) \cdot L(\sigma',R_{H},1).
$$
Hence the term in the sum \eqref{sum} corresponding to a
particular $^\la H$ reads
\begin{equation}    \label{individual}
\sum_{\sigma': W_F \to {}^\lambda H} N_{\sigma'} \cdot
L(\sigma',\on{ad}_{^\la H},1)^{-1} \cdot
L(\sigma',R_{H},1)^{-1} \cdot |S_\sigma|^{-1}.
\end{equation}
The difference between \eqref{individual} and the left hand side of
\eqref{relative} is that we have the additional weighting factor
\begin{equation}    \label{extra}
L(\sigma',R_{H},1)^{-1}
\end{equation}
(note that it corresponds to the conormal bundle to $\on{Loc}_{^\la
  H}$ in $\on{Loc}_{^L G}$).

Thus, the formula which we want to prove is
\begin{multline}    \label{want}
\sum_\sigma N_\sigma \cdot L(\sigma,\on{ad},1)^{-1}
|S_\sigma|^{-1} = \\ \sum_{^\lambda H \subset {}^L G} \; \;
\sum_{\sigma': W_F \to {}^\lambda H} N_{\sigma'} \cdot
L(\sigma',\on{ad}_{^\la H},1)^{-1} \cdot
L(\sigma',R_{H},1)^{-1} |S_\sigma|^{-1}.
\end{multline}
The summation is over the $^\lambda H$ described above.

\medskip

In the same way as for the ordinary trace formula (as explained
in \cite{FLN}), our goal is
to obtain the sum decomposition \eqref{want} using the {\em right
  hand side} of formula \eqref{relative}.

Geometrically, we would like to have a decomposition of the right hand
side of \eqref{geom relative} as a direct sum of spaces labeled by the
$^\la H$. We will also need to construct the geometric counterpart of
the extra factor \eqref{extra}.

We now look more closely at what happens in one example.

\subsection{The case of $GL_2$}

Let $\rho$ be the second symmetric power of the defining
representation of $GL_2$. We want to decompose the vector space
\begin{equation}    \label{GL2}
\on{RHom}_{\Bun_{GL_2}}\left(\Psi,{\mathbb K}_{d,\rho}(\Psi) \right),
\end{equation}
into a direct sum of subspaces labeled by the groups $H$, which in
this case are the unramified tori (split and non-split) in
$GL_2$. They are labeled by the unramified double covers of $X$.

Let us look at the vector space \eqref{GL2} more
closely. The sheaf $\Psi$ is supported on the locus of rank two
bundles which have the form \eqref{exten}:
\begin{equation}    \label{exten1}
0 \to K_X^{1/2} \to \V \to K_X^{-1/2} \to 0.
\end{equation}
The support of $\Psi$ consists of bundles $\V$ of this form. We will
fix a point $ \infty \in X$ and identify rank two bundles ${\mc V}$
and ${\mc V}(n [\infty])$ for all $n \in \Z$. In other words, we
replace $\Bun_{GL_2}$ by its quotient by the group $\Z$ generated by
tensoring rank two vector bundles with the line bundle
$\OO([\infty])$.

When we apply our averaging functor $\K_{d,\rho}$, we average over the
Hecke modifications of the bundle $\V$: for each divisor $D = \sum n_i
[x_i]$ of degree $d$ we apply all possible Hecke modifications
corresponding to the representation $\on{Sym}^{n_i} \rho$ at the point
$x_i$, for all $i$, and then take the resulting bundle $\wt\V$ and
replace it by $\wt\V(-d[\infty])$. This bundle should then
again be of the form \eqref{exten1}. Thus, the space \eqref{GL2} can
be interpreted as the cohomology of a certain sheaf over the moduli
space ${\mc W}_{d,\rho}$ of data $(\V,\V',D,\phi)$, where $D = \sum n_i
[x_i]$, $\V$ and $\V'$ are rank two vector bundles on $X$ the form
\eqref{exten1}, together with an embedding
$$
\phi: \V \hookrightarrow \V'(d [\infty]),
$$
such that the quotient
$$
\V'(d[\infty])/\V \simeq \bigoplus_i {\mc T}_{x_i},
$$
where ${\mc T}_{x_i}$ are torsion sheaves supported at
$x_i$:
$$
{\mc T}_{x_i} \simeq {\mc O}_{k_i x_i} \oplus {\mc O}_{(2n_i-k_i)x_i},
\qquad \text{for some } k_i=0,\ldots,n_i
$$
(these correspond to the strata in the affine Grassmannian which
lie in the closure of $\on{Gr}_{n_i\rho}$).

We need to compute the cohomology of this sheaf and relate it to the
cohomologies corresponding to the contributions of $^\la H$ in the sum
on the right hand side of \eqref{want}. We hope that we could use the
methods of \cite{Ngo} and \cite{FGV:w} to do this.

\subsection{The abelian case}    \label{ex}

Finally, we consider the case of the group $GL_1$ and trivial
$\rho$. We wish to check the statement of \secref{geom} that the
cohomology \eqref{rhs} is really the geometric incarnation of the sum
appearing on the left hand side of \eqref{relative}.

Since $\rho$ is trivial, the symmetric power of the curve $X_d$
decouples from the formulas, and we might as well set $\K = \on{Id}$,
and so $N_\sigma = 1$ in formula \eqref{relative}. Note that the
adjoint representation of $GL_1$ is trivial, and so
$$
q^{-d_G} L(\sigma,\on{ad},s) = q^{-(g-1)} \zeta(s),
$$
where $\zeta(s)$ is the zeta-function of $X$. However, since $GL_1$ is
not simple, we have to make an adjustment and replace $\zeta(1)$ by
its residue at $s=1$, which we denote by $\wt\zeta(1)$. In other
words, we remove the factor corresponding to $H^2(X,\Ql)$, which is
independent of $\sigma$. (One can argue that this factor decouples
because the structure sheaf of the full derived stack of rank one local
systems on $X$ contains the exterior algebra of $H^2(X,\Ql)$ placed in
the cohomological degree $1$.)

Thus, we have to compute the sum
\begin{equation}    \label{sum ab}
\sum_\sigma q^{-(g-1)} \wt\zeta(1)^{-1} = q^{-(g-1)} \wt\zeta(1)^{-1}
\cdot \# \{ \sigma: W_F \to GL_1 \},
\end{equation}
where $\sigma$ runs over all unramified homomorphisms $W_F \to GL_1$.

Let $\al_i, i=1,\ldots,2g$, be the eigenvalues of the Frobenius on
$H^1(X,\Ql)$. We have $\al_i = \beta_i q^{1/2}$, where $|\beta_i|=1$
and for each $\beta_i$ there is $\beta_j = \beta_i^{-1}$. Therefore
$$
\wt\zeta(1)^{-1} = \frac{\on{det}(1-q^{-1}
\on{Fr},H^0(X,\Ql))}{\on{det}(1-q^{-1} \on{Fr},H^1(X,\Ql))} = q^{g-1}
\frac{q-1}{\prod_{i=1}^{2g} (1-\al_i)}.
$$
On the other hand, by the abelian class field theory, the number of
unramified $\sigma: W_F \to GL_1$ is equal to the dimension of the
space of unramified automorphic forms on $GL_1(\AD)$, which is the
number of $\Fq$-points of the Jacobian of $X$. This number is in turn
equal to the alternating sum of the traces of the Frobenius of the
cohomology of the Jacobian, which is $\Lambda^\bullet
H^1(X,\Ql)$. Hence the answer is
$$
\prod_{i=1}^{2g} (1-\al_i).
$$
Thus, the sum \eqref{sum ab} is equal to
$q-1$.

Its geometric incarnation is the vector space
\begin{equation}    \label{space}
H^\bullet(\on{Loc}_{GL_1},\OO)
\end{equation}
(since we have replaced $\K_{d,\rho}$ with the identity functor, we
replace the bundle ${\mc T}_{d,\rho}$ with the trivial bundle $\OO$).

It is easy to see that $H^0(\on{Loc}_{GL_1},\OO)$ is one-dimensional,
and $H^i(\on{Loc}_{GL_1},\OO) = 0$ for all $i>0$, so that the total
space \eqref{space} is one-dimensional, situated in cohomological
degree 0.  This means that the (non-existent) Frobenius should act on
it by $q-1$. Let us compare this number with the number appearing on
the right hand side.

First, let us look at the right hand side \eqref{lhs2} of the relative
geometric trace formula \eqref{geom relative}. As above, we replace
$\K_{d,\rho}$ by the identify functor. The resulting vector space is
\begin{equation}    \label{lhs3}
\on{RHom}_{D(\on{Pic}^0(X))}(\Psi,\Psi).
\end{equation}
Here $\Pic^0(X)$ is the Picard moduli stack of line bundles on $X$
of degree $0$ and
$\Psi = p_!(\Ql)$, where $p: \on{pt} \to \on{Pic}^0(X)$ corresponds to
the trivial line bundle. This is also a one-dimensional vector space
situated in cohomological degree zero. So we indeed have an
isomorphism between \eqref{space} and \eqref{lhs3} which is to be
expected because in the case of $G=GL_1$ the categorical Langlands
correspondence is a theorem of \cite{Laumon,Roth}.

Now let us compute the action of the Frobenius on the
one-dimensional vector space \eqref{lhs3}. In the case when $X$ is
defined over $\Fq$, the Frobenius is
well-defined.  Since the trivial line bundle
has the group of automorphisms ${\mathbb G}_m$, the trace of the
Frobenius on \eqref{lhs3} is the same as that on $\phi^!\phi_!(\Ql)$,
where $\phi: \on{pt} \to \on{pt}/{\mathbb G}_m$, which is equal to $\#
{\mathbb G}_m(\Fq) = q-1$. This matches the above calculation for
\eqref{space}.

Finally, let us verify  formula \eqref{trivial case} in the case $G=GL_1$.
Since the tangent bundle to $\Loc_{GL_1}$ is trivial,
with the fiber $H^1(X,\C) = H^1(X,\OO) \oplus H^0(X,K_X)$, and
$H^\bullet(\on{Loc}_{GL_1},\OO) = \C$,
the right hand side of \eqref{trivial case}
is equal to $\Lambda^\bullet(H^1(X,\C))$, which is also isomorphic to
the left hand side of \eqref{trivial case}.

\end{document}